\documentclass[10pt,letterpaper,oneside]{article}
\usepackage{REUReport}

\newcommand{\ope}[1]{\operatorname{#1}}
\newcommand{\boxe}{\ \Box\ }
\newcommand{\suma}{\leq_\oplus}

\title{\bfseries Homotopical Approach to Tensor Products of Modules}
\author{Karthik Boyareddygari\\ Mentor: Alex Sorokin}
\date{}

\begin{document}


\maketitle

\thispagestyle{first}

\tableofcontents

\newpage

\section{Introduction}
Category theory provides a means through which many far-ranging fields of mathematics can be related by their similar structure. In a paper by Robinson \cite{Robinson1981}, this interconnectivity afforded by categorical perspectives allowed for the realization of torsion products as the homotopy groups of a topological space, which is itself constructed for this express purpose. However, even stating this result formally requires a multitude of preliminaries in algebra, topology, and category theory.

The goal of this document is to present a self-contained guide to the fundamental concepts and results, with few proofs, required to do work with this kind of mathematics in hopes of making the field of homotopical algebra more accessible. We only assume familiarity with topological spaces and groups, so it is approachable from an undergraduate level. This project culminates in a discussion of the result of Robinson mentioned above along with a computation as a proof of concept. 

\newpage

\section{Categories}

\subsection{Basic Notions}
\begin{defin}
	A \textit{category} $\Cs{C}$ consists of the listed data:
	\begin{itemize}
		\item a class of \textit{objects} $\mathrm{Obj}(\Cs{C});$
		\item for any two objects $A,B \in \mathrm{Obj}(\Cs{C})$ there is a set $\Mor{\Cs{C}}{A}{B}$ of \textit{morphisms} from $A$ to $B$ such that
		\[
			\Mor{\Cs{C}}{A}{B} \cap \Mor{\Cs{C}}{A'}{B'} = \varnothing,
		\]
		if $A \ne A'$ or $B \ne B';$
		\item for any objects $A,B,C \in \ope{Obj}(\Cs{C})$, there is defined a \textit{composition} operation for morphisms defined by
		\[
			\circ:
			\begin{cases}
				\Mor{\Cs{C}}{B}{C} \times \Mor{\Cs{C}}{A}{B} \to \Mor{\Cs{C}}{A}{C} \\
				(g,f) \mapsto g \circ f
			\end{cases}
		\]
		with the condition that this operation is associative and that for any object $A \in \mathrm{Obj}(\Cs{C})$ there is a morphism $\mathds{1}_A \in \Mor{\Cs{C}}{A}{A}$ called the \textit{identity} such that for any object $B,$ $f \in \Mor{\Cs{C}}{A}{B},$ and $g \in \Mor{\Cs{C}}{B}{A}$:
		\[
			f \circ \mathds{1}_A = f, \qquad \mathds{1}_A \circ g = g.
		\]
\end{itemize}
\end{defin}

\begin{rem}
	In general $\mathrm{Obj}(\Cs{C})$ is a class, not necessarily a set; whenever $\mathrm{Obj}(\Cs{C})$ is a set, $\Cs{C}$ is called a \textit{small category}. Small categories may be visualized as graphs in which vertices correspond to objects and edges correspond to morphisms. For notational simplicity, it is common to write $X \in \Cs{C}$ in place of $X \in\mathrm{Obj}(\Cs{C})$ whenever it does not cause confusion. In addition, it is common to write $\Cs{C}(A,B)$ or even $(A,B)$ in place of $\Mor{\Cs{C}}{A}{B}$ for $A,B \in \Cs{C}$ when the category is understood.
\end{rem}

\begin{defin}
	Let $\Cs{C}$ be a category. Define an \textit{opposite category} $\Cs{C}^{op}$ as having the same objects---but denoting object $A$ by $A^{op}$---and $$\Mor{\Cs{C}^{op}}{A^{op}}{B^{op}} = \Mor{\Cs{C}}{B}{A}$$ with the composition rule $$f^{op} \circ g^{op} = (g \circ f)^{op}.$$
\end{defin}

\begin{theo}[Duality Principle]
	If $P$ is a true statement within a category $\Cs{C}$, then the dual statement $P^{op}$ obtained by reversing all arrows is true about category $\Cs{C}^{op}$.
\end{theo}

\begin{exam}
	The following are examples of some common categories:
	\begin{itemize}
		\item category $\Cb{Set}$ of sets and functions;
		\item category $\Cb{Top}$ of topological spaces and continuous functions;
		\item category $\Cb{Top}_*$ of pointed topological spaces with base point-preserving continuous functions;
		\item category $\Cb{Ab}$ of abelian groups and group homomorphisms;
		\item category $\Cb{Ring}$ of unital rings and ring homomorphisms;
		\item category $\Cb{R\text{-}Mod}$ of left $R$-modules and homomorphisms of left $R$-modules.
	\end{itemize}
\end{exam}

\begin{exam} \label{pos}
	A partially ordered set (poset) $(I,\leq)$ may be viewed as a category $\Cs{P}$ by regarding
	\begin{itemize}
		\item objects $\mathrm{Obj}(\Cs{P}) = \{i\ |\ i \in I\}$,
		\item morphisms $\Mor{\Cs{P}}{i}{j} =
		\begin{cases}
			f_{ij}:i \to j, & \text{if }i \leq j \\
			\varnothing, & \text{otherwise.}
		\end{cases}$
	\end{itemize}
	An especially important poset is $\Delta[n] = \{0 \leq 1 \leq \cdots \leq n\}$, whose category representation may be illustrated diagrammatically as
	\[
		\begin{tikzcd}
			0 \ar[r] & 1 \ar[r] & \cdots \ar[r] & n,
		\end{tikzcd}
	\]
	after omitting the arrows derived from transitivity and identities. Furthermore, $\Delta[n]$ represents the ordinal with $n+1$ elements, which will become important later on in the discussion of simplicial sets.
\end{exam}

\begin{defin}
	A morphism $f$ is called a \textit{monomorphism} (monic) if it is left cancellable; that is for any morphisms $g_1,g_2$, $$fg_1 = fg_2 \Longrightarrow g_1 = g_2.$$ Dually, a morphism $f$ is called an \textit{epimorphism} (epic) if it is right cancellable; that is for any morphisms $h_1,h_2$, $$h_1f = h_2f \Longrightarrow h_1 = h_2.$$ A morphism which is both epic and monic is called a \textit{bimorphism}.
\end{defin}

\begin{rem}
	The definition of a monomorphism is based on the condition for a function between sets to be injective, and in fact monomorphisms do coincide with injective functions in the category of sets. On the other hand, epimorphisms are related to surjections, but do not coincide with them in all categories.
\end{rem}

\begin{prop}
	Epimorphisms are surjections in both categories $\Cb{Set}$ and $\Cb{R\text{-}Mod},$ but they are not necessarily surjective in $\Cb{Ring}.$
\end{prop}

\begin{proof}
	Let us start in $\Cb{Set}.$ Suppose that $f:X \to Y$ is an epimorphism and that $g_1,g_2:Y \to Z$ are any functions such that $g_1f = g_2f.$ This indicates that $g_1$ and $g_2$ coincide on the image of $f,$ but we can say no more because elements of $Y$ outside the image of $f$ may be mapped to anywhere. But our hypothesis implies that $g_1 = g_2$ on all of $Y,$ which is only possible if there are no elements of $Y$ not in the image of $f.$ Hence, $f$ must be surjective.
	
	Now let us investigate the category $\Cb{Ring}.$ It will suffice to give a counterexample, so consider the inclusion $i:\Z \to \Q.$ Again, let $g_1,g_2$ be any ring homomorphisms with source $\Q$ such that $g_1i = g_2i.$ This tells us that $g_1$ coincides with $g_2$ on the integers. Because ring homomorphisms preserve the multiplicative structure, we can deduce where all the reciprocals of the integers are mapped to. Extending this further, we see that every rational is a product of an integer and the reciprocal of another integer, so any homomorphism from $\Q$ is completely determined by where it sends $\Z.$ Therefore, $g_1i = g_2i \Longrightarrow g_1 = g_2,$ but clearly this is not a surjection.
	
	Finally, we examine $\Cb{R\text{-}Mod}.$ Let $f:X \to Y$ be an epimorphism. Now consider the natural projection map $\pi:Y \to Y/\Im(f)$ and the zero map $z:Y \to Y/\Im(f).$ Clearly $\pi f = zf,$ so $\pi = z$ by hypothesis. Projection maps are always surjective, so this implies that $\Im(f) = Y.$ Therefore, $f$ is a surjection.
\end{proof}

\begin{defin} \label{iso}
	A morphism $f$ is called a \textit{retraction} if it is right invertible while a morphism $g$ is called a \textit{section} if it is left invertible; that is, if $$fg = \mathds{1},$$ then $f$ is a retraction, and $g$ is a section. Notice that retractions and sections come in pairs. A morphism $h:A \to B$ which is both a retraction and a section is called an \textit{isomorphism}. In this case, the objects $A$ and $B$ are said to be isomorphic, and we write $A\cong B$.
\end{defin}

\begin{lem}
	Let $\Cs{C}$ be a category. All the following maps are morphisms within $\Cs{C}$ with sources and targets implicit from the allowed compositions.
	\begin{enumerate}[label = \emph{(\roman*)}]
		\item If $s_1$ and $s_2$ are sections, then $s_2s_1$ is a section.
		\item If $r_1$ and $r_2$ are retractions, then $r_1r_2$ is a retraction.
		\item If $s_2s_1$ is a section, then $s_1$ is a section.
		\item If $r_1r_2$ is a retraction, then $r_1$ is a retraction.
	\end{enumerate}
\end{lem}

\begin{proof}
	The proofs for sections will be formally stated but not for retractions, as they are very similar.
	\begin{enumerate}[label = (\roman*)]
		\item Because $s_1,s_2$ are sections, there exist maps $r_1,r_2$ such that $$r_1s_1 = \mathds{1}, ~~~~ r_2s_2 = \mathds{1}.$$ Because we may form the composition $s_2s_1,$ the composition $r_1r_2$ may be formed since the source of $r_1$ is the target of $r_2.$ Then we find that $$(r_1r_2)(s_2s_1) = r_1(r_2s_2)s_1 = r_1\mathds{1}s_1 = r_1s_1 = \mathds{1},$$ so $r_1r_2$ is a left inverse for $s_2s_1.$ Therefore, $s_2s_1$ is a section.
		\item This proof follows from switching the roles of the hypothesis and constructed maps in (i).
		\item Because $s_2s_1$ is a section, there exists a map $r$ such that $rs_2s_1 = \mathds{1}.$ But then $rs_2$ is a left inverse for $s_1,$ so $s_1$ is a section.
		\item Similar to (iii).
	\end{enumerate}
\end{proof}

\begin{defin}
	Consider a category $\Cs{C}$. An object $I \in \Cs{C}$ is called an \textit{initial object} if there is exactly one morphism from $I$ to any object in the category. Dually, an object $T\in\Cs{C}$ is called a \textit{terminal object} if there is exactly one morphism from any object in the category into $T$. An object which is both initial and terminal is called a \textit{zero object}, which is commonly denoted by 0. A category $\Cs{C}$ with a zero object is called a \textit{pointed category}.
\end{defin}

\begin{rem}
	Once cokernels are introduced, it will become evident that the proof from $\Cb{R\text{-}Mod}$ may be recycled for most all pointed categories which have cokernels to show that epimorphisms coincide with surjections. On the other hand, the notions were not equivalent in $\Cb{Ring}$ by exploiting that there exist subrings which generate their field of fractions. As such, it is not difficult to see that epimorphisms are not surjections for commutative monoids.
\end{rem}

\begin{exam}
	Outlined below are some categories with initial, terminal, and zero objects indicated when applicable.
	\begin{itemize}
		\item In $\Cb{Set},$ the empty set is an initial object while any singleton is a terminal object. The same is true in $\Cb{Top}$ when the empty set and singleton are regarded as topological spaces.
		\item In the poset category $(I,\leq),$ the bottom (resp. top) element is an initial (resp. terminal) object, should it exist.
		\item In $\Cb{Ring},$ $\Z$ is the initial object while the zero ring is a terminal object. In fact, the image of $\Z$ in a ring $R$ is called the \textit{prime ring} of $R$ and is used to define the characteristic of a ring.
		\item $\Cb{Ab}$ is a pointed category with zero object the trivial group.
		\item $\Cb{Top}_*$ is a pointed category with zero objects being any singleton.
	\end{itemize}
\end{exam}

\begin{rem}
	Initial, terminal, and zero objects are unique up to unique isomorphism if they exist. Notice that initial objects in $\Cs{C}$ are terminal objects in $\Cs{C}^{op}$ and vice versa. In a pointed category $\Cs{C}$, between any two objects $A,B \in \Cs{C}$ there exists a unique morphism $0_{A,B}:A \to B$, which is defined as the composition $$A \to 0 \to B.$$ Because the zero map is unique, the subscripts in $0_{A,B}$ are usually omitted so that we may instead use 0 to denote both the zero object and the zero map.
\end{rem}

\begin{defin}
	Let $\Cs{C}$ and $\Cs{D}$ be categories. A \textit{covariant functor} $F:\Cs{C}\to\Cs{D}$ is a pair of correspondences between classes of objects and morphisms in $\Cs{C}$ and $\Cs{D}$ defined by
	\[
		F:
		\begin{cases}
			\mathrm{Obj}(\Cs{C}) \to \mathrm{Obj}(\Cs{D}) \\
			A  \mapsto F(A)
		\end{cases} ~~
		F_{A,B}:
		\begin{cases}
			\Mor{\Cs{C}}{A}{B} & \to \Mor{\Cs{D}}{F(A)}{F(B)} \\
			(f: A \to B) & \mapsto (F(f):F(A) \to F(B))
		\end{cases}
	\]
	such that for any composable morphisms $f$ and $g$ and any object $A \in \Cs{C}$: $$F(fg) = F(f)F(g),$$ $$F(\mathds{1}_A) = \mathds{1}_{F(A)}.$$
\end{defin}

\begin{rem}
	Covariant functors preserve the direction of morphisms, and they can be thought of as mappings that respect the category structure (i.e. preserving the identities and compositions of morphisms). Covariant functors are simply referred to as functors.
\end{rem}

\begin{defin}
	A \textit{diagram} in a category $\Cs{C}$ is a functor $D:\Cs{D} \to \Cs{C}$, where $\Cs{D}$ is a small category. Objects are represented by vertices, and morphisms are represented by arrows.
\end{defin}

\begin{exam}
	Recalling that posets may be viewed as categories, consider the poset
	\[
		I = \left\{
		\begin{tikzcd}
			& \bullet \ar[d] \\
			\bullet \ar[r] & \bullet
		\end{tikzcd}
		\right\}
	\]
	Let $\Cs{C}$ be a category and $F:I \to \Cs{C}$ be a functor. Then $F$ is a diagram in $\Cs{C}$ which may be represented as
	\begin{center}
		\begin{tikzcd}
			& A \ar[d,"f"] \\
			B \ar[r,"g"'] & C,
		\end{tikzcd}
	\end{center}
	where $A,B,C \in \Cs{C}$ and $f,g$ are morphisms.
\end{exam}

\begin{exam}
	A very important and commonly encountered type of functor is the covariant representable functor. Let $\Cs{C}$ be a category and $f:A \to B$ be a morphism in $\Cs{C}.$ Define
	\[
		\Mor{\Cs{C}}{X}{-}:
		\begin{cases}
			\Cs{C} \to \Cb{Set} \\
			A \mapsto (X,A) \\
			(A \stackrel{f}{\to} B) \mapsto ((X,A) \stackrel{f_*}{\to} (X,B))
		\end{cases}
	\]
	where the map $f_*$ is defined as
	\[
		\Mor{\Cs{C}}{X}{f} = f_*:
		\begin{cases}
			(X,A) \to (X,B) \\
			\alpha \mapsto f\alpha
		\end{cases}
	\]
	Such defined correspondence is a covariant functor, and it is usually called a covariant representable functor. If it is clear what the source category is, then we write this functor as $(X,-)$. Notice that the definition of $f_*$ as $\alpha \mapsto f\alpha$ is forced by the sources and targets of $f$ and $\alpha$.
\end{exam}

\begin{defin}
	Let $\Cs{C}$ and $\Cs{D}$ be categories. Then a \textit{contravariant functor} $F:\Cs{C}\to\Cs{D}$ is a pair of correspondences between classes of objects and morphisms in $\Cs{C}$ and $\Cs{D}$ defined by
	\[
		F:
		\begin{cases}
			\mathrm{Obj}(\Cs{C}) \to \mathrm{Obj}(\Cs{D}) \\
			A \mapsto F(A)
		\end{cases} ~~
		F_{A,B}:
		\begin{cases}
			\Mor{\Cs{C}}{A}{B} & \to \Mor{\Cs{D}}{F(B)}{F(A)} \\
			(f:A \to B) & \mapsto (F(f):F(B) \to F(A))
		\end{cases}
	\]
	such that for any composable morphisms $f$ and $g$ and any object $A \in \Cs{C}$: $$F(fg) = F(g)F(f),$$ $$F(\mathds{1}_A) = \mathds{1}_{F(A)}.$$
\end{defin}

\begin{rem}
	A contravariant functor $F:\Cs{C} \to \Cs{D}$ can be considered as a covariant functor $F':\Cs{C}^{op} \to \Cs{D}$. Therefore, contravariant functors are covariant functors acting either to or from an opposite category. Recalling the remark about covariant functors preserving category structure, the involvement of an opposite category produces a sort of reversing of arrows.

	Within a diagram of a category, this translates to all arrows being turned backwards while the objects are left stationary. Again, the term functor usually refers to covariant functor unless it is specified that it is a contravariant functor.
\end{rem}

\begin{defin}
	Given functors $F:\Cs{C} \to \Cs{D}$ and $G:\Cs{D} \to \Cs{E}$, we define the \textit{composition of functors} $G \circ F$ as $$(G \circ F)(X) = G(F(X)),$$ $$(G \circ F)(f) = G(F(f)),$$ where $f$ is a morphism in $\Cs{C}$ and $X \in \Cs{C}$. Note that these functors may be of opposite variance.
\end{defin}

\begin{exam}
	Let's consider the contravariant representable functor. Let $\Cs{C}$ be a category and $X$ be an object in $\Cs{C}.$ Define
	\[
		\Mor{\Cs{C}}{-}{X}:
		\begin{cases}
			\Cs{C} \to \Cb{Set} \\
			A \mapsto (A,X) \\
			(A \stackrel{f}{\to} B) \mapsto ((B,X) \stackrel{f^*}{\to} (A,X))
		\end{cases}
	\]
	where $f^*$ is defined as,
	\[
		\Mor{\Cs{C}}{f}{X} = f^*:
		\begin{cases}
			\Mor{\Cs{C}}{B}{X} \to \Mor{\Cs{C}}{A}{X} \\
			\beta \mapsto \beta f
		\end{cases}
	\]
	It is easy to see that it is a contravariant functor. Again if it is clear what the source category is, then we write this functor as $(-,X)$. Although less straightforward, the action of this functor on morphisms is still the only way that makes sense for correct composition.
\end{exam}

\begin{exam}
	Here we observe two simple but important functors with no distinction between contravariant and covariant. Let $\Cs{C}$ and $\Cs{D}$ be categories, $A \in \Cs{C}$ be any object in $\Cs{C}$, $f$ be a morphism in $\Cs{C}$, and $D \in \Cs{D}$ be some object in $\Cs{D}$. The identity functor $\mathds{1}_\Cs{C}$ on $\Cs{C}$ defined as $$\mathds{1}_\Cs{C}(A) = A, ~~ \mathds{1}_\Cs{C}(f) = f,$$ i.e. it is the functor which fixes all objects and morphisms.
	The constant functor with the value $D$ defined as
	\[
		D:
		\begin{cases}
			\Cs{C} & \to \Cs{D} \\
			A & \mapsto D \\
			f & \mapsto \mathds{1}_D
		\end{cases}
	\]
	is the functor which sends all objects to $D$ and all morphisms to the identity morphism of $D$.
\end{exam}


\subsection{Natural Transformations and Yoneda Lemma}
\begin{defin}
	Given categories $\Cs{C}$ and $\Cs{D}$ and covariant functors $F,G:\Cs{C} \to \Cs{D}$, a \textit{natural transformation} $\eta:F \to G$ is a one-parameter family of morphisms in $\Cs{D}$ $$\eta = \{\eta_X:F(X) \to G(X)\}_{X \in \Cs{C}},$$ such that the following diagram commutes for every morphism $f \in \Mor{\Cs{C}}{X}{Y}$,
	\begin{center}
		\begin{tikzcd}
			F(X) \ar[d,"F(f)"'] \ar[r,"\eta_X"] & G(X) \ar[d,"G(f)"] \\
			F(Y) \ar[r,"\eta_Y"'] & G(Y).
		\end{tikzcd}
	\end{center}
\end{defin}

A natural transformation is a means of relating functors. Natural transformations also exist for contravariant functors, and the definition is the same except the commuting diagram has its vertical arrows pointing up. For now, we describe when a natural transformation provides an equivalence.

\begin{defin} \label{nequ}
	Let $F,G:\Cs{C} \to \Cs{D}$ be functors. A natural transformation $\eta:F \to G$ is called a \textit{natural isomorphism} if $\eta_X$ is an isomorphism in $\Cs{D}$ for each $X \in \Cs{C}$. In this case the functors $F$ and $G$ are said to be isomorphic, and we write $F \approx G$.
\end{defin}

\begin{defin}
	Given natural transformations $\eta:F \to G$ and $\tau:G \to H$ of the functors $F,G,H:\Cs{C} \to \Cs{D}$, we may define the \textit{composition of natural transformations} $$\tau\eta = \{\tau_X \eta_X:F(X) \to H(X)\}_{X \in \Cs{C}}$$ as the composition of the component morphisms $\tau_X$ and $\eta_X$ for any object $X$.
\end{defin}

\begin{defin}
	Given categories $\Cs{C}$ and $\Cs{D}$, a functor $F:\Cs{C}\to\Cs{D}$ is said to be an \textit{equivalence} if there exists a functor $G:\Cs{D}\to\Cs{C}$ such that $$FG \approx \mathds{1}_\Cs{D}, ~~ GF \approx \mathds{1}_\Cs{C}.$$ In this case the categories $\Cs{C}$ and $\Cs{D}$ are said to be equivalent, and we write $\Cs{C} \approx \Cs{D}$. It is clear that $G$ would then also be an equivalence of the categories $\Cs{C}$ and $\Cs{D}$.
\end{defin}

Given a functor $F$, it may not be easy to demonstrate a compatible functor $G$ to establish that $F$ is an equivalence. The following theorem provides a means for checking whether a functor is an equivalence based on its own properties.

\begin{theo}
	A functor $F:\Cs{C} \to \Cs{D}$ is an equivalence of categories if and only if it is:
	\begin{itemize}
		\item faithful (injective on morphisms): $F_{X,Y}:\Mor{\Cs{C}}{X}{Y} \to \Mor{\Cs{D}}{F(X)}{F(Y)}$ is an injective function;
		\item full (surjective on morphisms): $F_{X,Y}:\Mor{\Cs{C}}{X}{Y} \to \Mor{\Cs{D}}{F(X)}{F(Y)}$ is a surjective function;
		\item essentially surjective (surjective on objects): there exists $C \in \Cs{C}$ such that $F(C) = D$ for each $D \in \Cs{D}$.
	\end{itemize}
\end{theo}

For the next lemma, we must introduce a piece of notation. For functors $F,G:\Cs{C}\to\Cs{D}$, let $(F,G)$ be the collection of all natural transformations between the functors $F$ and $G$.

\begin{lem}[Yoneda Lemma]
	Let $\Cs{C}$ be a category and $C \in \Cs{C}$ be any object. If $F:\Cs{C} \to \Cb{Set}$ is a functor, then there is a bijection 
	\[
		y:
		\begin{cases}
			((C,-),F) \to F(C) \\
			\eta \mapsto \eta_C(\mathds{1}_C).
		\end{cases}
	\]
\end{lem}

\begin{lem}[Contravariant Yoneda Lemma]
	Let $\Cs{C}$ be a category and $C \in \Cs{C}$ be any object in $C$. If $F:\Cs{C} \to \Cb{Set}$ is a contravariant functor, then there is a bijection 
	\[
		y:
		\begin{cases}
			((-,C),F) \to F(C) \\
			\eta \mapsto \eta_C(\mathds{1}_C).
		\end{cases}
	\]
\end{lem}

With these two versions of the Yoneda lemma, we obtain the following corollary.

\begin{cor} \label{bij}
	For any $A,B\in\Cs{C}$, there are bijections $$((A,-),(B,-)) \leftrightarrow (B,A),$$ $$((-,A),(-,B)) \leftrightarrow (A,B).$$
\end{cor}

To obtain the above bijections, we apply the Yoneda Lemma with the functor $F$ being $(B,-)$ and the contravariant Yoneda Lemma with the functor $F$ being $(-,B).$

\begin{defin}
	Let $\Cs{C}$ be a small category and $\Cs{D}$ be any category. Then the \textit{functor category} $\Cs{D}^\Cs{C},$ also denoted $(\Cs{C},\Cs{D}),$ is the category whose objects are functors and morphisms are natural transformations, and the composition is defined as a composition of natural transformations:
	\begin{itemize}
		\item $\mathrm{Obj}(\Cs{D}^\Cs{C}) = \{F:\Cs{C} \to \Cs{D}\ |\ F \text{ is a functor}\};$
		\item $\Mor{\Cs{D}^\Cs{C}}{F}{G} = (F,G)$.
	\end{itemize}
\end{defin}

Now with the necessary terminology in place, we may proceed to the following very important corollary. In essence, it asserts that any small category is a subcategory of $\Cb{Set}$ for which all morphism sets are left completely intact. Said differently, small categories are full subcategories of $\Cb{Set}$.

\begin{cor}[Yoneda Embedding]
	If $\Cs{C}$ is a small category, then the functor
	\[
		Y:
		\begin{cases}
			\Cs{C}^{op} \to \Cb{Set}^\Cs{C} \\
			C \mapsto (C,-)
		\end{cases}
	\]
	is injective on objects, full, and faithful. This functor is called the Yoneda Embedding.
\end{cor}

\begin{proof}
	To show that $Y$ is injective on objects, consider the fact that for objects $A,B \in \Cs{C}$ such that $A \neq B$, $(A,-)$ and $(B,-)$ will have disjoint images by the definition of a category. As such, $(A,-) \neq (B,-)$ so that $Y$ is injective on objects. Now recall from \thref{bij} that there is a bijection between $\Mor{\Cs{C}}{A}{B} = \Mor{\Cs{C}^{op}}{B}{A}$ and $((B,-),(A,-))$. Because natural transformations are morphisms in $\Cb{Set}^\Cs{C}$, this shows that $Y$ is full and faithful.
\end{proof}

\begin{cor}
	Let $A,B \in \Cs{C}$ be objects of $\Cs{C}$. If for every $C \in \Cs{C}$ there is a bijection $$(A,C) \leftrightarrow (B,C),$$ then $A \cong B$ in $\Cs{C}$. The converse also holds.
\end{cor}

Because bijections are isomorphisms in $\Cb{Set}$, this corollary amounts to saying that if $Y(A) \approx Y(B)$, then $A \cong B$. Seeing as $\approx$ indicates isomorphism in $\Cb{Set}^\Cs{C}$ and that $Y$ is injective on objects, this corollary closely resembles the condition for injectivity of functions.

The reader may have noted the many contexts where the notation $(-,-)$ appears. For example replacing the dashes with objects gives a collection of morphisms, categories a collection of functors, functors (with identical source and target) a collection of natural transformations, and functors (with identical target) a collection of comma categories. While it may seem the notation is overused, the usage is quite similar across the cases because they all take structures and return a collection of things which relate those structures. Care need only be taken in paying attention to what the input structures are.


\subsection{Limits and Colimits}
\begin{defin}
	Let $\Cs{C},\Cs{D},\Cs{E}$ be categories and $F:\Cs{C} \to \Cs{E},G:\Cs{D} \to \Cs{E}$ be functors. Then a \textit{comma category} $F/G $ is a category whose:
	\begin{itemize}
		\item objects are triples $(C,D,f)$ where $C \in \Cs{C}$, $D \in \Cs{D}$, and $f:F(C) \to G(D)$ is a morphism in $\Cs{E}$;
		\item morphisms between $(C_1,D_1,f_1)$ and $(C_2,D_2,f_2)$ are pairs $(\alpha,\beta)$, where 
		\[
			\alpha:F(C_1) \to F(C_2), \quad \beta:G(D_1) \to G(D_2)
		\]
		are morphisms in $\Cs{E}$ such that $f_2\alpha = \beta f_1$. It can be expressed diagrammatically as
		\begin{center}
			\begin{tikzcd}
				F(C_1) \ar[d,"f_1"'] \ar[r,"\alpha"] & F(C_2) \ar[d,"f_2"] \\
				G(D_1) \ar[r,"\beta"'] & G(D_2).
			\end{tikzcd}
		\end{center}
	\end{itemize}
\end{defin}

\begin{defin}
	In a category $\Cs{C}$, consider a family of objects and morphisms between them $$(\{C_i\},\{f_{ij}:C_i \to C_j\})_{i,j \in I}.$$ A \textit{cone} over this family is a pair $(C,\{f_i:C \to C_i\}_{i \in I})$ such that the following diagram commutes for any $i,j \in I$:
	\begin{center}
		\begin{tikzcd}
			& \ar[ld,"f_i"'] C \ar[rd,"f_j"] & \\
			C_i \ar[rr,"f_{ij}"'] && C_j.
		\end{tikzcd}
	\end{center}
	Similarly, a \textit{cocone} below this family is a pair $(C,\{f_i:C_i \to C\}_{i \in I})$ such that the following diagram commutes for any $i,j \in I$:
	\begin{center}
		\begin{tikzcd}
			C_i \ar[rd,"f_i"'] \ar[rr,"f_{ij}"] && \ar[ld,"f_j"] C_j \\
			& C. &
		\end{tikzcd}
	\end{center}
	The maps $f_i$ associated with the cone and cocone are called \textit{structural maps}.
\end{defin}

We can define a preorder between cones that will be useful in stating the primary structure of this section.

\begin{defin}
	A cone $(C,\{f_i\}_{i \in I})$ is \textit{dominated} by a cone $(C',\{f'_i\}_{i \in I})$ if there exists a \textit{domination map} $g:C' \to C$ such that $$f'_i = f_ig$$ for all $i \in I$. Dually, a cocone $(C,\{f_i\}_{i \in I})$ is \textit{dominated} by a cocone $(C',\{f_i'\}_{i \in I})$ if there exists a \textit{domination map} $g':C \to C'$ such that $$f'_i = g'f_i$$ for all $i \in I$.
\end{defin}

Below is a remark with diagrams to help demonstrate cones and domination visually. The diagrams will illustrate domination for cones, and the corresponding diagrams for cocones may be obtained by reversing all structural maps.

\begin{rem}
	Let the family of objects and morphisms make up the base of a geometric cone. Then the cone (resp. cocone) is the object at the tip of the geometric cone, resulting in a diagram with the following form:
	
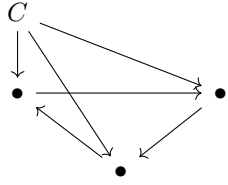
\begin{center}
		\begin{tikzcd}
			\ar[d] C \ar[rdd] \ar[rrd] & & \\
			\bullet \ar[rr] & & \ar[ld] \bullet \\
			& \ar[lu] \bullet &
		\end{tikzcd}
	\end{center}	
	
	When a cone $C'$ dominates a cone $C$, a new tip is introduced with an arrow $g$ pointing from $C'$ to $C$ (resp. from $C$ to $C'$), and $C',$ and the new tip becomes the highest point. Now we form a new cone (resp. cocone) from the new base and tip by letting the descending arrows from tip to base (resp. base to tip) be the hypotenuses of the right triangle with legs $g$ and $f_i$, resulting in a diagram with the following form:
	\begin{center}
		\begin{tikzcd}
			& & \ar[lldd,dashed] \ar[lld,"g"'] \ar[lddd,dashed] C' \ar[dd,dashed] \\
			\ar[d] C \ar[rdd] \ar[rrd] & & \\
			\bullet \ar[rr] & & \ar[ld] \bullet \\
			& \ar[lu] \bullet &
		\end{tikzcd}
	\end{center}
\end{rem}

\begin{defin}
	A \textit{limit} of a functor $F \in \Cs{C}^I$, where $(I,\leq)$ is a poset category (refer to \ref{pos}), is a cone $(C,\{f_i\}_{i \in I})$ over the family $(\{F(i)\},\{F(i) \to F(j)\})_{i \leq j \in I}$ which is dominated by every other cone, such that the domination maps are unique. We write $C = \varprojlim F$. So if $(C',\{f'_i\}_{i \in I})$ is a cone over the family, then 
	
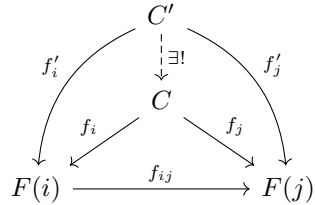
\begin{center}
		\begin{tikzcd}
			& \ar[ldd,bend right,"f'_i"'] \ar[d,dashed,"\exists!"] C' \ar[rdd,bend left,"f'_j"] & \\
			& \ar[ld,"f_i"'] C \ar[rd,"f_j"] & \\
			F(i) \ar[rr,"f_{ij}"] & & F(j)
		\end{tikzcd}
	\end{center}
	commutes. If the domination maps are not unique, it is called a \textit{weak limit}. Limits are also called projective limits or inverse limits.

	Similarly, a \textit{colimit} of a functor $F \in \Cs{C}^I$, where $(I,\leq)$ is a poset, is a cocone $(C,\{f_i\}_{i \in I})$ below the family $(\{F(i)\},\{F(i) \to F(j)\})_{i,j \in I}$ which is dominated by every other cocone, such that the domination maps are unique. We write $C = \varinjlim F$. So if $(C',\{f'_i\}_{i \in I})$ is a cocone under the family, then 
	
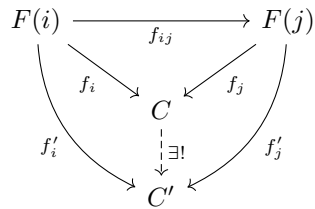
\begin{center}
		\begin{tikzcd}
			F(i) \ar[rd,"f_i"'] \ar[rdd,bend right,"f'_i"'] \ar[rr,"f_{ij}"'] & & \ar[ld,"f_j"] \ar[ldd,bend left,"f'_j"] F(j) \\
			& \ar[d,dashed,"\exists!"] C & \\
			& C' &
		\end{tikzcd}
	\end{center}
	commutes. If the domination maps are not unique, it is called a \textit{weak colimit}. Colimits are also called injective limits or direct limits.
\end{defin}

\begin{defin}
	A category $\Cs{C}$ is said to be (\textit{co})\textit{complete} if it admits all small (co)limits (i.e. there exists a (co)limit for any diagram $F$ in $\Cs{C}$). Similarly, $\Cs{C}$ is said to be \textit{finitely (co)complete} if it admits all finite (co)limits. That is, for any finite poset $\Cs{D}$ and functor $F:\Cs{D} \to \Cs{C}$, there exists a (co)limit of $F$.
\end{defin}


\subsection{Adjunction of Functors}
Here we observe how a relation between functors resulting from simple changes of perspective produce a litany of useful and important information about the functors. One such motivating example originates from a function of 2-variables
\[
	f:
	\left\{
	\begin{array}{rcl}
		X \times Y & \longrightarrow & Z \\
		(x,y) & \longmapsto & z,
	\end{array}
	\right.
\]
which we may write alternatively as $f(x,y) = z.$ It is common in multivariable calculus courses to parameterize such multivariable functions, each time increasing the number of functions but simultaneously decreasing the number of variables of these functions. In particular, we may instead write $f_x(y) = f(x,y)$ to define a family of maps, one for each value $x.$ Formally, this function may be written
\[
	f_x:
	\left\{
	\begin{array}{rcl}
		Y & \longrightarrow & Z \\
		y & \longmapsto & z
	\end{array}
	\right.
	.
\]
What we find is that every 2-variable function corresponds to a family of 1-variable functions, expressed symbolically as
\[
	(X \times Y,Z) \longleftrightarrow (X,(Y,Z)).
\]
By defining mappings $L = - \times Y$ and $R = (Y,-),$ we may rewrite the above bijection as
\[
	(L(X),Z) \longleftrightarrow (X,R(Z)),
\]
which no longer explicitly mentions $Y.$ 

There is yet another motivating example which instead comes from the scalar product of vectors. Let $u,v \in \mathbb{R}^n$ be any vectors and $O \in (\mathbb{R}^n,\mathbb{R}^n)$ be some linear operator (equivalently an $n \times n$ matrix). The scalar product of $x$ and $y$ is given by $x\cdot y = x^Ty,$ which is nothing more than the familiar dot product from linear algebra. Introducing the linear transformation, we notice that
\[
	x \cdot Oy = x^TOy = x^T(O^T)^Ty = (O^Tx)^Ty = O^Tx \cdot y.
\]
Intuitively, the scalar product measures the ``$x$ character of $y,$'' which most often is geometrically interpreted as how parallel $y$ is to $x.$ By operating on $y$ before taking the scalar product, we actually measure how parallel the image of $y$ (under $O$) is to $x.$ The important change of perspective demonstrated by the above chain of equalities is that this measurement is equivalent to first transforming $x$ in a ``reverse'' fashion and then measuring it against $y,$ where reverse is used very loosely to designate how transpose matrices act geometrically. If we change the notation to, $(x,y) = x \cdot y,$ then the equality may be written suggestively as
\[
	(x,Oy) = (O^Tx,y)
\]
so the operator may now go without being explicitly mentioned. In the special case where the operator is orthogonal (geometrically representing a rotation), operating on $y$ and then measuring gives the same result as operating on $x$ in reverse and then measuring. Alternatively, it is like thinking of a rotation not as a transformation of vectors but instead as a reverse transformation of the space within which the vectors sit while leaving the vectors to fall in their new spots.

\begin{defin}
	Let $\Cs{C},\Cs{D}$ be categories and $L:\Cs{C} \to \Cs{D},R:\Cs{D} \to \Cs{C}$ be functors. If for any $C \in \Cs{C}, D \in \Cs{D}$ there is a bijection
	\[
		\varphi_{C,D}:\Mor{\Cs{D}}{L(C)}{D} \to \Mor{\Cs{C}}{C}{R(D)}
	\]
	which is natural in both $C$ and $D,$ then we say that the functors $L$ and $R$ are \textit{adjoint}, $(L,R)$ is an \textit{adjoint pair}, $L$ is \textit{left adjoint} to $R,$ or $R$ is \textit{right adjoint} to $L;$ we write this symbolically as $L \dashv R$ or as
	\[
		\begin{tikzcd}
			\Cs{C} \ar[r, "L"{name = L}, bend left = 25] &\Cs{D}. \ar[l, "R"{name = R}, bend left = 25] \ar[phantom, from = L, to = R, "\dashv" rotate = -90]
		\end{tikzcd}
	\]
\end{defin}

\begin{rem}
	What is meant by saying that $\varphi_{C,D}$ is natural in $C$ and $D$ is that for all $C_1,C_2 \in \Cs{C},$ $D_1,D_2 \in \Cs{D},$ $f \in \Mor{\Cs{C}}{C_1}{C_2},$ and $g \in \Mor{\Cs{D}}{D_1}{D_2},$ the diagram
	\[
		\begin{tikzcd}
			\ar[d,"\Mor{\Cs{D}}{Lf}{g} = g \circ \Box \circ Lf"'] \Mor{\Cs{D}}{L(C_1)}{D_1} \ar[r,"\varphi_{C_1,D_1}"] & \ar[d, "\Mor{\Cs{C}}{f}{Rg} = Rg \circ \Box \circ f"] \Mor{\Cs{C}}{C_1}{R(D_1)} \\
			\Mor{\Cs{D}}{L(C_2)}{D_2} \ar[r,"\varphi_{C_2,D_2}"'] & \Mor{\Cs{C}}{C_2}{R(D_2)}
		\end{tikzcd}
	\]
	commutes. The word ``natural'' does in fact relate to natural transformations, and it is useful to stop and think about explicitly how these tie together.
\end{rem}

Intuitively, given a pair of objects in separate categories, a pair of adjoint functors provides rules for finding another pair of category-separated objects such that there is a one-to-one correspondence between the morphisms for objects in the same category. As we will see more in later discussions, bijections between sets of morphisms often indicate some type of universal arrow mechanism; this holds true here as well, which is made precise in the following theorems.

\begin{theo}
	Let $L:\Cs{C} \to \Cs{D}$ be a functor. Then $L$ has a right adjoint if and only if there exists a functor $R:\Cs{D} \to \Cs{C}$ such that $L/D$ has terminal object $(R(D),D,L(R(D)) \stackrel{\varepsilon_D}{\longrightarrow} D)$ for any $D \in \Cs{D}.$ Under such conditions, $L \dashv R.$
\end{theo}

\begin{theo}
	Let $R:\Cs{D} \to \Cs{C}$ be a functor. Then $R$ has a left adjoint if and only if there exists a functor $L:\Cs{C} \to \Cs{D}$ such that $C/R$ has initial object $(C,L(C),C \stackrel{\mu_C}{\longrightarrow} R(L(C)))$ for any $C \in \Cs{C}.$ Under such conditions, $L \dashv R.$
\end{theo}

The universal arrows in these theorems arise from the unique arrows originating from initial objects (resp. ending at terminal objects). These universal arrows, denoted by $\varepsilon_D$ and $\mu_C$ in the above theorems, are in fact component maps of natural transformations which form the foundation of the next equivalent formulation of adjunction.

\begin{defin}
	Let $F:\Cs{C} \to \Cs{D}$ and $G:\Cs{D} \to \Cs{C}$ be a pair of functors. The natural transformations (if they exist)
	\[
		\mu:\mathds{1}_{\Cs{C}} \to GF, \quad \varepsilon:FG \to \mathds{1}_{\Cs{D}}
	\]
	are called \textit{unit} and \textit{counit}, respectively.
\end{defin}

\begin{theo}
	Let $F:\Cs{C} \to \Cs{D}$ and $G:\Cs{D} \to \Cs{C}$ be a pair of functors. There exists a unit-counit pair
	\[
		\mu:\mathds{1}_{\Cs{C}} \to GF, \quad \varepsilon:FG \to \mathds{1}_{\Cs{D}}
	\]
	satisfying the ``counit-unit equations''
	\begin{itemize}
		\item $\mathds{1}_L = \varepsilon L \circ L\mu,$
		\item $\mathds{1}_R = R\varepsilon \circ \mu R,$
	\end{itemize}
	if and only if $L \dashv R.$
\end{theo}

\begin{cor}
	Adjoints are defined up to natural isomorphism, and the composition of adjunctions is an adjunction.
\end{cor}

\begin{exam} \label{adj}
	While the reader may not be familiar with all of the examples to be stated here, there will be brief descriptions of functors which are not discussed further in this document. Some adjunctions are discussed at length elsewhere and so will be omitted from this list. Regardless, some of these examples will help to demonstrate just how ubiquitous adjunctions are.
	\begin{itemize}
		\item $(-)_1:\Cb{Rng} \dashv \Cb{Ring}:F,$ where $F$ forgets about the role of the identity and $(-)_1$ adjoins the identity. If $R$ is a rng, then $(R)_1 = R \times \Z$ with the operation $(r,n)(s,m) = (rs+rm+ns,nm)$ and identity element $(0,1);$
		\item $K_0:\Cb{CMon} \dashv \Cb{Ab}:F,$ where $K_0$ adjoins inverses through a particular quotient construction (study Grothendieck group for reference) and $F$ forgets about the role of inverses;
		\item Group ring: $\Cb{Grp} \dashv \Cb{Ring}:$ Group of units, where a group is sent to its group ring over $\Z$ and a ring is sent to its group of units;
		\item $Frac(-):$ Domains $\dashv$ Fields: $F,$ where $Frac(R)$ takes the domain $R$ to its field of fractions and $F$ forgets the role of multiplicative inverses;
		\item $\Sigma:\Cb{Top}_* \dashv \Cb{Top}_*:\Omega,$ where $\Sigma$ and $\Omega$ produce the reduced suspension and loop space of a pointed topological space, respectively;
		\item $\beta:\Cb{Top}_2 \dashv \Cb{KHaus}:\iota,$ where $\beta$ is the Stone-\v{C}ech compactification, i.e. smallest compact space in which we may embed the given Hausdorff space, and $\iota$ is the inclusion of compact Hausdorff spaces into Hausdorff spaces;
		\item If $F:\Cs{C} \stackrel{\approx}{\longrightarrow} \Cs{D}:G$ is an equivalence of categories, then $F \dashv G,$ $F \vdash G,$ and $\varepsilon,\mu$ are natural isomorphisms.
	\end{itemize}
\end{exam}

Lastly, there is a very desirable property of adjoint pairs.

\begin{defin}
	A functor $F:\Cs{C} \to \Cs{D}$ is called \textit{continuous} if
	\[
		F(\varprojlim_iC_i) \cong \varprojlim_iF(C_i)
	\]
	and cocontinuous if
	\[
		F(\varinjlim_iC_i) \cong \varinjlim_iF(C_i).
	\]
\end{defin}

\begin{theo}
	If $L \dashv R$ is an adjoint pair of functors, then $L$ is cocontinuous and $R$ is continuous.
\end{theo}


\newpage

\section{Rings and Modules}

\subsection{Fundamentals of Rings and Modules}
Although rings are common knowledge to many mathematicians, we start by setting straight what conventions will be followed throughout this document.
\begin{defin}
	A nonempty set $R$ with two binary operations
	\[
		+:
		\left\{
		\begin{array}{rcl}
			R \times R & \longrightarrow & R \\
			(a,b) & \longmapsto & a+b
		\end{array}
		\right.
		, \qquad \qquad
		\cdot:
		\left\{
		\begin{array}{rcl}
			R \times R & \longrightarrow & R \\
			(a,b) & \longmapsto & ab
		\end{array}
		\right.
	\]
	subject to the conditions
	\begin{enumerate}[label = (\roman*)]
		\item $(R,+)$ is an abelian group (with identity element 0);
		\item there exists an element $1 \in R$ such that $1a = a1 = a$ for any $a \in R$ (existence of multiplicative identity);
		\item $(ab)c = a(bc)$ for any $a,b,c \in R$ (associativity of multiplication);
		\item $a(b+c) = ab+ac$ for any $a,b,c \in R$ (left distributivity);
		\item $(a+b)c = ac+bc$ for any $a,b,c \in R$ (right distributivity);
	\end{enumerate}
	is called a unital \textit{ring}. The operations are referred to as \textit{addition} ($+$) and multiplication ($\cdot$). Imposing only conditions (ii) and (iii) on $\cdot$ makes $(R,\cdot)$ a \textit{monoid}, i.e. a group without the condition of having inverses for every element.
\end{defin}

\begin{defin}
	By adding (or removing) further axioms, we may define other related structures. Let $R$ be a ring.
	\begin{itemize}
		\item If there is not a multiplicative identity, then $R$ is called a \textit{rng}.
		\item If $ab = ba$ for any $a,b \in R,$ then $R$ is called a \textit{commutative ring}.
		\item If there are non-zero elements $a,b \in R$ such that $ab = 0,$ then we call $a$ a \textit{left zero-divisor} and $b$ a \textit{right zero-divisor}. If $R$ has no left and right zero-divisors, we call $R$ a \textit{domain}.
		\item If $R$ is a commutative domain, it is called an \textit{integral domain}.
		\item If for every $a \in R-\{0\}$ there exists an element $a^{-1} \in R$ such that $aa^{-1} = a^{-1}a = 1,$ then $R$ is called a \textit{division ring}. Clearly, $R$ must be a unital ring first in order for 1 to exist. We call $a^{-1}$ the \textit{multiplicative inverse} of $a.$ Informally, a division ring is a ring for which every non-zero element has a multiplicative inverse.
		\item If $R$ is a commutative division ring, then it is a \textit{field}.
	\end{itemize}
\end{defin}

The multiplicative identity is usually denoted by 1, and it is often called the \textit{unity} for a ring. As such, rngs are then just rings without unity. For those uncomfortable with the idea of zero-divisors, consider that non-zero matrices can multiply to the zero matrix and that a zero-divisor of this kind is dependent on which side multiplication occurs. Also, the reader may have noticed that there is no mention of zero-divisors with respect to fields, which leads us to this simple and illustrative result.

\begin{lem}
	If $R$ is a field, then $R$ is an integral domain.
\end{lem}

\begin{proof}
	Suppose that $R$ is a field and that there are non-zero elements $a,b \in R$ such that $ab = 0.$ Then we have that
	\[
		ab = 0 \Longrightarrow a^{-1}ab = a*0 \Longrightarrow 1*b = 0 \Longrightarrow b = 0,
	\]
	but this contradicts our assumption that $b$ is non-zero.
\end{proof}

There are many more properties which may be obtained from the axioms for any of these structures, but these may be found easily from many other sources. We will now provide some other basic but important definitions.

\begin{defin}
	Let $R$ be a ring. Then a subset $I \subseteq R$ is a \textit{left ideal} if it satisfies the following axioms:
	\begin{enumerate}[label = (\roman*)]
		\item if $i_1,i_2 \in I,$ then $i_1-i_2 \in I.$
		\item for any $r \in R$ and $i \in I,$ $ri \in I$ ($RI \subseteq I$).
	\end{enumerate}
\end{defin}

\begin{rem}
	Axiom (i) establishes that $I$ is a subgroup of $(R,+)$ while axiom (ii) establishes that $I$ is closed under multiplication by ring elements.
\end{rem}

There are also right ideals and two-sided ideals. In commutative rings, all ideals are clearly two-sided, so mention of sidedness disappears. We now mention how to generate ideals because it will later mimic how to generate modules. Although this definition is stated for commutative rings, the definitions are similar for non-commutative rings but with the inclusion of handedness.

\begin{defin}
	Let $S \subseteq R$ be a subset of a commutative ring, then the \textit{ideal generated by $S$}, denoted by $\langle S\rangle$, is defined as
	\[
		R\langle S\rangle = \langle S\rangle = \left\{\left.\sum_{i = 1}^n r_i s_i\right|r_i \in R,s_i \in S, n \geq 1\right\}.
	\]
	If $S$ is finite, then $\langle S\rangle$ is said to be a \textit{finitely generated ideal}. If $S = \{s\}$, then $\langle S\rangle = sR$ is a \textit{principal ideal} generated by $s$.
\end{defin}

\begin{rem}
	With groups, there were distinguished subgroups (normal subgroups) for which the set of cosets of the subgroup took on a group structure. With all subrings (subsets which are rings sharing the same identity as the ambient ring), there is a group structure on these cosets, but not necessarily a ring structure. Ideals are of particular importance because the cosets of ideals have a well-defined ring structure.
\end{rem}

\begin{defin}
	Let $R$ and $S$ be rings. A function $\varphi:R \to S$ satisfying
	\begin{itemize}
		\item $\varphi(r_1+r_2) = \varphi(r_1)+\varphi(r_2),$
		\item $\varphi(r_1r_2) = \varphi(r_1)\varphi(r_2),$
		\item $\varphi(1) = 1,$
	\end{itemize}
	is called a \textit{ring homomorphism}. A bijective ring homomorphism is called a \textit{ring isomorphism}.
\end{defin}

It may seem odd that it is included as a separate axiom that the multiplicative identity is preserved, but it is not implied by the other axioms due to the lack of multiplicative inverses. In addition, this is a very important condition to impose as it entirely changes the category of study if it is not present, which will be seen in later sections. Now we observe some distinguished subsets constructed from ring homomorphisms.

\begin{defin}
	Let $\varphi:R \to S$ be a ring homomorphism. The \textit{kernel} of $\varphi,$ denoted by $\Ker(\varphi)$ is defined as
	\[
		\Ker(\varphi) = \{r \in R\ |\ \varphi(r) = 0 \in S\}.
	\]
	Furthermore, the \textit{image} of $\varphi,$ denoted by $\Im(\varphi)$ is defined as
	\[
		\Im(\varphi) = \{\varphi(r) \in S\ |\ r \in R\}.
	\]
\end{defin}

\begin{defin}
	Let $S$ be a ring. A subset $T \subseteq S$ is called a \textit{subring} if there exists a ring $R$ and ring homomorphism $\varphi:R \to S$ such that $T = \Im(\varphi).$
\end{defin}

\begin{rem}
	Kernels of ring homomorphisms are two-sided ideals. Furthermore, images of ideals are not necessarily ideals unless the ring homomorphism is surjective.
\end{rem}

\begin{defin}
	Let $I \subseteq R$ be a left ideal of the ring $R.$ Then the \textit{quotient ring} $R/I$ is defined as the set of cosets
	\[
		R/I = \{r+I\ |\ r \in R\}
	\]
	with addition and multiplication given by
	\[
		(r+I)+(r'+I) = (r+r')+I, \qquad (r+I)(r'+I) = rr'+I,
	\]
	for any $r,r' \in R$ and $r \in R.$
\end{defin}

Vector spaces are familiar objects of study from linear algebra courses, and they have a litany of useful properties and applications. To recall, a vector space is an abelian group with an operation defined involving a field, which we call scalar multiplication. If we relax these conditions so that scalars instead originate from a ring, we obtain an important structure which is less well-behaved than vector spaces in that there is no longer a strict notion of basis.

\begin{defin}
	Let $R$ be a ring. A \textit{left $R$-module} is an abelian group $(M,+)$ with an operation
	\[
		\cdot:
		\left\{
		\begin{array}{rcl}
			R \times M & \to & M \\
			(r,m) & \mapsto & rm
		\end{array}
		\right.
	\]
	subject to the conditions
	\begin{itemize}
		\item $1m = m$ for all $m \in M;$
		\item $(r_1r_2)m = r_1(r_2m)$ for all $r_1,r_2 \in R$ and $m \in M;$
		\item $r(m_1+m_2) = rm_1+rm_2$ for all $r \in R$ and $m_1,m_2 \in M;$
		\item $(r_1+r_2)m = r_1m+r_2m$ for all $r_1,r_2 \in R$ and $m \in M.$
	\end{itemize}
	We call this operation \textit{scalar multiplication}. An additive subgroup $S$ of $(M,+)$ is called a \textit{submodule} if it is closed under scalar multiplication, i.e. $rs \in S$ for any $r \in R$ and $s \in S,$ and we denote this by $S \leq M.$ We call $S$ a \textit{proper submodule} if $S \neq M.$
\end{defin}

A right $R$-module just has multiplication on the right as opposed to on the left. When the ring is commutative, we dispense with any mention of handedness. If the handedness of a module over a non-commutative ring is not mentioned, assume that it applies to the left (resp. right) case. Note that $R$-modules are not necessarily isomorphic to some direct sum of copies of $R$ as in the case of vector spaces. Those that do satisfy this property constitute an important family of $R$-modules.

\begin{defin}
	Let $X$ be a nonempty set and $R$ be a ring. A left $R$-module
	\[
		R\ev{X} = \left\{ \left. \sum_{i = 1}^nr_ix_i \right| n \geq 0,r_i \in R,x_i \in X \right\},
	\]
	consisting of formal linear combinations, is called a \textit{free left $R$-module generated by $X.$} $X$ is called the \textit{basis} of $R\ev{X},$ and its cardinality is called the \textit{rank} of $R\ev{X}.$
\end{defin}

The notation for a free right $R$-module generated by a set $X$ is $\ev{X}R.$ In many cases, it is unnecessary to mention handedness because it is apparent from the handedness of other modules.

\begin{rem}
	Considering that sets are bijective (isomorphic) if and only if they have the same cardinality, it is reasonable to wonder what relation there is between free $R$-modules generated by bijective sets. Harkening back to vector spaces, we may alternatively define
	\[
		R\ev{X} := \bigoplus_{x \in X}R_x,
	\]
	where $R_x = R$ for all $x \in X.$ From this definition, it is evident that if sets $X,Y$ are bijective, then $R\ev{X} \cong R\ev{Y}.$ If $X$ is finite of cardinality $n,$ then we represent the free $R$-module generated by $X$ by $R^n.$
\end{rem}

Instead of constructing a free module, we may go backwards and define when an existing module is a free module. Note that the sum mentioned is an actual sum within an abelian group as opposed to a formal sum.

\begin{defin}
	A left $R$-module $F$ is said to be \textit{free} if there exists a set $X \subseteq F$ such that for any $a \in F,$ there exists a unique collection $\{r_x\}_{x \in X}$ of scalars such that
	\[
		a = \sum_{x \in X}r_xx,
	\]
	where almost all (i.e. all but finitely many) $r_x = 0.$ If it exists, $X$ is called a \textit{basis}.
\end{defin}

In linear algebra, the action of a linear transformation was defined entirely by where it sent the basis vectors. Free modules behave similarly, and a formal description of this attribute also describes an important universal property. 

\begin{defin}
	Let $X$ be a nonempty set, $R$ be a ring, and $\iota:X \hookrightarrow R\ev{X}$ be the inclusion. For any function $f:X \to Y,$ where $Y$ is a left $R$-module, there is a unique $R$-homomorphism $g:R\ev{X} \to Y$ making the diagram
	
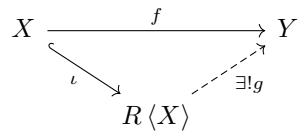
\[
		\begin{tikzcd}
			X \ar[rd,hook,"\iota"'] \ar[rr,"f"] & & Y \\
			& R\ev{X} \ar[ru,dashed,"\exists!g"'] &
		\end{tikzcd}
	\]
	commute. The universal arrow $g$ is obtained by \textit{extending $f$ by linearity}. Conversely, any left $R$-module $F$ with map $\iota':X \to F$ satisfying the universal property described above is called a \textit{free $R$-module generated by $X.$}
\end{defin}

Extending by linearity allows us to construct a canonical isomorphism between free $R$-modules generated by equinumerous sets. As such, we may talk about \emph{the} free $R$-module generated by a set $X,$ up to isomorphism. Before leaving the notion of freeness, we would be remiss to not mention how this notion generalizes outside of modules.

\begin{defin}
	Let $S:\Cs{C} \to \Cs{D}$ be a forgetful functor between appropriate categories $\Cs{C}$ and $\Cs{D}$ and $X \in \Cs{D}$ be any object. A \textit{free $\Cs{C}$ object} on $X$ is a pair $(Y,\iota),$ where $Y \in \Cs{C}$ and $\iota:X \to S(Y),$ such that for any $Z \in \Cs{C}$ and morphism $f:X \to S(Z),$ there exists a unique morphism $g:Y \to Z$ making the diagram
	\[
		\begin{tikzcd}
			X \ar[rd,"\iota"'] \ar[rr,"f"] & & S(Z) \\
			& S(Y) \ar[ru,dashed,"\exists!S(g"'] &
		\end{tikzcd}
	\]
	commute. In other words, a \textit{free $\Cs{C}$ object} on $X$ is an initial object of the comma category $X/S,$ where $X:\Cs{D} \to \Cs{D}$ is the constant functor. If $S$ has a left adjoint $F:\Cs{D} \to \Cs{C},$ then we call $F$ a \textit{free functor}, and we could alternatively define a \textit{free $\Cs{C}$ object} on $X$ to be $F(X)$ with a suitable morphism $\iota:X \to S(F(X)).$
\end{defin}

\begin{rem}
	An \textit{endomorphism} of a group $G$ is a homomorphism from $G$ to itself. The set of all endomorphisms of a group form a ring under the operations of point-wise addition and composition of functions; we call this ring the \textit{endomorphism ring of $G$} and denote it by $\ope{End}(G).$
	
	Scalar multiplication may be equivalently defined through an alternative construction mirroring the permutation representation of a group action. Let $R$ be a ring. A scalar multiplication in an $R$-module $M$ may then be defined as a ring homomorphism $R \to \ope{End}(M).$
\end{rem}

\begin{exam}
	Here are some examples of familiar structures which may be viewed as modules.
	\begin{itemize}
		\item Abelian groups may be regarded as $\Z$-modules, with submodules being subgroups.
		\item Vector spaces over a field $k$ are $k$-modules, with submodules being subspaces.
		\item Any ring $R$ may be considered as (left or right) $R$-module with scalar multiplication as ring multiplication, with submodules being ideals.
	\end{itemize}
\end{exam}

By virtue of the fact that we call modules as left/right $R$-modules, it is redundant to explicitly say that $R$ is a ring. As such, we will resort to naming the ring through its relation to a specific module. Unless stated otherwise, lack of a mention of handedness when considering non-commutative rings means that inserting either ``left'' or ``right'' will leave a true statement. Now we mention very simple methods of constructing submodules.

\begin{defin}
	Let $S \subseteq M$ be a subset of a left $R$-module, then the \textit{submodule generated by $S$}, denoted by $\langle S\rangle$, is defined as
	\[
		R\ev{S} = \langle S\rangle = \left\{\left.\sum_{i = 1}^n r_i s_i\right|r_i \in R,s_i \in S, n \geq 1\right\}.
	\]
	If $S = \{m\}$ for some $m \in M,$ then $\langle S\rangle$ is called the \textit{cyclic submodule generated by $m.$}
\end{defin}

\begin{defin}
	Let $S,T \leq M$ be submodules of a left $R$-module $M.$ Then the \textit{sum of submodules} $S$ and $T,$ denoted by $S+T,$ is defined as
	\[
		S+T = \{s+t\ |\ s \in S, t \in T\}.
	\]
	Furthermore, the \textit{intersection of submodules} $S$ and $T,$ denoted by $S \cap T,$ is defined as
	\[
		N = S \cap T = \{n\ |\ n \in S, n \in T\}.
	\]
\end{defin}

Sums and intersections of additive subgroups which are not submodules often are not submodules because subgroups are not necessarily closed under scalar multiplication. Now we present a means of imposing a module structure on a set constructed from a given submodule. This type of construction is fundamental within much of algebra.

\begin{defin}
	Let $N \leq M$ be a submodule of of the left $R$-module $M.$ Then the \textit{quotient module} $M/N$ is defined as the set of cosets
	\[
		M/N = \{m+N\ |\ m \in M\}
	\]
	with addition and scalar multiplication given by
	\[
		(m+N)+(m'+N) = (m+m')+N, \qquad r(m+N) = rm+N,
	\]
	for any $m,m' \in M$ and $r \in R.$
\end{defin}

\begin{defin}
	Let $M,N$ be left $R$-modules. A function $\varphi:M \to N$ satisfying
	\begin{itemize}
		\item $\varphi(m_1+m_2) = \varphi(m_1)+\varphi(m_2),$
		\item $\varphi(rm) = r\varphi(m),$
	\end{itemize}
	for any $r \in R$ and $m,m_1,m_2 \in M$ is called an \textit{$R$-homomorphism}. A bijective $R$-homomorphism is called an \textit{$R$-isomorphism}.
\end{defin}

It may be apparent that there is a lot of parallelism to these definitions. Explicitly, there is some type of object with structure along with functions which relate these objects by preserving structure. More importantly, these maps relating structured objects give rise to important substructures.

\begin{defin}
	Let $\varphi:M \to N$ be an $R$-homomorphism of left $R$-modules. The following distinguished subsets are, in fact, submodules.
	\begin{itemize}
		\item The \textit{kernel} of $\varphi$ is denoted by $\Ker(\varphi)$ and defined as $\Ker(\varphi) = \{m \in M\ |\ \varphi(m) = 0 \in N\}.$
		\item The \textit{image} of $\varphi$ is denoted by $\Im(\varphi)$ and defined as $\Im(\varphi) = \{\varphi(m)\ |\ m \in M\}.$
		\item The \textit{cokernel} of $\varphi$ is denoted by $\Coker(\varphi)$ and defined as $\Coker(\varphi) = N/\Im(\varphi).$
	\end{itemize}
\end{defin}

Now, we should discuss a means of combining modules. There are two ways of doing this which may seem different but which are inherently tied to one another.

\begin{defin}
	If $M,N$ are left $R$-modules, then their \textit{external direct sum}, denoted by $M \boxplus N,$ is their cartesian product $M \times N$ with coordinate-wise operations
	\begin{align*}
		(m,n)+(m',n') & = (m+m',n+n'), \\
		r(m,n) & = (rm,rn),
	\end{align*}
	where $r \in R,$ $m,m' \in M,$ and $n,n' \in N.$ Under these operations, $T = M \boxplus N$ becomes a left $R$-module.
\end{defin}

\begin{prop} \label{summ}
	The following statements are equivalent for left $R$-modules $M_k$ and $T,$ where $1 \leq k \leq n$ for some $n \in \mathbb{N}:$
	\begin{enumerate}[label = \emph{(\roman*)}]
		\item $M_1 \boxplus M_2 \boxplus \ldots \boxplus M_n \cong T;$
		\item there exist injective $R$-homomorphisms $\iota_k:M_k \to T$ such that
		\[
			T = \sum_{k = 1}^n \Im(\iota_k) \quad\text{and}\quad \Im(\iota_k) \cap \Im(\iota_j) = 0,
		\]
		where $1 \leq k \neq j \leq n;$
		\item there exist $R$-homomorphisms $\iota_k:M_k \to T$ such that, for every $t \in T,$ there are unique $m_k \in M_k$ such that
		\[
			t = \sum_{k = 1}^n \iota_k(m_k);
		\]
		\item there exist $R$-homomorphisms $\iota_k:M_k \to T$ (called injections) and $\pi_k:T \to M_k$ (called projections) such that
		\begin{itemize}
			\item $\pi_k\iota_k = \mathds{1}_{M_k},$
			\item $\sum_{k = 1}^n \iota_k\pi_k = \mathds{1}_T,$
			\item $\pi_k\iota_j = 0$ for $1 \leq k \neq j \leq n.$
		\end{itemize}
		\item there exist $R$-homomorphisms $\pi_k:T \to M_k,$ such that the map
		\[
			\psi:
			\begin{cases}
				T \to M_1 \boxplus  M_2 \boxplus \ldots \boxplus M_n \\
				t \to (\pi_1(t),\pi_2(t),\ldots,\pi_n(t))
			\end{cases}
		\]
		is an isomorphism.
	\end{enumerate}
\end{prop}

\begin{defin}
	If $M,N$ are submodules of a left $R$-module $T$ such that each $t \in T$ has a unique expression of the form $t = m+n$ for $m \in M$ and $n \in N,$ then we call $T$ the \textit{internal direct sum} of $M$ and $N,$ denoted as $T = M \oplus N.$
\end{defin}

\begin{cor}
	Let $T$ be a left $R$-module having submodules $M_k,$ where $1 \leq k \leq n$ for some $n \in \mathbb{N}.$ Then 
	\[
		T = \bigoplus_{k = 1}^nM_k
	\]
	if and only if $\sum_{k = 1}^n M_k = T$ and $M_k \cap M_j = 0$ for $1 \leq k \neq j \leq n.$ In particular,
	\[
		\bigoplus_{k = 1}^n M_k \cong M_1 \boxplus M_2 \boxplus \ldots \boxplus M_n.
	\]
\end{cor}

This corollary tells us that the external direct sum of modules is in fact an internal direct sum of them (up to isomorphism) when regarded as submodules of the resulting external direct sum. Consequently, we will dispense with the notation $\boxplus$ and let $\oplus$ represent both internal and external direct sums. Now we present a condition for when a ``summing'' collection of submodules constitute a direct sum.

\begin{prop} \label{spa}
	Let $M = S_1+\ldots+S_n$ in which $S_i \leq M$ for all $i.$ Then $M = S_1 \oplus \ldots \oplus S_n$ if and only if
	\[
		S_i \cap \left(S_1+\ldots+\hat{S_i}+\ldots+S_n\right) = 0
	\]
	for all $i,$ where $\hat{S_i}$ indicates that $S_i$ is omitted from the sum.
\end{prop}

\begin{exam}
	Let $B$ be a 2-dimensional vector space over field $k$ with basis $\{x,y\}.$ It is then a $k$-module and may be written as the direct sum $V = \ev{x} \oplus \ev{y}.$ Observe that
	\[
		\ev{x+y} \cap \ev{x} = 0 = \ev{x+y} \cap \ev{y},
	\]
	but it is not true that $V = \ev{x+y} \oplus \ev{x} \oplus \ev{y}$ because $\ev{x+y} \cap \ev{\ev{x}+\ev{y}} \neq 0.$ But note that we may write $V = \ev{x+y} \oplus \ev{x} = \ev{x+y} \oplus \ev{y},$ which goes to show that bases are not unique and that alternative bases may be easily constructed from the standard basis.
\end{exam}

\begin{exam}
	If $V$ is an $n$-dimensional vector space over field $k$ with basis $\{x_1,\ldots,s_n\},$ then $V = \ev{x_1} \oplus \ldots \oplus \ev{x_n}.$ This is the familiar statement from linear algebra that every vector may be represented uniquely by a linear combination of the basis vectors. Generalizing the previous example, we may not add any other basis vectors because we would end up representing 0 in multiple ways, implying that we cannot have any fewer basis vectors. Consequently, $V$ is a direct sum of $n$ 1-dimensional vector spaces if and only if $\ope{dim}(V) = n.$
\end{exam}

\begin{defin}
	A submodule $M$ of a left $R$-module $T$ is a \textit{direct summand} of $T$ if there exists a submodule $N$ for which $M \oplus N = T.$ If such $N$ exists, it is called the \textit{complement} of $M$ in $T.$ To denote a module as a direct summand, we write $M \suma T.$
\end{defin}

It is obvious that inclusion of submodules is a transitive relation, but it is also true that direct summands are also transitive. If $S \suma N$ and $N \suma M$ be $R$-modules, then $S \suma M$ (this can be shown by writing $M$ as a direct sum of $N$ and its complement before writing $N$ as the direct sum of $S$ and its complement).

\begin{rem}
	For a direct summand $S$ of a left $R$-module $M,$ complements of $S$ are not necessarily unique. Let $V$ be a vector space over the field $k$ with basis $\{a,b\}.$ Then $\langle\alpha a+b\rangle$ is a complement of $\langle a\rangle$ for any $\alpha \in k.$ While they may not be unique, all complements are isomorphic. This may seem strange, but the definition of internal direct sum involves an equality, not an isomorphism.
\end{rem}

There exists a special type of homomorphism which aids in detecting direct summands. This was mentioned earlier on in the category theory section and also in the splitting lemma. A proof will be provided since the result of the following lemma will be used to simplify the proof of the splitting lemma.

\begin{defin}
	A submodule $M$ of a left $R$-module $T$ is a \textit{retract} of $T$ if there exists an $R$-homomorphism $r:T \to M,$ called a \textit{retraction}, such that $r(m) = m$ for all $m \in M.$ Equivalently, $r$ is a retraction if and only if the inclusion $i:M \to T$ composes with $r$ to the identity, i.e. $ri = \mathds{1}_M.$
\end{defin}

\begin{lem} \label{ret}
	A submodule $M$ of a left $R$-module $T$ is a direct summand if and only if there exists a retraction $r:T \to M.$
\end{lem}

\begin{proof}
	Suppose that $M \leq T$ are $R$-modules such that there is a retraction $r:T \to M;$ in particular, $ri = \mathds{1}_M$ for $i:M \hookrightarrow T$ the inclusion. We claim that $T = K \oplus M,$ where $K = \Ker(r).$ We may show this using statement (ii) from \thref{summ}. Let $j:K \hookrightarrow T$ be the inclusion. Because $ri = \mathds{1}_M,$ 
	\[
		\Im(i) \cap \Im(j) = 0
	\]
	because otherwise there would be multiple elements of $M$ mapped to zero by $r$ by the definition of the kernel. On the other hand, we may write any element $t \in T$ as $t = (t-r(t))+r(t).$ Notice that $r(t) \in M$ and that $r(t) = i(r(t)).$ Furthermore $r^2 = r,$ so $t-r(t) \in K$ since $K$ is the kernel of $r.$ Similarly, we may then say that $t-r(t) = j(t-r(t)),$ so
	\[
		T = \Im(i) \oplus \Im(j)
	\]
	and we are done.
	
	Suppose that $T = M \oplus K$ is a direct sum of $R$-modules. Then for any $t \in T,$ there exist unique $m \in M$ and $k \in K$ such that $t = m+k.$ We then define
	\[
		r:
		\begin{cases}
			T \longrightarrow M \\
			m+k \mapsto m
		\end{cases}.
	\]
	For $m \in M,$ we may write $m = m+0,$ so $r(m) = m.$ As such, $r$ is a retraction.
\end{proof}

\begin{cor}
	Let $T = S \oplus N.$
	\begin{enumerate}[label = \emph{(\roman*)}]
		\item If there exists some submodule $M$ such that $S \leq M \leq T,$ then $M = S \oplus (N \cap M).$
		\item If there exists some submodule $S' \leq S,$ then $T/S' = S/S' \oplus (N+S')/S'.$
	\end{enumerate}
\end{cor}

\begin{proof}
	\thref{ret} allows us to demonstrate suitable retractions to produce the desired decompositions into direct summands.
	\begin{enumerate}[label = (\roman*)]
		\item Because we may write each $t \in T$ uniquely as a sum $t = s+n$ for $s \in S$ and $n \in N,$ we may define a map $r:T \to S$ defined by $s+n \mapsto s.$ This is clearly a retraction and remains so even if we restrict the domain to obtain the map $r|_M:M \to S.$ Restricting the domain similarly restricts the kernel, so we obtain that $\Ker(r|_M) = N \cap M,$ and we are done.
		\item Consider the map on quotient modules $\hat{r}:T/S' \to S/S'$ defined similarly to $r$ but with elements being representatives of cosets rather than individual elements. This is clearly a retraction, so we have only left to determine the kernel. Note that $\Ker(\hat{r})$ consists of all representatives $t+S'$ such that
		\[
			\hat{r}(t+S') = r(t)+S' = 0+S' \Longleftrightarrow r(t) \in S'.
		\]
		From this it is clear that all elements of the kernel take the form $(n+s')+S'$ for some $s' \in S'$ and $n \in N,$ so $\Ker(\hat{r}) = (N+S')/S'$ and we are done.
	\end{enumerate}
\end{proof}

Thus far, we have only been considering methods of combining finitely many modules together to form a new module. When our family of modules is infinite, there are two primary methods for combining the modules.

\begin{defin}
	Let $\{S_i\}_{i \in I}$ be an indexed family of sets. An \textit{$I$-tuple} is a function $f:I \to \bigcup_{i \in I} S_i$ where $f(i) \in S_i$ for all $i \in I.$ The \textit{cartesian product} $\prod_{i \in I} S_i$ is the set of all $I$-tuples.
\end{defin}

\begin{defin}
	Let $\{M_i\}_{i \in I}$ be an indexed family of left $R$-modules. The \textit{direct product} $\prod_{i \in I} M_i$ is the cartesian product $\prod_{i \in I} M_i$ with component-wise addition and scalar multiplication. The \textit{direct sum} $\bigoplus_{i \in I} M_i$ consists of all those $I$-tuples having 0 for almost all (i.e. all but finitely many) components, also with component-wise addition and scalar multiplication. The \textit{$k$-th projection} is the map
	\[
		\pi_k:
		\left\{
		\begin{array}{rcl}
			\prod_{i \in I} M_i & \longrightarrow & M_k \\
			f & \longmapsto & f(k)
		\end{array}
		\right.,
	\]
	and the \textit{$k$-th injection} is the map
	\[
		\iota_k:
		\left\{
		\begin{array}{rcl}
			M_k & \longrightarrow & \bigoplus_{i \in I} M_i \\
			m_k & \longmapsto & f_k
		\end{array}
		\right.,
	\]
	where $f_k(k) = m_k$ and 0 otherwise.
\end{defin}

\begin{rem}
	The direct sum is always a submodule of the direct product since it is a subset. Direct sums and direct products of finitely indexed families of modules are isomorphic, motivating our decision to omit mention of direct products up until this point. Unsurprisingly, direct sums over finitely indexed families coincide with the finite direct sums discussed earlier in this subsection.
\end{rem}

There is another result which generalizes \thref{spa} by using submodules generated by unions as opposed to sums. In particular, it asserts that a collection of submodules whose union spans the ambient module constitutes a direct sum precisely when it is impossible to write elements of a particular submodule as linear combinations of elements from other submodules.

\begin{prop}
	Let $\{M_i\}_{i \in I}$ be a family of submodules of a left $R$-module $T.$ If $T = \ev{\bigcup_{i \in I} S_i},$ then the following are equivalent:
	\begin{enumerate}[label = \emph{(\roman*)}]
		\item  $T = \bigoplus_{i \in I} M_i;$
		\item every $t \in T$ has a unique expression of the form $t = m_{i_1}+m_{i_2}+\cdots+m_{i_n},$ where $m_{i_j} \in M_{i_j};$
		\item $M_i \cap \ev{\bigcup_{i \neq j} M_j} = 0$ for each $i \in I.$
	\end{enumerate}
\end{prop}


\subsection{Correspondence and Isomorphism Theorems}
Correspondence theorems provide a bijective correspondence between ``containers'' of a modded out submodule and the submodules of the resulting quotient object.

\begin{theo}[Correspondence Theorem]
	Suppose that $M \leq T$ is a submodule of the left $R$-module $T.$ Let
	\[
		\mathcal{S} := \{S \leq T\ |\ M \leq S\}, \qquad \mathcal{S}/\mathcal{M} := \{S/M \leq T/M\ |\ M \leq S\}
	\]
	be the sets of intermediate submodules in $T$ and in $T/M,$ respectively. Then the function
	\[
		\phi:
		\begin{cases}
			\mathcal{S} \longrightarrow \mathcal{S}/\mathcal{M} \\
			S \longmapsto S/M
		\end{cases}
	\]
	is a bijection.
\end{theo}

The correspondence theorem for rings (resp. abelian groups) is just a corollary if you replace $T$ by a ring viewed as a left $T$-module (resp. group viewed as a $\Z$-module) and $M$ by a left ideal (resp. subgroup). We now present the primary isomorphism theorem which ties together quotient objects, kernels, and images of homomorphisms. 

\begin{theo}[First Isomorphism Theorem] \label{firso}
	If $\varphi:M \to N$ is an $R$-homomorphism of left $R$-modules, then the function
	\[
		\overline{\varphi}:
		\left\{
		\begin{array}{rcl}
			M/\Ker(\varphi) & \longrightarrow & \Im(\varphi) \\
			m+\Ker(\varphi) & \longmapsto & \varphi(m)
		\end{array}
		\right.
	\]
	is an $R$-isomorphism.
\end{theo}

\begin{rem}
	This theorem not only provides an isomorphism, but it gives a means of factoring $\varphi,$ namely as the composition $\varphi = \iota\overline{\varphi}\pi,$ where $\iota:\Im(\varphi) \hookrightarrow N$ is the inclusion and $\pi:M \to M/\Ker(\varphi)$ is the canonical projection. In particular, the projection map is determined by the kernel, so $\varphi$ is determined by its kernel.
\end{rem}

This theorem also allows us to prove a statement about all $R$-modules which serves to highlight a similarity with groups, specifically how they are all derived from free groups by adding relations.

\begin{theo} \label{res}
	Every left or right $R$-module is a quotient module of some free $R$-module.
\end{theo}

\begin{proof}
	Let $M$ be a left $R$-module and $\varphi:R\ev{M} \to M$ be the extension of the identity function $\mathds{1}:M \hookrightarrow M$ by linearity. Explicitly, $\varphi$ sends formal linear combinations in $R\ev{M}$ to actual linear combinations in $M.$ Clearly $\varphi$ is surjective, so \thref{firso} states that $M \cong F/\Ker(\varphi).$
\end{proof}

The proof shows a particular way to write $M$ as a quotient of a free module, but any surjective $R$-homomorphism $F \to M,$ for $F$ a free $R$-module, will suffice to prove the result. The following two isomorphism theorems follow from applying the first isomorphism theorem to certain maps.

\begin{cor}[Second Isomorphism Theorem]
	Suppose that $M,N \leq T$ are submodules of a left $R$-module $T.$ If $\pi:T \to T/M$ is the canonical projection, let $\pi|N:N \to T/M$ be the restriction to the submodule $N.$ Applying \thref{firso} to $\pi|N$ produces the $R$-isomorphism
	\[
		\overline{\pi|N}:N/(N \cap M) \longrightarrow (N+M)/M.
	\]
\end{cor}

\begin{cor}[Third Isomorphism Theorem]
	If $S \leq M \leq T$ is a chain of submodules of a left $R$-module $T,$ then applying \thref{firso} to the map $f:T/S \to T/M$ produces the $R$-isomorphism 
	\[
		\overline{f}:(T/S)/(M/S) \longrightarrow T/M.
	\]
\end{cor}

Again, there are analogues of all the isomorphism theorems for rings as well as abelian groups (and general groups) obtained by identifying the module structure of rings and abelian groups.


\subsection{Exact Sequences}
\begin{defin}
	A sequence of left $R$-modules and $R$-module homomorphisms
	\[
		\begin{tikzcd}
			\cdots \ar[r] & M_{n+1} \ar[r,"f_{n+1}"] & M_n \ar[r,"f_n"] & M_{n-1} \ar[r] & \cdots
		\end{tikzcd}
	\]
	is called an \textit{exact sequence} if $\Im(f_{n+1}) = \Ker(f_n)$ for all $n.$
\end{defin}

There are some specific examples of sequences which may be considered exact based on simple criteria.

\begin{prop}
	Let $A,B,C$ be left $R$-modules and $f,g,h$ be $R$-module homomorphisms.
	\begin{itemize}
		\item A sequence $0 \to A \stackrel{f}{\to} B$ is exact if and only if $f$ is a monomorphism.
		\item A sequence $B \stackrel{g}{\to} C \to 0$ is exact if and only if $g$ is an epimorphism.
		\item A sequence $0 \to A \stackrel{h}{\to} B \to 0$ is exact if and only if $h$ is an isomorphism.
	\end{itemize}
\end{prop}

\begin{defin}
	A \textit{short exact sequence} of left $R$-modules is an exact sequence of the form
	\[
		\begin{tikzcd}
			0 \ar[r] & A \ar[r,"f"] & B \ar[r,"g"] & C \ar[r] & 0,
		\end{tikzcd}
	\]
	i.e. $f$ is monic, $g$ is epic, and $\Im(f) = \Ker(g).$
\end{defin}

Short exact sequences may be very useful for stating many definitions and theorems, as will become apparent throughout this document. For now, we will state the definitions of finitely generated, finitely presented, and coherent modules as well as restate the isomorphism theorems for modules.

\begin{defin} \label{pres}
	Let $M$ be a left $R$-module. A \textit{presentation} of $M$ is an exact sequence
	\[
		\begin{tikzcd}
			\bigoplus_{\alpha \in A}R \ar[r] & \bigoplus_{\beta \in B}R \ar[r] & M \ar[r] & 0,
		\end{tikzcd}
	\]
	where $A,B$ are index sets.
\end{defin}

\begin{rem}
	Such a presentation always exists on account of \thref{res}.
\end{rem}

\begin{defin}
	With reference to the exact sequence in definition \ref{pres}, a left $R$-module $M$ is
	\begin{enumerate}
		\item \textit{finitely generated} if there exists a presentation in which $B$ is finite;
		\item \textit{finitely presented} if there exists a presentation in which both $A$ and $B$ are finite;
		\item \textit{coherent} if $M$ is finitely generated and if every finitely generated submodule $N \leq M$ is finitely presented. Alternatively $M$ is coherent if, for any submodule $N \leq M,$ $B$ being finite implies that $A$ is finite in any presentation of $N.$
	\end{enumerate}
\end{defin}

\begin{rem}
	It should be clear that $3 \Longrightarrow 2 \Longrightarrow 1$ from the definitions. It is also interesting to note the parallelism to presentations of groups. The direct sum indexed by $B$ corresponds to the generators of the free group while the direct sum indexed by $A$ corresponds to the relations which generate a normal subgroup to be modded out by.
\end{rem}

\begin{theo}
	Let all modules be $R$-modules and all maps be $R$-module homomorphisms.
	\begin{itemize}
		\item If
		\[
			\begin{tikzcd}
				0 \ar[r] & A \ar[r,"f"] & B \ar[r,"g"] & C \ar[r] & 0
			\end{tikzcd}
		\]
		is a short exact sequence, then $A \cong \Im(f)$ and $B/\Im(f) \cong C.$
		\item If $A \subseteq B \subseteq C$ is a chain of submodules, then there is a short exact sequence
		\[
			\begin{tikzcd}
				0 \ar[r] & B/A \ar[r] & C/A \ar[r] & C/B \ar[r] & 0.
			\end{tikzcd}
		\]
	\end{itemize}
\end{theo}

We know in general that the last module in a short exact sequence is isomorphic to the quotient of the second module by the image of the first module ($C \cong B/A$ based on the short exact sequence in the above definition). This brings us to a special class of short exact sequences.

\begin{defin}
	A short exact sequence
	\[
		\begin{tikzcd}
			0 \ar[r] & A \ar[r] & B \ar[r] & C \ar[r] & 0.
		\end{tikzcd}
	\]
	is called \textit{split exact} if $B \cong A \oplus C.$
\end{defin}

Split exact sequences often make computations more manageable, so it would be nice to be able to identify them. Luckily, there is a lemma which gives us some criterion for determining when a short exact sequence is split exact.

\begin{lem}[Splitting Lemma]
	Let $R$ be a ring. For a short exact sequence
	\[
		\begin{tikzcd}
			0 \ar[r] & A \ar[r,"f"] & B \ar[r,"g"] & C \ar[r] & 0.
		\end{tikzcd}
	\]
	of $R$-modules, the following statements are equivalent:
	\begin{enumerate} [label = \emph{(\roman*)}]
		\item $f$ is a section;
		\item $g$ is a retraction;
		\item there is an isomorphism $B \cong A \oplus C$ making a commutative diagram as below, where the maps are inclusions and projections.
	\end{enumerate}
	\[
		\begin{tikzcd}
			& & B \ar[dd,"\cong"] \ar[rd,"g"] & & \\
			0 \ar[r] & A \ar[ru,"f"]  \ar[rd,"\iota_1"'] & & C \ar[r] & 0 \\
			& & A \oplus C \ar[ru,"\pi_2"'] & &
		\end{tikzcd}
	\]
\end{lem}

\begin{proof}
	This proof will be completed in 3 parts.
	\begin{itemize}
		\item (i) $\Longrightarrow$ (iii): Because $f$ is a section, there exists a section $p:B \to A$ such that $pf = \mathds{1}_A.$ Define 
		\[
			h:
			\begin{cases}
				B & \to A \oplus C \\
				b & \mapsto (p(b),g(b)).
			\end{cases}
		\]
		The map $p$ is clearly epic as a left divisor of an epimorphism. Because both its components are epimorphisms, $h$ is an epimorphism. Additionally,
		\[
			\Ker(h) = \Ker(p) \cap \Ker(g) = \Ker(p) \cap \Im(f).
		\]
		Choose some $a \in A$ such that $f(a) \in \Ker(p).$ Because $pf = \mathds{1}_A,$ $p(f(a)) = 0$ if and only if $a = 0.$ Therefore, $\Ker(h) = \{f(0)\} = 0$ so that $h$ is a monomorphism. Therefore $h:B \to A \oplus C$ is an isomorphism. It is then easy to check that the triangles commute:
		\[
			\pi_2 h(b) = \pi_2(p(b),g(b)) = g(b),
		\]
		\[
			hf(a) = (pf(a),gf(a)) = (a,0) = \iota_1(a).
		\]
		
		\item (iii) $\Longrightarrow$ (i) Let $h:B \to A \oplus C$ be an isomorphism satisfying $\iota_1 = hf,$ which is guaranteed by hypothesis. Define
		\[
			p:
			\begin{cases}
				B & \to A \\
				b & \mapsto \pi_1h(b),
			\end{cases}
		\]
		where $\pi_1:A \oplus C \to A$ is the projection onto the first component. Then observe that
		\[
			pf(a) = \pi_1(hf(a)) = \pi_1\iota_1(a) = \pi_1(a,0) = a,
		\]
		so $p$ is a retraction.
		
		\item (ii) $\iff$ (iii): Because $g$ is a retraction, \thref{ret} states that
		\[
			B \cong C \oplus \Ker(g) = C \oplus \Im(f) \cong C \oplus A,
		\]
		where the equality is due to exactness and the last isomorphism follows from \thref{firso}. There must exist a section $s:C \to B$ such that $gs = \mathds{1}_C$ because $g$ is a retraction. Similar to the proof of \thref{ret}, we may write $b = (b-sg(b))+sg(b),$ where $b-sg(b) \in \Ker(g) = \Im(f).$ We then define
		\[
			h:
			\begin{cases}
				B & \to A \oplus C \\
				b & \mapsto (f^{-1}(b),g(b)),
			\end{cases}
		\]
		and it is not difficult to show that the triangles commute:
		\[
			hf(a) = (f^{-1}f(a),gf(a)) = (a,0) = \iota_1(a),
		\]
		\[
			\pi_2h(b) = \pi_2(f^{-1}(b),g(b)) = g(b).
		\]
		
		Conversely, start with $h$ as defined, which satisfies (iii), and $\iota_2:C \to A \oplus C$ the inclusion. Note that 
		\[
			gh^{-1}\iota_2 = \pi_2hh^{-1}\iota_2 = \pi_2\iota_2 = \mathds{1}_C,
		\]
		so $g$ is a retraction.
	\end{itemize}
\end{proof}

Whenever there is a sequence of objects and morphisms within a particular category $\Cs{C},$ we may apply any functor with source $\Cs{C}$ and observe how the properties of the sequence in the image differ from those of the original sequence. In particular, loss or preservation of exactness is often of importance. Here we classify some ways in which exactness may be lost or preserved within the setting of short exact sequences because it extends to exact sequences of arbitrary length.

\begin{defin}
	Let 
	\[
		\begin{tikzcd}
			0 \ar[r] & A \ar[r] & B \ar[r] & C \ar[r] & 0
		\end{tikzcd}
	\]
	be a short exact sequence in $\Cb{R\text{-}Mod}.$ A functor $F:\Cb{R\text{-}Mod} \to \Cb{Ab}$ is called
	\begin{itemize}
		\item \textit{right exact} if
		\[
			\begin{tikzcd}
				F(A) \ar[r] & F(B) \ar[r] & F(C) \ar[r] & 0
			\end{tikzcd}
		\]
		is exact,
		\item \textit{left exact} if
		\[
			\begin{tikzcd}
				0 \ar[r] & F(A) \ar[r] & F(B) \ar[r] & F(C)
			\end{tikzcd}
		\]
		is exact,
		\item \textit{half exact} if
		\[
			\begin{tikzcd}
				F(A) \ar[r] & F(B) \ar[r] & F(C)
			\end{tikzcd}
		\]
		is exact,
		\item \textit{exact} if
		\[
			\begin{tikzcd}
				0 \ar[r] & F(A) \ar[r] & F(B) \ar[r] & F(C) \ar[r] & 0
			\end{tikzcd}
		\]
		is exact.
	\end{itemize}
\end{defin}

Unsurprisingly, there are corresponding notions for contravariant functors.

\begin{defin}
	A contravariant functor $F:\Cb{R\text{-}Mod} \to \Cb{Ab}$ is called
	\begin{itemize}
		\item \textit{right exact} if
		\[
			\begin{tikzcd}
				F(C) \ar[r] & F(B) \ar[r] & F(A) \ar[r] & 0
			\end{tikzcd}
		\]
		is exact,
		\item \textit{left exact} if
		\[
			\begin{tikzcd}
				0 \ar[r] & F(C) \ar[r] & F(B) \ar[r] & F(A)
			\end{tikzcd}
		\]
		is exact,
		\item \textit{half exact} if
		\[
			\begin{tikzcd}
				F(C) \ar[r] & F(B) \ar[r] & F(A)
			\end{tikzcd}
		\]
		is exact,
		\item \textit{exact} if
		\[
			\begin{tikzcd}
				0 \ar[r] & F(C) \ar[r] & F(B) \ar[r] & F(A) \ar[r] & 0
			\end{tikzcd}
		\]
		is exact.
	\end{itemize}
\end{defin}

\begin{rem}
	It is fairly easy to keep track of what exact and half exact functors are in the two cases, but it is a bit more tricky with left and right exact functors. The simplest way to remember is to recognize that the ``functored'' exact sequence, when written with arrows going from left to right, is exact on the side suggested by the name. However, there is another interpretation which gets more to the heart of the matter.
	
	For functors, left exact refers to preservation of injectivity; on the other hand (pun intended), right exact refers to preservation of surjectivity. Dually, left exact contravariant functors convert epimorphisms into monomorphisms while the right exact ones convert monomorphisms into epimorphisms. Yet another means of describing this distinction involves kernels and cokernels, and this may be deduced easily by the reader by considering what kernels of monomorphisms and cokernels of epimorphisms are.
\end{rem}

In the next subsections, we will discuss some important functors which provide important examples of how exactness is lost or preserved.

\subsection{Hom}
For $R$-modules, Mor sets are instead called Hom sets because morphisms are $R$-module homomorphisms. In fact categories of groups, rings, or other familiar algebraic structures have the same name replacement because morphisms are just different types of homomorphisms. Let $M,N$ be left $R$-modules. Then we write $f \in \Hom{R}{M}{N}$ to mean that $f:M \to N$ is a homomorphism for which $f(rm) = rf(m)$ for any $r \in R$ and $m \in M$, so the subscript indicates from which ring scalars are taken for multiplication by elements of $M$. Something similar may be said when $M,N$ are right $R$-modules. It is interesting to note that $\Hom{R}{M}{N}$ has additional structure when $M,N$ are modules over $R$.

\begin{theo}
	Let $R$ be a ring and $M,N$ be left $R$-modules. Suppose that $f,g \in \Hom{R}{M}{N}$ are any $R$-module homomorphisms, then $\Hom{R}{M}{N}$ forms an abelian group under the binary operation $(f+g)(m) = f(m)+g(m)$. The same is also true when $M,N$ are right $R$-modules.
\end{theo}

\begin{rem}
	Associativity and commutativity follow from those properties in $N$ while the identity is the zero homomorphism (as $\Cb{R\text{-}Mod}$ is a pointed category). For a map $f$, its inverse $-f$ may be defined by $-f(m) = f((-m))$ and clearly is an $R$-module homomorphism. Therefore, Hom is just a specific kind of representable functor, the same notation will be used when scalars are unimportant or obvious.
\end{rem}

What this theorem tells us is that Hom is a functor from $\Cb{R\text{-}Mod}$ to $\Cb{Ab}.$ Even more, this functor has some useful properties which are worth discussing, especially considering that these properties manifest in other situations. But before going onto new properties, let us briefly return to something familiar.

\begin{theo}
	Both the functor $(X,-)$ and contravariant functor $(-,Y)$ are left exact.
\end{theo}

The proof is just a routine diagram chase. However, showing that these functors are not right exact requires counterexamples.

\begin{exam}
	Consider the projection map $\pi:\Z \to \Z/n\Z,$ which is clearly a surjection. Apply the functor $(\Z/n\Z,-)$ to obtain the induced map $(\Z/n\Z,\pi):(\Z/n\Z,\Z) \to (\Z/n\Z,\Z/n\Z).$ Because $\Z$ has no non-zero elements of finite order, $(\Z/n\Z,\Z) \cong 0.$ Clearly then $(\Z/n\Z,\pi)$ cannot be a surjection without also asserting that $(\Z/n\Z,\Z/n\Z) \cong 0,$ which is clearly false considering that the identity map is a non-zero element. Therefore, $(X,-)$ does not preserve surjections and therefore cannot be right exact.
	
	Consider the map $f:\Z \to \Z$ defined by $n \mapsto 2n,$ which is clearly injective. Apply $(-,\Z)$ to obtain the induced map $(f,\Z):(\Z,\Z) \to (\Z,\Z).$ Consider an element $\beta$ of the target $(\Z,\Z)$ defined as multiplication by 3. For there to exist some element $\alpha$ of the source $(\Z,\Z)$ such that $(f,\Z)(\alpha) = \beta,$ there would need to be an integer solution to the equation $3 = 2x,$ which is clearly a contradiction.
\end{exam}

\begin{defin}
	Let $F:\Cb{R\text{-}Mod} \to \Cb{Ab}$ be a covariant or contravariant functor. $F$ is called \textit{pointed} if $F(0) = 0.$
\end{defin}

\begin{exam}
	Let $M$ be any $R$-module. Because $(M,0)$ is the trivial abelian group, the functor $(M,-)$ is pointed. Similarly, $(-,M)$ is pointed. However, $M \oplus -$ is not a pointed functor because $M \oplus 0 \cong M \ncong 0,$ where the isomorphisms are between the underlying abelian groups.
\end{exam}

\begin{defin} \label{adf}
	Let $F:\Cb{R\text{-}Mod} \to \Cb{Ab}$ be a covariant or contravariant functor. $F$ is called \textit{additive} if for any $X,Y \in \Cb{R\text{-}Mod}$ and $f,g \in \Hom{R}{X}{Y},$ $F(f+g) = F(f)+F(g),$ i.e. $F_{X,Y}:\Hom{R}{X}{Y} \to \Hom{R}{F(X)}{F(Y)}$ is a group homomorphism.
\end{defin}

\begin{theo}
	Let $F:\Cb{R\text{-}Mod} \to \Cb{Ab}$ be a covariant or contravariant functor. The following are equivalent:
	\begin{enumerate}[label = \emph{(\roman*)}]
		\item $F$ is additive;
		\item for any $X,Y \in \Cb{R\text{-}Mod},$ $F(X \oplus Y) \cong F(X) \oplus F(Y);$
		\item $F$ preserves split short exact sequences. 
	\end{enumerate}
\end{theo}

\begin{rem}
	Clearly Hom is an additive functor according to how its induced maps are defined. Furthermore, there is a hierarchy between all these properties of functors. It is given by exact $\Longrightarrow$ left/right exact $\Longrightarrow$ half exact $\Longrightarrow$ additive $\Longrightarrow$ pointed.
\end{rem}

So Hom sets for $R$-modules have an abelian group structure, but it is interesting to investigate when Hom sets can have greater structure. This leads us to the following definition and theorem.

\begin{defin}
	Let $R,S$ be rings. An abelian group $M$ is said to be an $(R,S)$\textit{-bimodule} if
	\begin{itemize}
		\item $M \in \Cb{R\text{-}Mod}$,
		\item $M \in \Cb{Mod\text{-}S}$,
		\item $(\forall m \in M,\forall s \in S,\forall r \in R)$, $(rm)s = r(ms).$
	\end{itemize}
	That $M$ is an $(R,S)$-bimodule may also be conveyed by writing $_RM_S$, $M \in \Cb{R\text{-}Mod\text{-}S}$, or $M \in \Cb{(R,S)\text{-}Mod}$.
\end{defin}

\begin{exam}
	Every ring $R$ can be considered as an $(R,R)$-bimodule. Scalar multiplication on either side may simply defined as ring multiplication on either side. In this case, the compatibility condition between left and right scalar multiplication is satisfied by the associativity of ring multiplication.
\end{exam}

\begin{theo}
	Let $R,S,T$ be rings.
	\begin{itemize}
		\item $\Hom{R}{_RM}{{}_RN_S} \in \Cb{Mod\text{-}S}$ where $(f*s)(rm) = f(rm)s$.
		\item $\Hom{R}{_RM_S}{{}_RN} \in \Cb{S\text{-}Mod}$ where $(s*f)(rm) = f(rms)$.
		\item $\Hom{S}{M_S}{{}_RN_S} \in \Cb{R\text{-}Mod}$ where $(r*f)(ms) = rf(ms)$.
		\item $\Hom{S}{_RM_S}{N_S} \in \Cb{Mod\text{-}R}$ where $(f*r)(ms) = f(rms)$.
		\item $\Hom{S}{_RM_S}{{}_TN_S} \in \Cb{T\text{-}Mod\text{-}R}$ where $(t*f*r)(ms) = tf(rms)$.
		\item $\Hom{R}{_RM_S}{{}_RN_T} \in \Cb{S\text{-}Mod\text{-}T}$ where $(s*f*t)(rm) = f(rms)t$.
	\end{itemize}
\end{theo}

This theorem outlines how various combinations of bimodules can be used to impose additional structure on Hom sets. But thinking categorically, this provides a means for creating functors between module and bimodule categories. This also provides motivation for speaking about the covariant and contravariant Hom functors, which are nothing more than important instances of the covariant and contravariant representable functors.

\begin{exam}
	Let $R,S$ be rings. Then we have the following functors:
	\begin{itemize}
		\item $\Hom{R}{-}{{}_RN_S}:\Cb{R\text{-}Mod} \to \Cb{Mod\text{-}S}$
		\item $\Hom{R}{_RM_S}{-}:\Cb{R\text{-}Mod} \to \Cb{S\text{-}Mod}$
		\item $\Hom{S}{-}{{}_RN_S}:\Cb{S\text{-}Mod} \to \Cb{R\text{-}Mod}$
		\item $\Hom{S}{_RM_S}{-}:\Cb{Mod\text{-}S} \to \Cb{Mod\text{-}R}$
	\end{itemize}
\end{exam}


\subsection{Tensor Product}
\begin{defin}
	Let $R$ be a ring, $A \in \Cb{Mod\text{-}R}$, $B \in \Cb{R\text{-}Mod}$, and $G$ be any abelian group. A function $f:A \times B \to G$ between abelian groups is called
	\begin{itemize}
		\item \textit{biadditive} if $$f(a+a',b) = f(a,b)+f(a',b),$$ $$f(a,b+b') = f(a,b)+f(a,b')$$ for any $a,a' \in A$ and $b,b' \in B$;
		\item \textit{balanced} if $$f(ar,b) = f(a,rb)$$ for any $r \in R$, $a \in A$, and $b \in B$.
	\end{itemize}
\end{defin}

\begin{defin} \label{tens}
	Let $R$ be a ring, $A_R$ and $_R B$ be right and left $R$-modules, respectively. Then their \textit{tensor product} is a pair $(A \otimes_R B,\tau)$ where $A \otimes_R B$ is an abelian group and $\tau:A \times B \to A \otimes_R B$ a biadditive and balanced function such that for any abelian group $X$ and anyc biadditive and balanced function $f:A \times B \to X$, there exists a unique homomorphism of abelian groups $f':A \otimes_R B \to X$ making the following diagram commute,
	\[
		\begin{tikzcd}
			A \times B \ar[rr,"f"] \ar[rd,"\tau"'] && X \\
			& A \otimes_R B. \ar[ru,dashed,"\exists! f'"'] &
		\end{tikzcd}
	\]
\end{defin}

The elements of the tensor product take the form of finite sums of the expressions $a \otimes b = \tau(a,b)$, where $a \in A$ and $b \in B$. The tensor product seems to resemble a colimit of sorts, but it is not because the morphisms are special functions but not necessarily group homomorphisms. A colimit would have structural maps which are morphisms within the category of the diagram. We show that the tensor product of any two right and left $R$-modules exists and is unique up to isomorphism of abelian groups.

\begin{theo}
	Let $A$ and $B$ be right and left $R$-modules, respectively. Then $A \otimes_R B$ exists.
\end{theo}

\begin{proof}
	Let $R$ be a ring, $A_R$ and $_R B$ be right and left $R$-modules, respectively. Let $F$ be the free abelian group with basis $A \times B$. We define a subgroup $S \subseteq F$ consisting of all elements of the type $$(a,b+b')-(a,b)-(a,b'),$$ $$(a+a',b)-(a,b)-(a',b),$$ $$(ar,b)-(a,rb),$$ where $a,a' \in A$, $b,b' \in B$, and $r \in R$. This particular subgroup was selected because taking a quotient by this subgroup will ensure that all elements of those forms will become equivalent to 0. With this said, we define $A \otimes_R B = F/S$ and
	\[
		\tau:
		\begin{cases}
			A \times B \to A \otimes_R B \\
			(a,b) \mapsto a \otimes b.
		\end{cases}
	\]
	Making elements of the form found in $S$ equivalent to 0 then forces $\tau$ to be biadditive and balanced. Now consider the following diagram:
	
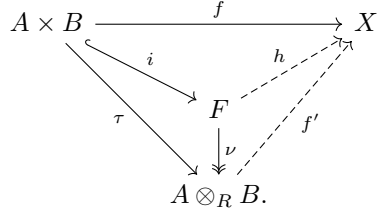
\begin{center}
		\begin{tikzcd}
			A\times B \ar[rr,"f"] \ar[rd,hook,"i"]  \ar[rdd,"\tau"'] & & X \\
			& F \ar[d,two heads,"\nu"] \ar[ru,dashed,"h"] & \\
			& A \otimes_R B. \ar[ruu,dashed,"f'"'] & 
		\end{tikzcd}
	\end{center}
	Here $i:A \times B \to F$ is the inclusion, $\nu:F \to A \otimes_R B$ the natural map from a group into its quotient, $X$ an abelian group, and $f:A \times B \to X$ an arbitrary biadditive and balanced function. Because the domain of $f$ is the set of generators of $F$, we may extend $f$ by linearity to induce a map $h:F \to X$. We may define a map $f':A \otimes_R B \to X$ by $$f'(a \otimes b) = h(a,b) = f(a,b).$$ Since $f$ is biadditive and balanced, $S \subseteq \Ker(h)$ so that $f'$ is well-defined. Note that $$h(a,b) = f'(a \otimes b) = f'((a,b) + S) = f'\nu(a,b),$$ so $h = f'\nu$. Lastly, by the definition of $\nu$, $\nu(a,b) = a \otimes b = \tau(a,b)$, so we obtain that the outer triangle is commutative, as desired.
\end{proof}

\begin{theo}
	For $A$ and $B$ right and left $R$-modules, respectively, $A \otimes_R B$ is unique up to isomorphism.
\end{theo}

\begin{proof}
	Suppose that there exists a pair $(Y,g)$ where $Y$ is an abelian group and $g:A \times B \to Y$ a biadditive and balanced function also satisfying definition \ref{tens}. Because $g$ is biadditive and balanced, there exists a unique group homomorphism $g':A \otimes_R B \to Y$, since $A \otimes_R B$ is a tensor product, making the following diagram commute:
	
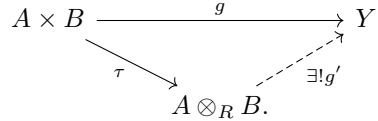
\begin{center}
		\begin{tikzcd}
			A \times B \ar[rd,"\tau"'] \ar[rr,"g"] & & Y \\
			& A \otimes_R B. \ar[ur,dashed,"\exists! g'"'] &
		\end{tikzcd}
	\end{center}
	By a similar argument, we have a unique group homomorphism $\tau':Y \to A \otimes_R B$ making the following diagram commute:
	
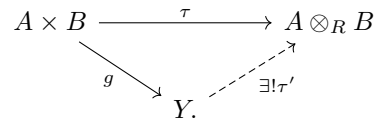
\begin{center}
		\begin{tikzcd}
			A \times B \ar[rd,"g"'] \ar[rr,"\tau"] & & A \otimes_R B \\
			& Y. \ar[ur,dashed,"\exists! \tau'"'] &
		\end{tikzcd}
	\end{center}
	By ``gluing'' the triangles together, we obtain another diagram. This diagram commutes because its component triangles individually commute. This commuting allows us to simplify the diagram as such:
	
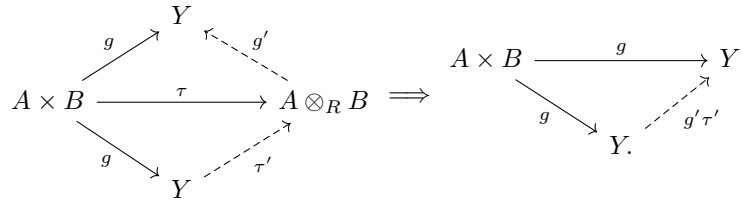
\[
		\begin{tikzcd}
			& Y & \\
			A \times B \ar[ru,"g"] \ar[rd,"g"'] \ar[rr,"\tau"] & & \ar[lu,dashed,"g'"'] A \otimes_R B \\
			& Y \ar[ru,dashed,"\tau'"'] &
		\end{tikzcd}
		\Longrightarrow
		\begin{tikzcd}
			A \times B \ar[rd,"g"'] \ar[rr,"g"] & & Y \\
			& Y. \ar[ur,dashed,"g'\tau'"'] &
		\end{tikzcd}
	\]
	But this diagram is just another tensor product triangle, and the induced map between a tensor product and itself is just the identity. Consequently, $g'\tau' = \mathds{1}_Y$ so that $g'$ is a retraction and $\tau'$ is a section. By gluing the initial triangles together across $g$ instead, we find that $\tau' g' = \mathds{1}_{A \otimes_R B}$ so that $\tau'$ is a retraction and $g'$ is a section. Then by \thref{iso}, both $\tau'$ and $g'$ are isomorphisms of abelian groups, so we are done.
\end{proof}

\begin{rem}
	Just as with Hom, we can give tensor products additional structure when we include bimodules. Let $R,S,T$ be rings. Using an obvious multiplication definition, we obtain
	\[
		_RM_S \otimes_S {}_SN_T = {}_R(M \otimes_S N)_T,
	\]
	which implies that
	\begin{itemize}
		\item $({}_RM_S \otimes_S -):\Cb{S\text{-}Mod} \to \Cb{R\text{-}Mod},$
		\item $(- \otimes_S {}_SN_T):\Cb{Mod\text{-}S} \to \Cb{Mod\text{-}T}.$
	\end{itemize}
\end{rem}

There are some other ways in which we may uncover new structure which are slightly less obvious.

\begin{defin}
	Let $f:R \to S$ be a ring homomorphism, and let $M$ be a left $S$-module. Then $M$ is a left $R$-module with multiplication defined as
	\[
		rm = f(r)m.
	\]
	We say that $_RM$ is obtained from $_SM$ through \textit{restriction of scalars}.
\end{defin}

\begin{defin}
	Let $R,S$ be rings, $f:R \to S$ be a ring homomorphism, and $M$ be a right $R$-module. We define $S$ to have the structure of a left $R$-module according to the rule
	\[
		rs = f(r)s \Longrightarrow M_R \otimes_R {}_RS_S \in \Cb{Mod\text{-}S},
	\]
	so we obtain a right $S$-module structure from an $R$-module and say that $(M \otimes_R S)_S$ is obtained from $M$ through \textit{extension of scalars}.
\end{defin}

\begin{rem}
	If $f:M_R \to M'_R,g:{}_RN \to {}_RN'$ are $R$-homomorphisms, then we may define a homomorphism of abelian groups:
	\[
		f \otimes g:M \otimes_R N \to M' \otimes_R N'
	\]
	if we using the biadditive and balanced maps guaranteed by the definition of a tensor product like so
	\[
		\begin{tikzcd}
			A \times B \ar[rr,"\tau'"] \ar[rd,"\tau"'] && M' \otimes_R N' \\
			& A \otimes_R B. \ar[ru,dashed,"\exists! f \otimes g'"'] &
		\end{tikzcd}
	\]
\end{rem}

Now that we have explained what tensor products and Hom are and why they are significant, it is fitting to conclude with a deep connection between the two functors which strongly resembles one of the examples which motivated the concept of adjunction.

\begin{theo}[Hom-Tensor Adjunction]
	Let $R,S$ be rings and $M_R,{}_RX_S,N_S$ be modules. Then there is a natural isomorphism of abelian groups:
	\[
		\begin{tikzcd}
			\Hom{S}{M \otimes_R X}{N} \ar[r,"\tau_{M,X,N}"] & \Hom{R}{M}{\Hom{S}{X}{N}}.
		\end{tikzcd}
	\]
	In particular,
	\[
		\tau_{M,X,N}(M \otimes_R X \stackrel{f}{\longrightarrow} N) =
		\left\{
		\begin{array}{rcl}
			M & \longrightarrow & \Hom{S}{X}{N} \\
			m &\longmapsto & f(m \otimes -)
		\end{array}
		\right.
	\]
	As such, there is the adjunction:
	\[
		\begin{tikzcd}
			\Cb{Mod\text{-}R} \ar[r, "(- \otimes_R X_S)"{name = L}, bend left = 25] &\Cb{Mod\text{-}S}. \ar[l, "\Hom{S}{_RX_S}{-}"{name = R}, bend left = 25] \ar[phantom, from = L, to = R, "\dashv" rotate = -90]
		\end{tikzcd}
	\]
\end{theo}


\subsection{Limits and Colimits in the Category of Modules}
	In this section, we will go over notable limits. This will be accomplished by varying the form of the family (by changing the poset) over (resp. beneath) which the limit takes form. Because $\Cb{R\text{-}Mod}$ is finitely complete, these examples will also be mentioned. When they exist, limits provide solutions to universal mapping problems, which is to say that the limit is a cone through which all other cones factor uniquely. The dual approach, obtained by adding co- wherever possible in the above statements, will also be discussed.

\subsubsection{Product and Coproduct}
Let $\Cs{C}$ be a category and $A,B \in \Cs{C}$.

\begin{defin}
	The (binary) \textit{product} of $A$ and $B$ is a triple $(A \sqcap B,\pi_A,\pi_B)$ in which $A \sqcap B \in \Cs{C}$ and $\pi_A:A \sqcap B \to A$, $\pi_B:A \sqcap B \to B$ are the \textit{projection} morphisms such that for any $X \in \Cs{C}$ and any pair of morphisms $g_A:X \to A$ and $g_B:X \to B$, there exists a unique $h:X\to C$ making the following diagram commute:
	\begin{center}
		\begin{tikzcd}
			& \ar[ld,"g_A"'] X \ar[d,dashed,"\exists! h"] \ar[rd,"g_B"] \\
			A & \ar[l,"\pi_A"] A \sqcap B \ar[r,"\pi_B"'] & B.
		\end{tikzcd}
	\end{center}
\end{defin}

\begin{defin}
	The (binary) \textit{coproduct} of $A$ and $B$ is a triple $(A \sqcup B,\iota_A,\iota_B)$ in which $A \sqcup B \in \Cs{C}$ and $\iota_A:A \to A \sqcup B$, $\iota_B:B \to A \sqcup B$ are the \textit{injection} morphisms such that for any $Y \in \Cs{C}$ and any pair of morphisms $g_A:A \to Y$ and $g_B:B \to Y$, there exists a unique $h:A \sqcup B \to Y$ making the following diagram commute:
	\begin{center}
		\begin{tikzcd}
			A \ar[r,"\iota_A"]  \ar[rd,"g_A"'] & A \sqcup B \ar[d,dashed,"\exists! h"] & \ar[l,"\iota_B"'] \ar[ld,"g_B"] B \\
			& Y &
		\end{tikzcd}
	\end{center}
\end{defin}

The definitions of (binary) product and coproduct may be obtained much more succinctly by means of the limit and colimit constructions, respectively. Let $I = \{\bullet,*\}$ be a poset with two incomparable elements and $F:I \to \Cs{C}$ be a functor such that $F(\bullet) = A$ and $F(*) = B$. Then the pair $(A,B)$ can be identified with the image of the functor $F$ so that $$A \sqcap B = \varprojlim F,$$ $$A \sqcup B = \varinjlim F.$$ Additionally, products and coproducts may be generalized to families of more than just two objects by considering a trivial poset with arbitrarily many elements.

\begin{exam}
	In $\Cb{R\text{-}Mod}$, the (co)product is called the direct product (sum).
\end{exam}

\subsubsection{Pullback and Pushout}
Consider posets $I$ and $J$ defined as,
\[
	I = \left\{
	\begin{tikzcd}
		& \bullet \ar[d] \\
		\bullet \ar[r] & \bullet
	\end{tikzcd}
	\right\}
	~~
	J = \left\{
	\begin{tikzcd}
		\ar[d] \bullet \ar[r] & \bullet \\
		\bullet & 
	\end{tikzcd}
	\right\}.
\]
Let $F:I \to \Cs{C}$ and $G:J \to \Cs{C}$ be functors whose images are defined, respectively, by

\begin{tabu} to \textwidth {XX}
	\begin{equation} \label{one}
		\begin{tikzcd}[ampersand replacement = \&]
			\& A \ar[d,"f"] \\
			B \ar[r,"g"'] \& C
		\end{tikzcd}
	\end{equation}
	&
	\begin{equation} \label{two}
		\begin{tikzcd}[ampersand replacement = \&]
			\ar[d,"g'"'] C' \ar[r,"f'"] \& A' \\
			B' \& .
		\end{tikzcd}
	\end{equation}
\end{tabu}

An image of a functor $F:I \to \Cs{C}$ is called a \textit{cospan}, and an image of a functor $G:J \to \Cs{C}$ is called a \textit{span}. In general, $I$ and $J$ are called cospan and span diagrams, respectively.

\begin{defin}
	Given a diagram of objects and morphisms in $\Cs{C}$ of the form \eqref{one}, its \textit{pullback} (or \textit{fibered product}) is a triple $(D,\alpha,\beta)$ where $D \in \Cs{C}$, $\alpha:D \to A$, and $\beta:D \to B$ such that 
	\begin{equation} \label{thre}
		f\alpha = g\beta. 
	\end{equation}
	For every similar triple $(D',\alpha',\beta')$ also satisfying \eqref{thre} (with appropriate substitutions), there exists a unique $h:D' \to D$ making the diagram
	\begin{equation} \label{fou}
		\begin{tikzcd}
			\ar[rdd,bend right,"\beta'"'] D' \ar[rd,dashed,"\exists! h"] \ar[rrd,bend left, "\alpha'"] && \\
			& \ar[d,"\beta"'] D \ar[r,"\alpha"] & A \ar[d,"f"] \\
			& B \ar[r,"g"'] & C.
		\end{tikzcd}
	\end{equation}
	commute. This diagram is referred to as a \textit{pullback square}.
\end{defin}

\begin{defin}
	Given a diagram of objects and morphisms in $\Cs{C}$ of the form \eqref{two}, its \textit{pushout} (or \textit{fibered sum}) is a triple $(D,\alpha,\beta)$ where $D \in \Cs{C}$, $\alpha:A \to D$, and $\beta:B \to D$ such that 
	\begin{equation} \label{fiv}
		\alpha f = \beta g. 
	\end{equation}
	For every similar triple $(D',\alpha',\beta')$ also satisfying \eqref{fiv} (with substitutions), there exists a unique $h:D \to D'$ making the diagram
	\begin{center}
		\begin{tikzcd}
			\ar[d,"g"'] C \ar[r,"f"] & \ar[d,"\alpha"] A \ar[rdd,bend left,"\alpha'"] & \\
			B \ar[r,"\beta"'] \ar[rrd,bend right,"\beta'"'] & D \ar[rd,dashed,"\exists! h"] & \\
			& & D'
		\end{tikzcd}
	\end{center}
	commute. This diagram is referred to as a \textit{pushout square}.
\end{defin}

As with products and coproducts, limits and colimits make the definitions of pullback and pushout, respectively, much more compact. In particular, pullbacks are limits of cospans ($\varprojlim F$). Dually, pushouts are colimits of spans ($\varinjlim G$). Additionally, pullbacks (resp. pushouts) may be generalized to families with more than just two legs on the cospan (resp. span).

\begin{theo}
	In $\Cb{R\text{-}Mod}$, pullbacks and pushouts exist.
\end{theo}

\begin{proof}
	The pullback case will be proven, but the pushout solution will be provided for verification by the reader. Let $A,B$ be $R$-modules, $f:A \to C$ and $g:B \to C$ be $R$-module homomorphisms. Define $$D = \{(a,b) \in A \oplus B\ |\ f(a) = g(b)\}.$$ Let $\alpha:D \to A$ be the projection $(a,b) \mapsto a$ restricted to domain $D$, and $\beta:D \to B$ to be the projection $(a,b) \mapsto b$ restricted to domain $D$. Let $d = (a_0,b_0)\in D$ be an element of $D$, then $$f\alpha(d) = f(a_0) = g(b_0) = g\beta(d),$$ so $(D,\alpha,\beta)$ makes the inner pullback square commute.

	Consider some $(D',\alpha',\beta')$ with sources and targets as indicated in \eqref{fou} such that the outer square commutes. Define a map $h:D' \to D$ by $$x \mapsto (\alpha'(x),\beta'(x)).$$ Then $$\alpha h(x) = \alpha(\alpha'(x),\beta'(x)) = \alpha'(x)\text{ and }\beta h(x) = \beta(\alpha'(x),\beta'(x)) = \beta'(x),$$ and so $(D,\alpha,\beta)$ is a pullback. Note that the conditions $\alpha' = \alpha h$ and $\beta' = \beta h$ uniquely determine the definition of $h$.
	
	Before writing the pushout, define $$X = \{(f(c),-g(c))\ |\ c \in C\}.$$ Then the pushout $D = (A \oplus B)/X$.
\end{proof}

\subsubsection{Kernel and Cokernel}
Let $\Cs{C}$ be a pointed category i.e. $\Cs{C}$ has a zero object. The kernel and cokernel are actually special cases of the pullback and pushout, respectively.

\begin{defin}
	If one of the legs of a cospan is a zero morphism, then the pullback is called the \textit{kernel} of the remaining morphism. In particular, the kernel of a map $f:A \to C$ is a pair $(\Ker(f),\ker(f))$ such that $$f \circ \ker(f) = 0,$$ and if $\alpha:X \to A$ is a map satisfying $$f \circ \alpha = 0,$$ then there exists a unique map $h:X \to \Ker(f)$ such that the following diagram commutes:
	\begin{center}
		\begin{tikzcd}
			\ar[rdd,bend right,"0"'] X \ar[rd,dashed,"\exists! h"] \ar[rrd,bend left, "\alpha"] && \\
			& \ar[d,"0"'] \Ker(f) \ar[r,"\ker(f)"] & A \ar[d,"f"] \\
			& 0 \ar[r,"0"'] & C.
		\end{tikzcd}
	\end{center}
\end{defin}

\begin{defin}
	If one of the legs of a span is a zero morphism, then the pushout is called the \textit{cokernel} of the remaining morphism. In particular, the cokernel of a map $f:C \to A$ is a pair $(\Coker(f),\coker(f))$ such that $$\coker(f) \circ f = 0,$$ and if $\alpha:A \to Y$ is a map satisfying $$\alpha \circ f = 0,$$ then there exists a unique map $h:\Coker(f) \to Y$ such that the following diagram commutes:
	\begin{center}
		\begin{tikzcd}
			\ar[d,"0"'] C \ar[r,"f"] & \ar[d,"\coker(f)"'] A \ar[rdd,bend left,"\alpha"] & \\
			0 \ar[r,"0"'] \ar[rrd,bend right,"0"'] & \Coker(f) \ar[rd,dashed,"\exists! h"] & \\
			& & Y
		\end{tikzcd}
	\end{center}
\end{defin}

\begin{exam}
	In $\Cb{R\text{-}Mod}$, this definition of kernels and cokernels yields the familiar submodules taught in algebra courses. In particular, $\Ker(f)$ is the submodule of $A$ mapped to 0 by $f$, and $\Coker(f)$ is the submodule of $A$ formed by the quotient $A/\Im(f)$. For a proof of these facts, alter the constructions of the pullback and pushout by incorporating the zero module and zero map.
\end{exam}

\subsubsection{Equalizer and Coequalizer}
\begin{defin}
	Let $f,g:A \to B$ be a pair of morphisms in $\Cs{C}$. The \textit{equalizer} is a pair $(\ope{Eq}(f,g),\alpha)$ where $\ope{Eq}(f,g) \in \Cs{C}$ and $\alpha:\ope{Eq}(f,g) \to A$ such that $$f\alpha = g\alpha.$$ For every similar pair $(X,\alpha')$, there exists a unique $h:X \to \ope{Eq}(f,g)$ making the following diagram commute:
	\begin{center}
		\begin{tikzcd}
			& \ar[ld,dashed,"\exists! h"'] X \ar[d,"\alpha'"] & \\
			\ope{Eq}(f,g) \ar[r,"\alpha"'] & A \ar[r,shift left,"f"] \ar[r, shift right,"g"'] & B.
		\end{tikzcd}
	\end{center}
\end{defin}

\begin{defin}
	Let $f,g:A \to B$ be a pair of parallel morphisms in $\Cs{C}$. The \textit{coequalizer} is a pair $(\ope{Coeq}(f,g),\beta)$ where $\ope{Coeq}(f,g) \in \Cs{C}$ and $\beta:B \to \ope{Coeq}(f,g)$ such that $$\beta f = \beta g.$$ For every similar pair $(Y,\beta')$, there exists a unique $h:\ope{Coeq}(f,g) \to Y$ making the following diagram commute:
	
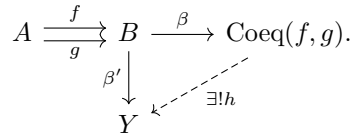
\begin{center}
		\begin{tikzcd}
			A \ar[r,shift left,"f"] \ar[r,shift right,"g"'] & \ar[d,"\beta'"'] B \ar[r,"\beta"] & \ar[ld,dashed,"\exists! h"] \ope{Coeq}(f,g). \\
			& Y &
		\end{tikzcd}
	\end{center}
\end{defin}

We may again retrieve the definitions concisely through the (co)limit construction. Consider the poset $I$ defined as,
\[
	I = \{
	\begin{tikzcd}
		\bullet \ar[r,shift left] \ar[r,shift right] & \bullet
	\end{tikzcd}
	\}.
\]
Let $F:I \to \Cs{C}$ be a functor whose image is defined by \eqref{three}.

\begin{equation} \label{three}
	\begin{tikzcd}
		A \ar[r,shift left,"f"] \ar[r,shift right,"g"'] & B
	\end{tikzcd}
\end{equation}
	Then we obtain $$\ope{Eq}(f,g) = \varprojlim F,$$ $$\ope{Coeq}(f,g) = \varinjlim F.$$ In addition, equalizers and coequalizers may be generalized to families with more than just two parallel morphisms between $A$ and $B$.

\begin{exam}
	In $\Cb{R\text{-}Mod}$, there is a simple description of the equalizer and coequalizer. Namely, $$\ope{Eq}(f,g) = \{a \in A\ |\ f(a) = g(a)\} = \Ker(f-g).$$ For the coequalizer, define $$D = \{a \in A\ |\ f(b) = a = g(b)\text{ for some }b \in B\};$$ then $\ope{Coeq}(f,g) = A/D$.
\end{exam}

Now that the most common manifestations of limits and colimits have been treated, it is worthwhile to investigate when they even exist. Recalling the definition of (co)completeness, it would be useful if there was a simpler means for checking whether a category is (co)complete. The following theorem provides an alternate method for checking finite (co)completeness.

\begin{theo} \label{comp}
	The following are equivalent for a category $\Cs{C}$:
	\begin{enumerate}[label = \emph{(\roman*)}]
		\item $\Cs{C}$ is finitely (co)complete;
		\item $\Cs{C}$ has a terminal (resp. initial) object and binary pullbacks (resp. pushouts);
		\item $\Cs{C}$ has a terminal (resp. initial) object, binary (co)products, and (co)equalizers.
	\end{enumerate}
\end{theo}

To prove this theorem, one might proceed to show that the second two criteria imply the first since the first clearly implies everything else. A proof would then assume either of the second two criteria and show that you can construct a (co)limit from an arbitrary finite diagram. When binary pullbacks exist, the unique maps to the terminal object set up the conditions for an initial pullback square. When products and equalizers exist, products are the first step before potentially introducing equalizers in case there are morphisms between the objects in the diagram. The dual approach would proceed similarly.


\newpage

\section{Special Modules and Purity}

\subsection{Projective Modules}
Recall that the covariant Hom functor $(X,-)$ in general is only left exact, so it does not preserve epimorphisms. This may make the reader curious if there is ever a case in which $(X,-)$ is an exact functor. Answering this question distinguishes a certain collection of modules.
\begin{defin} \label{proje}
	A left $R$-module $P$ is a \textit{projective module} if for any left modules $M$ and $N$, any epimorphism $f:M \to N$, and any homomorphism $g:P \to N$, there exists a homomorphism $h:P \to M$ making the following diagram commute:
	\[
		\begin{tikzcd}
			& \ar[ld,dashed,"h"'] P \ar[d,"g"]  & \\
			M \ar[r,"f"'] & N \ar[r] & 0.
		\end{tikzcd}
	\]
	Said differently, a module is projective if every homomorphism into a quotient module $N \cong M/\Ker(f)$ \textit{lifts} to a homomorphism into the base module $M.$
\end{defin}

Upon examination of this definition, it seems evident that it was formulated with the condition of making covariant Hom an exact functor in mind. For example, the bottom row of the diagram is essentially the right portion of an exact sequence. Then surjectivity of the Hom sets follows immediately because of the induced map. This, along with other equivalent formulations, is made more formal in the following theorem.

\begin{theo} \label{proj}
	For a left $R$-module $P$, the following statements are equivalent:
	\begin{enumerate}[label = \emph{(\roman*)}]
		\item $P$ is a projective module;
		\item the covariant Hom functor $\Hom{R}{P}{-}$ is exact:
		\[
			\begin{tikzcd}[cramped,sep = small]
				0 \ar[r] & A \ar[r] & B \ar[r] & C \ar[r] & 0
			\end{tikzcd} ~~
			\Longrightarrow ~~
			\begin{tikzcd}[cramped,sep = small]
				0 \ar[r] & (P,A) \ar[r] & (P,B) \ar[r] & (P,C) \ar[r] & 0;
			\end{tikzcd}
		\]
		\item any short exact sequence
		\[
			\begin{tikzcd}[cramped,sep = small]
				0 \ar[r] & A \ar[r] & B \ar[r] & P \ar[r] & 0
			\end{tikzcd}
		\]
		splits, i.e. the map $B \to P$ is a retraction so that $B \cong A \oplus P;$
		\item there exists a free left $R$-module $F$ and submodule $Q \subseteq F$ such that $F = P \oplus Q;$ informally, every projective module is a direct summand of a free module;
		\item there exists a system $(\{e_i\},\{f_i\})_{i \in I}$ of elements $e_i \in P$ and $f_i \in P^* = \Hom{R}{P}{R}$ such that, for any $x \in P,$ only finitely many of the $f_i(x)$ are nonzero and $$x = \sum_i f_i(x)e_i.$$
	\end{enumerate}
\end{theo}

\begin{cor}
	From the equivalent definitions of projective modules, the following facts follow immediately.
	\begin{itemize}
		\item Direct summands of projective modules are projective.
		\item Direct sums of projective modules are projective.
		\item All free modules are projective.
		\item Because every module is a quotient of some free module, every module is a quotient of some projective module.
	\end{itemize}
\end{cor}

\begin{defin}
	The system $(\{e_i\},\{f_i\})_{i \in I}$ mentioned in \thref{proj} is called a \textit{projective basis of $P.$} The $e_i$ are referred to as \textit{basis vectors} and the $f_i$ as \textit{coordinate functions}. Notice that if $P$ is generated by $\{e_i\}$ then $P^* = \Hom{}{P}{R}$ is generated by $\{f_i\}.$ In accordance with general modules, a projective module whose set of basis vectors is finite is \textit{finitely generated}.
\end{defin}

\begin{exam}
	We will consider a case where we have a finitely generated module which is not projective as well as a projective module which is not free. The same abelian group will be projective in one case and not in another based on the selection of ring.
	\begin{itemize}
		\item If we let our ring $R = \Z$, then a free module would be of the form $\Z^n.$ Multiplying by a positive number gives rise to the exact sequence
		\[
			\begin{tikzcd}
				0 \ar[r] & \Z \ar[r,"*5"] & \Z \ar[r] & \Z/5\Z \ar[r] & 0.
			\end{tikzcd}
		\]
		Notice this exact sequence demonstrates that $\Z/5\Z$ is finitely generated. If $\Z/5\Z$ was projective, then this exact sequence would split, which would imply that $\Z \cong \Z \oplus \Z/5\Z,$ but this is clearly not true. As a result, $\Z/5\Z$ is not projective as a $\Z$-module.
		\item If we let our ring $R = \Z/10\Z$, then a free module would be of the form $(\Z/10\Z)^n.$ Because $$\Z/5\Z \oplus \Z/2\Z \cong \Z/10\Z,$$ $\Z/5\Z$ is projective as a direct summand of a free module. However, $\Z/5\Z$ cannot be free because it has only 5 elements, where free modules would need to have $10^n$ elements. As such, $\Z/5\Z$ is projective but not free.
	\end{itemize}
\end{exam}

The next important theorem makes us of countably generated projective modules, which are defined in much a similar fashion to finitely generated projective modules.

\begin{theo}[Kaplansky]
	Every projective module is a direct sum of countably generated projective modules.
\end{theo}

\begin{lem}[Schanuel's Lemma]
	Let $R$ be a ring. If $M$ is a module, $P$ and $P'$ are projective modules, and
	\[
		\begin{tikzcd}
			0 \ar[r] & K \ar[r] & P \ar[r] & M \ar[d,equal] \ar[r] & 0 \\
			0 \ar[r] & K' \ar[r] & P' \ar[r] & M \ar[r] & 0,
		\end{tikzcd}
	\]
	is a diagram with exact rows. Then $K \oplus P' \cong K' \oplus P.$
\end{lem}

We will now see that a presentation of a projective module indicating that it is finitely generated also indicates that it is finitely presented.

\begin{cor}
	Every finitely generated projective module $P$ is finitely presented.
\end{cor}

\begin{proof}
	Because $P$ is finitely generated, we denote its basis vector set by $\{e_i\}_{i \leq n}$. We may construct a surjective map $\varphi:R^n \to P$ defined by sending the basis elements $\varepsilon_i \in R^n$ to the basis vectors of $P$ as in $\varepsilon_i \mapsto e_i$. This gives rise to the exact sequence
	\[
		\begin{tikzcd}
			0 \ar[r] & \Ker(\varphi) \ar[r] & R^n \ar[r,"\varphi"] & P \ar[r] & 0.
		\end{tikzcd}
	\]
	Since $P$ is a projective module, this is a split exact sequence. As a result, $R^n \cong \Ker(\varphi) \oplus P$. Because $R^n$ is finitely generated, $\Ker(\varphi)$ must also be finitely generated. In fact, $\Ker(\varphi)$ is projective as a direct summand of a free $R$-module. Therefore, $P$ is finitely presented.
\end{proof}

There is an additional fact about projective modules which is interesting to note. If $P$ is a finitely generated projective right $R$-module, then there exists a square matrix $E = E^2 \in M_n(R)$ such that $$P \cong R^n/ER^n = \Coker(E) = \{x + \Im(E)\ |\ x \in R^n\}.$$

\begin{defin}
	Let $M$ be an $R$-module and $K \subseteq M$ be any submodule $K \neq M$. If $N \subseteq M$ is a submodule such that $N \oplus K \neq M$ for all $K$, then $N$ is called \textit{small in $M$} or \textit{superfluous}. We denote this by $N \ll M.$
\end{defin}

\begin{defin}
	Let $M$ be a left $R$-module. If there is an exact sequence
	\[
		\begin{tikzcd}
			0 \ar[r] & \Ker(\pi) \ar[r] & P \ar[r,"\pi"] & M \ar[r] & 0,
		\end{tikzcd}
	\]
	where $P$ is a projective $R$-module and $\Ker(\pi) \ll P,$ then the pair $(P,\pi)$ is called a \textit{projective cover of $M.$}
\end{defin}

\begin{defin}
	If every finitely generated module over $R$ has a projective cover, then $R$ is called \textit{semiperfect}. If every module $R/aR$ has a projective cover, then $R$ is called \textit{semiregular}.
\end{defin}


\subsection{Injective Modules}
Just as projective modules were motivated by the condition of making $(X,-)$ exact, there is a class of modules which make $(-,X)$ exact. Not only does this definition satisfy the same condition for the dual of $(X,-)$, but the definition itself is exactly the dual of that for projective modules.

\begin{defin}
	A left $R$-module $E$ is an \textit{injective module} if for any left modules $A$ and $B$, any monomorphism $f:A \to B$, and any homomorphism $g:A \to E$, there exists a homomorphism $h:B \to E$ making the following diagram commute.
	\[
		\begin{tikzcd}
			0 \ar[r] & A \ar[d,"g"'] \ar[r,"f"] & \ar[ld,dashed,"h"] B \\
			& E &
		\end{tikzcd}
	\]
	Said differently, a module $E$ is injective if every homomorphism from a submodule to $E$ can be \textit{extended} to a homomorphism from the whole module.
\end{defin}

It should be easy to see that this does indeed make $(-,X)$ exact by reasoning dual to that used for projective modules. Injective modules also have equivalent formulations which are important.

\begin{theo}
	For a left $R$-module $E$, the following statements are equivalent:
	\begin{enumerate}[label = \emph{(\roman*)}]
		\item $E$ is an injective module;
		\item the contravariant Hom functor $\Hom{R}{-}{E}$ is exact:
		\[
			\begin{tikzcd}[cramped,sep = small]
				0 \ar[r] & A \ar[r] & B \ar[r] & C \ar[r] & 0
			\end{tikzcd} ~~
			\Longrightarrow ~~
			\begin{tikzcd}[cramped,sep = small]
				0 \ar[r] & (C,E) \ar[r] & (B,E) \ar[r] & (A,E) \ar[r] & 0;
			\end{tikzcd}
		\]
		\item any short exact sequence
		\[
			\begin{tikzcd}[cramped,sep = small]
				0 \ar[r] & E \ar[r] & A \ar[r] & B \ar[r] & 0
			\end{tikzcd}
		\]
		splits, i.e. the map $E \to A$ is a section so that $A \cong E \oplus B;$
		\item any short exact sequence
		\[
			\begin{tikzcd}[cramped,sep = small]
				0 \ar[r] & E \ar[r] & A \ar[r] & C \ar[r] & 0,
			\end{tikzcd}
		\]
		where $C$ is a cyclic module, splits;
		\item Baer Criterion: for any left ideal $I$ of $R$ and any homomorphism $f: I \to E$, there exists a homomorphism $\hat{f}:R \to E$ making the following diagram commute.
		\[
			\begin{tikzcd}
				0 \ar[r] & I \ar[d,"f"'] \ar[r] & \ar[ld,dashed,"\hat{f}"] R \\
				& E &
			\end{tikzcd}
		\]
	\end{enumerate}
\end{theo}

One of the examples of injective modules requires the following important definition.

\begin{defin}
	Let $R$ be an integral domain. Consider the equivalence relation $\sim$ on $R \times R$ defined by
	\[
		(n,m) \sim (x,y) \text{ iff } ny = mx.
	\]
	Then the \textit{field of fractions} $\ope{Frac}(R)$ is
	\[
		\ope{Frac}(R) := R \times R/\sim.
	\]
	Informally the ordered pairs represent fractions, making the equivalence relation equate fractions which are considered equivalent by cross multiplication.
\end{defin}

\begin{exam}
	The following are injective modules:
	\begin{itemize}
		\item $\Q$ as a $\Z$-module,
		\item $\Q/\Z$ as a $\Z$-module,
		\item $\ope{Frac}(R)$ as an $R$-module, when $R$ is an integral domain,
		\item any vector space over the field of fractions $\ope{Frac}(R)$ of an integral domain $R$ when viewed as an $R$-module.
	\end{itemize}
\end{exam}

\begin{rem}
	There are some useful properties enjoyed by injective modules similar to those for projective modules, which in many cases are dual those for projective modules.
	\begin{itemize}
		\item Arbitrary direct products of injective modules are injective. 
		\item Finite direct sums of injective modules are injective.
		\item Injective modules are direct summands of all modules containing them.
		\item Direct summands of injective modules are injective.
		\item Every module is a submodule of some injective module.
	\end{itemize}
\end{rem}

The dual of Schanuel's Lemma also holds.

\begin{lem}[Schanuel's Lemma]
	Let $R$ be a ring. If $M$ is a module, $E$ and $E'$ are injective modules, and the rows in the diagram
	\[
		\begin{tikzcd}
			0 \ar[r] & M \ar[d,equal] \ar[r] & E \ar[r] & L \ar[r] & 0 \\
			0 \ar[r] & M \ar[r] & E' \ar[r] & L' \ar[r] & 0
		\end{tikzcd}
	\]
	are exact. Then $E \oplus L' \cong E' \oplus L.$
\end{lem}

And unsurprisingly, there is a dual notion to projective covers.

\begin{defin}
	Let $M$ be a nonzero $R$-module. A nonzero submodule $N \subseteq M$ is said to be \textit{essential in $M$} (denoted as $N \triangleleft M$) if it has nonzero intersection with each nonzero submodule of $M.$
\end{defin}

\begin{defin}
	Let $M$ and $E$ be $R$-modules and $\alpha:M \to E$ be a homomorphism such that $\Im(\alpha) \triangleleft E.$ Then we say that $E$ is an \textit{essential extension of $M.$} If $E$ is also injective, then $E$ is called an \textit{injective envelope of $M.$} We denote this by $\ope{Env}(M).$
\end{defin}

\begin{theo}[Eckmann-Sch\"opf]
	Let $M$ be a left $R$-module. Then $\ope{Env}(M)$ always exists. Additionally if $E,E'$ are two injective envelopes of $M,$ then there exists an isomorphism $\eta:E \to E'$ for which $\eta|_M = \mathds{1}_M.$
\end{theo}


\subsection{Flat Modules}
Projective and injective modules were discovered in trying to find modules which made Hom functors exact. We saw earlier that the tensor product functor is right exact but not left exact, so it is natural to wonder if there exists some class of modules for which tensor products preserve exactness. This leads us to the following definition.

\begin{defin}
	A right $R$-module $A$ is called \textit{flat} if the exactness of the sequence of left $R$-modules
	\[
		\begin{tikzcd}
			0 \ar[r] & X \ar[r] & Y \ar[r] & Z \ar[r] & 0
		\end{tikzcd}
	\]
	implies the exactness of the sequence
	\[
		\begin{tikzcd}
			0 \ar[r] & A \otimes_R X \ar[r] & A \otimes_R Y \ar[r] & A \otimes_R Z \ar[r] & 0.
		\end{tikzcd}
	\]
\end{defin}

\begin{theo}
	For a right $R$-module $A$, the following are equivalent:
	\begin{enumerate}[label = \emph{(\roman*)}]
		\item $A$ is a flat module;
		\item $A \otimes_R -$ is an exact functor;
		\item
		\begin{tikzcd}[cramped]
			A \otimes_R I \ar[r,"\mathds{1}_A \otimes i"] & A \otimes_R R
		\end{tikzcd}
		is a monomorphism for any left ideal $I$ of $R;$
		\item
		\begin{tikzcd}[cramped]
			A \otimes_R I \ar[r,"\mathds{1}_A \otimes i"] & A \otimes_R R
		\end{tikzcd}
		is a monomorphism for any finitely generated left ideal $I$ of $R;$
		\item for any $n$-tuple $a^T = (a_1,a_2,\ldots,a_n)$ of elements in $A$ and any $X \in M_{n,d}(R)$ such that $$a^T X = 0^T,$$ there exists an $m$-tuple $\alpha^T = (\alpha_1,\alpha_2,\ldots,\alpha_m)$ of elements of $A$ and $Y \in M_{m,n}(R)$ such that
		\[
			\alpha^T Y = a^T, ~~ YX = 0.
		\]
	\end{enumerate}
\end{theo}

It should be clear that $R$ as a right $R$-module is flat because tensoring over $R$ with $R$ produces an isomorphic right $R$-module. While arbitrary direct sums of projective modules are projective if and only if each summand is projective, there is a similar statement for flat modules.

\begin{theo} \label{fla}
	A direct sum $M = \bigoplus_{k \in K} M_k$ of right $R$-modules is flat if and only if each $M_k$ is flat.
\end{theo}

\begin{cor}
	Every projective right $R$-module is flat.
\end{cor}

\begin{proof}
	Any projective right $R$-module $P$ is isomorphic to a direct summand of a free $R$-module $R^{(I)}$ for $I$ some set. Because $R$ is flat, $R^{(I)}$ must be flat by \thref{fla}. But $R^{(I)}$ can only be flat if all its direct summands is flat, so $P$ must also be flat.
\end{proof}

The reader may wonder whether the converse is true. The following example shows that it is not.

\begin{exam}
	If we let our ring be $R = \Z$, then we know that tensor products with $\Q$ over $\Z$ are trivial. For example, tensoring any exact sequence over $\Z$ by $\Q/\Z$ will annihilate the entire sequence, which is trivial case of an exact sequence. Consequently, $\Q/\Z$ is a flat $\Z$-module. Because it is a quotient, there is an obvious short exact sequence that we may form, namely
	\[
		\begin{tikzcd}
			0 \ar[r] & \Z \ar[r] & \Q \ar[r] & \Q/\Z \ar[r] & 0.
		\end{tikzcd}
	\]
	If $\Q/\Z$ was projective, then we would have that $$\Q/\Z \oplus \Z \cong \Q.$$ But because every nonzero element of $\Q/\Z$ has finite order, $(q,0) \in \Q/\Z \oplus \Z$ would have finite order whenever $q \neq 0.$ But $\Q$ has no elements of finite order, so $\Q/\Z$ cannot be projective.
\end{exam}

\begin{theo}
	Let $M$ be a right $R$-module. If every finitely generated submodule of $M$ is flat, then $M$ is flat.
\end{theo}

\begin{theo}
	Let $M$ be a finitely generated right $R$-module. Then $M$ is projective if and only if $M$ is finitely presented and flat.
\end{theo}

There are also connections between flat and injective modules. Unsurprisingly, this comes out of taking a sort of dual.

\begin{defin}
	If $M$ is a right $R$-module, we define its \textit{character module} $M^+$ as the left $R$-module $$M^+ = \Hom{\Z}{M}{\Q/\Z}.$$
\end{defin}

\begin{lem}
	A sequence of right $R$-modules
	\[
		\begin{tikzcd}
			A \ar[r,"f"] & B \ar[r,"g"] & C
		\end{tikzcd}
	\]
	is exact if and only if the induced sequence of character modules
	\[
		\begin{tikzcd}
			C^+ \ar[r,"g^+"] & B^+ \ar[r,"f^+"] & A^+
		\end{tikzcd}
	\]
	is also exact.
\end{lem}

\begin{theo}[Lambek]
	A right $R$-module $M$ is flat if and only if its character module $M^+$ is an injective left $R$-module.
\end{theo}


\subsection{Pure Exact Sequences}

In looking at projective, injective, and flat modules, our motivation was to make certain functors exact. However, we could ask the reversed question of what short exact sequences make the functors exact always. In particular, there is a class of exact sequences which make $X \otimes_R -$ exact for all right $R$-modules $X.$

\begin{defin}
	A short exact sequence
	\[
		\begin{tikzcd}
			0 \ar[r] & A \ar[r,"f"] & B \ar[r] & C \ar[r] & 0
		\end{tikzcd}
	\]
	of left $R$ modules is \textit{pure exact} if
	\[
		\begin{tikzcd}
			0 \ar[r] & X \otimes_R A \ar[r,"\mathds{1}_X \otimes f"] & X \otimes_R B \ar[r] & X \otimes_R C \ar[r] & 0
		\end{tikzcd}
	\]
	is exact for all right $R$-modules $X.$ In this case, we call $f(A) \subseteq B$ a \textit{pure submodule}.
\end{defin}

The following theorem formalizes the relationship between purity and flatness.

\begin{theo}
	A left $R$-module $Z$ is flat if and only if every exact sequence
	\[
		\begin{tikzcd}
			0 \ar[r] & X \ar[r] & Y \ar[r] & Z \ar[r] & 0
		\end{tikzcd}
	\]
	of left $R$-modules is pure exact.
\end{theo}

\begin{theo}
Every split short exact sequence is pure exact.
\end{theo}

\begin{proof}
	Let
	\[
		\begin{tikzcd}
			0 \ar[r] & A \ar[r,"i"] & B \ar[r,"p"] & C \ar[r] & 0
		\end{tikzcd}
	\]
	be a split short exact sequence. Then there exists a map $q:B \to A$ such that
	\[
		qi = \mathds{1}_A.
	\]
	Applying the functor $M \otimes_R -$ to the short exact sequence, we obtain
	\[
		\begin{tikzcd}
			M \otimes_R A \ar[r,"\mathds{1}_M \otimes i"] & M \otimes_R B \ar[r,"\mathds{1}_M \otimes p"] & M \otimes_R C \ar[r] & 0.
		\end{tikzcd}
	\]
	The composition
	\[
		(\mathds{1}_M \otimes q)(\mathds{1}_M \otimes i) = \mathds{1}_M \otimes qi = \mathds{1}_M \otimes \mathds{1}_A = \mathds{1}_{M \otimes_R A}
	\]
	shows that $\mathds{1}_M \otimes i$ is a section, so it is a monomorphism as the right divisor of an injective map (the identity). Then the tensored exact sequence is also exact, so the original exact sequence is pure.
\end{proof}

And the following theorem gives a condition for a submodule to be pure.

\begin{theo}
	Let $f:X \to Y$ be an injection of left $R$-modules. Then $f(X)$ is a pure submodule of $Y$ if and only if, given any commutative diagram with $F_0$ and $F_1$ finitely generated free $R$-modules, there is a map $F_1 \to X$ making the upper triangle commute.
	\[
		\begin{tikzcd}
			& \ar[d] F_0 \ar[r] & \ar[d] \ar[dl,dashed] F_1 \\
			0 \ar[r] & X \ar[r,"f"] & Y
		\end{tikzcd}
	\]
\end{theo}

The following lemma ties in exactness of $(X,-)$ with purity.

\begin{lem}
	Let
	\[
		\begin{tikzcd}
			0 \ar[r] & X \ar[r,"\iota"] & Y \ar[r,"\pi"] & Z \ar[r] & 0
		\end{tikzcd}
	\]
	be a pure exact sequence. From exactness, $Z \cong Y/X,$ so $\pi:Y \to Z$ is the canonical projection. If $M$ is a finitely presented left $R$-module, then $\pi_*:(M,Y) \to (M,Z)$ is surjective.
\end{lem}

\begin{rem}
	This lemma tells us that $(M,-)$ is an exact functor on pure exact sequences when $M$ is finitely presented.
\end{rem}

\begin{cor}
	Let
	\[
		\begin{tikzcd}
			0 \ar[r] & X \ar[r] & Y \ar[r] & Z \ar[r] & 0
		\end{tikzcd}
	\]
	be an exact sequence. If $Z$ is finitely presented, then this sequence is pure exact if and only if it is split.
\end{cor}

And lastly, we have a more lax condition for checking whether a sequence is pure exact.

\begin{theo}
	An exact sequence of left $R$-modules
	\[
		\begin{tikzcd}
			0 \ar[r] & X \ar[r] & Y \ar[r] & Z \ar[r] & 0
		\end{tikzcd}
	\]
	is pure exact if and only if it remains exact after tensoring by any finitely presented right $R$-module.
\end{theo}


\subsection{Torsion Product Functor}

\begin{defin}
	A pointed category $\Cs{C}$ satisfying the conditions:
	\begin{itemize}
		\item $\Hom{}{A}{B}$ is an abelian group for all $A,B \in \Cs{C},$
		\item given morphisms in $\Cs{C}$
		\[
			\begin{tikzcd}
				A \ar[r,"a"] & B \ar[r,shift left,"b"] \ar[r,shift right,"c"'] & C \ar[r,"d"] & D,
			\end{tikzcd}
		\]
		the distributive laws hold:
		\[
			d(b+c) = db+dc \quad (b+c)a = ba+ca,
		\]
		\item $\Cs{C}$ admits all finite products and coproducts,
	\end{itemize}
	is called an \textit{additive category}. A covariant or contravariant functor $F:\Cs{C} \to \Cs{D}$ between additive categories which respects the additive structure of Hom sets is called an \textit{additive functor} (cf. definition \ref{adf}).
\end{defin}

\begin{rem}
	In additive categories, finite products and coproducts coincide (up to isomorphism), so they are both written as the direct sum $\oplus.$ Conveniently, the direct sum has projection and injection maps along with conditions on their compositions in all additive categories, just as discussed in depth with modules. In addition, additive functors preserve direct sums.
\end{rem}

\begin{defin}
	An additive category $\Cs{C}$ is an \textit{abelian category} if
	\begin{itemize}
		\item all morphisms have kernels and cokernels,
		\item every monomorphism is a kernel and every epimorphism is a cokernel.
	\end{itemize}
\end{defin}

\begin{rem}
	By applying \thref{comp}, we find that all abelian categories are in fact both finitely complete and cocomplete, i.e. all finite limits and colimits exist. This makes abelian categories very convenient.
\end{rem}

\begin{exam}
	Here are some examples (and non-examples) of abelian categories.
	\begin{itemize}
		\item For any ring $R,$ both $\Cb{R\text{-}Mod}$ and $\Cb{Mod\text{-}R}$ are abelian categories. Notably, $\Cb{\Z\text{-}Mod} = \Cb{Ab}$ is the prototypical example of an abelian category.
		\item The full subcategory (i.e. restricted objects with preserved Hom sets) of $\Cb{Ab}$ consisting of all finitely generated abelian groups is abelian along with the full subcategory of all torsion abelian groups.
		\item The categories $\Cb{Grp}$ and $\Cb{ComRing}$ of groups and commutative rings, respectively, are both not abelian categories; in fact, neither of them are even additive.
		\item The full subcategory of $\Cb{Ab}$ consisting of all torsion-free abelian groups is additive but not abelian because some maps do not have cokernels, e.g. the inclusion $2\Z \hookrightarrow \Z$ has Cokernel $\Z/2\Z,$ which is not torsion-free.
	\end{itemize}
\end{exam}

This field of study heavily involves both homological and homotopical algebra, both of which grew out of algebraic topology. For the time being, we will focus on the homological aspect, which assigned a sequence of free abelian groups to a topological space that encoded specific information about the space. More generally, the sequence may be composed of free $R$-modules with little to no alteration to the subsequent theory. The goal of homological algebra in this respect was to identify and study functors that could extract invariants from these sequences, which could then be associated to a topological space and its characteristics.

The motivation for mentioning abelian categories prior to this lies in the fact that abelian categories are modeled after $\Cb{Ab},$ and in fact the Freyd-Mitchell embedding theorem asserts that any small abelian category is a full subcategory of $\Cb{R\text{-}Mod}$ for some ring $R,$ somewhat resembling the Cayley theorem of group theory which asserts that any group may be viewed as a group of symmetries. In other words, we may think of objects and morphisms in abelian categories as $R$-modules and $R$-homomorphisms. Consequently, we may study the objects in any small abelian category using the methods of homological algebra as long as we have a means of encoding information in similar sequences, which leads us to now present these ideas more rigorously.

\begin{defin}
	Let $\Cs{C}$ be an additive category. A sequence $(C_i,d_i)_{i \in \Z}$ of objects and morphisms in $\Cs{C}$ (called \textit{boundary maps})
        \[
        	(\Cb{C},d):
        	\begin{tikzcd}
			\cdots \ar[r] & C_{i+1} \ar[r,"d_{i+1}"] & C_{i} \ar[r,"d_{i}"] & C_{i} \ar[r] & \cdots
		\end{tikzcd}
        \]
        is called a \textit{chain complex} if
        \[
        	\Im(d_{i+1}) \subseteq \Ker(d_i)
        \]
        for all $i \in \Z.$ This is equivalent to saying that $d_id_{i+1} = 0$ for all $i \in \Z.$ If the boundary maps are obvious, we may notate the complex by just $\Cb{C}.$ The category of all chain complexes in $\Cs{C}$ is $\Cb{Comp}(\Cs{C}).$ Elements of $C_n$ are called \textit{$n$-chains}; elements of $\Ker(d_n)$ are called \textit{$n$-cycles}, and elements of $\Im(d_{n+1})$ are called \textit{$n$-boundaries}.
\end{defin}

\begin{prop}
	If $\Cs{C}$ is an abelian category, then $\Cb{Comp}(\Cs{C})$ is also an abelian category.
\end{prop}

\begin{rem}
	In algebraic topology, these chain complexes arise naturally from building free abelian groups from bases composed of continuous images of $n$-simplices. Finite linear combinations of simplices then may be interpreted as compiling a chain of simplices to traverse with particular orientations on each. The boundary maps are defined so that they produce lower dimensional chains which represent the boundaries of anything input into them. Most importantly, a chain which mapped to 0 through the boundary map had no boundary because it had formed a closed loop, or cycle. The condition that boundary maps compose to the zero map then has the topological significance that an $(n+1)$-simplex has a boundary which is a closed cycle and thus has no boundary of its own.
\end{rem}

\begin{defin}
	Let $(\Cb{C},d)$ and $(\Cb{C}',d')$ be chain complexes in $\Cs{C}.$ A family of morphisms $f = \{f_i|C_i \to C'_i\}_{i \in \Z}$ satisfying the relation
	\[
		f_{i-1}d_i = d'_if_i,
	\]
	which is to say that the maps commute with the boundary maps, is called a \textit{chain map}. Such a collection of maps constitutes an element $f \in \Mor{\Cb{Comp}(\Cs{C})}{\Cb{C}}{\Cb{C}'}.$
\end{defin}

\begin{rem}
	Because we may assign chain complexes to topological spaces and because topological spaces can have continuous maps between them, it is of interest how these maps manifest on the level of chain complexes. The continuous map induces a family of homomorphisms due to the continuous mapping of simplices from one space to another, and these maps commute with the boundary maps. As such, the motivation for the definition of chain maps also derives from algebraic topology.
\end{rem}

\begin{defin}
	Let $(\Cb{C},d)$ be a chain complex in $\Cs{C}.$ The quotient module
	\[
		H_n(\Cb{C}) := \frac{\Ker(d_n)}{\Im(d_{n+1})} := \frac{Z_n}{B_n}
	\]
	is called the \textit{$n$-th homology group} of $\Cb{C}.$
\end{defin}

\begin{prop}
	If $\Cs{C}$ is an abelian category, then $H_n:\Cb{Comp}(\Cs{C}) \to \Cs{C}$ is an additive functor for all $n \in \Z.$ In particular, if $f:(\Cb{C},d) \to (\Cb{C}',d')$ is a chain map, define
	\[
		H_n(f) = f_* :=
		\left\{
		\begin{array}{rcl}
			H_n(\Cb{C}) & \longrightarrow & H_n(\Cb{C}') \\
			z_n+B_n & \longmapsto & f_nz_n+B'_n
		\end{array}.
		\right.
	\]
\end{prop}

\begin{rem}
	Topologically, homology groups indicate when there is a cycle which is not itself a boundary of a higher dimensional simplex, which is in essence a hole in the topological space. Algebraically, homology groups are a measure of deviation from exactness. Notice that if the complex were in fact an exact sequence, then all the homology groups would be trivial.
\end{rem}

\begin{theo}[Long Exact Sequence in Homology]
	Let $\Cs{C}$ be an abelian category. If
	\[
		\begin{tikzcd}
			0 \ar[r] & \Cb{C}_1 \ar[r,"\iota"] & \Cb{C}_2 \ar[r,"\pi"] & \Cb{C}_3 \ar[r] & 0
		\end{tikzcd}
	\]
	is a short exact sequence in $\Cb{Comp}(\Cs{C}),$ then there is an exact sequence in $\Cs{C}$
	\[
		\begin{tikzcd}
			\cdots \ar[r] & H_{n+1}(\Cb{C}_3) \ar[r,"\partial_{n+1}"] & H_n(\Cb{C}_1) \ar[r,"\iota_*"] & H_n(\Cb{C}_2) \ar[r,"\pi_*"] & H_n(\Cb{C}_3) \ar[r,"\partial_n"] & H_{n-1}(\Cb{C}_1) \ar[r] & \cdots.
		\end{tikzcd}
	\]
	The morphisms $\partial$ are called \textit{connecting morphisms}.
\end{theo}

If there isn't enough information using just the definition of homology to carry out computations but the object of study has morphisms with better-understood objects (such as a quotient object and the corresponding Kernel), then the long exact sequence becomes an indispensable tool in completing computations. It is also useful to note that the connecting morphisms are natural with respect to complexes, which we state precisely now.

\begin{theo}
	Let $\Cs{C}$ be an abelian category. Given a commutative diagram in $\Cb{Comp}(\Cs{C})$ with exact rows,
	\[
		\begin{tikzcd}
			0 \ar[r] & \ar[d,"f"] \Cb{C}_1 \ar[r,"\iota"] & \ar[d,"g"] \Cb{C}_2 \ar[r,"\pi"] & \ar[d,"h"] \Cb{C}_3 \ar[r] & 0 \\
			0 \ar[r] & \Cb{A}_1 \ar[r,"\iota'"] & \Cb{A}_2 \ar[r,"\pi'"] & \Cb{A}_3 \ar[r] & 0,
		\end{tikzcd}
	\]
	there is a commutative diagram in $\Cs{C}$ with exact rows,
	\[
		\begin{tikzcd}
			\cdots \ar[r] & \ar[d,"h_*"] H_{n+1}(\Cb{C}_3) \ar[r,"\partial_{n+1}"] & \ar[d,"f_*"] H_n(\Cb{C}_1) \ar[r,"\iota_*"] & \ar[d,"g_*"] H_n(\Cb{C}_2) \ar[r,"\pi_*"] & \ar[d,"h_8"] H_n(\Cb{C}_3) \ar[r,"\partial_n"] & \ar[d,"f_*"] H_{n-1}(\Cb{C}_1) \ar[r] & \cdots \\
			\cdots \ar[r] & H_{n+1}(\Cb{A}_3) \ar[r,"\partial'_{n+1}"] & H_n(\Cb{A}_1) \ar[r,"\iota'_*"] & H_n(\Cb{A}_2) \ar[r,"\pi'_*"] & H_n(\Cb{A}_3) \ar[r,"\partial'_n"] & H_{n-1}(\Cb{A}_1) \ar[r] & \cdots
		\end{tikzcd}
	\]
\end{theo}

Now with the necessary terminology, we may go on to presenting a particularly useful tool for studying objects in abelian categories.

\begin{defin}
	Let $\Cs{C}$ be an abelian category and $C \in \Cs{C}$ be any object. An exact sequence of the form
	\[
		(\Cb{P},d):
		\begin{tikzcd}
			\cdots \ar[r] & P_2 \ar[r,"d_2"] & P_1 \ar[r,"d_1"] & P_0 \ar[r,"\varepsilon"] & C \ar[r] & 0,
		\end{tikzcd}
	\]
	where each $P_n$ is projective (cf. definition \ref{proje}), is called a \textit{projective resolution} of $C.$ The map $\varepsilon:P_0 \to C$ is called an \textit{augmentation}. Given a projective resolution $\Cb{P}$ of $C,$ the complex obtained by removing the augmentation
	\[
		(\Cb{P}_A,d):
		\begin{tikzcd}
			\cdots \ar[r] & P_2 \ar[r,"d_2"] & P_1 \ar[r,"d_1"] & P_0 \ar[r] & 0
		\end{tikzcd}
	\]
	is called a \textit{deleted projective resolution} of $C.$ When $\Cs{C} = \Cb{R\text{-}Mod},$ strengthening the condition that each $P_n$ is projective by requiring that each $P_n$ be free or flat results in \textit{free resolutions} and \textit{flat resolutions} of $C,$ respectively.
\end{defin}

\begin{rem}
	No information is lost in removing the augmentation because $C \cong \Coker(d_1)$ and $\varepsilon = \coker(d_1).$ We may think of resolutions as enhanced presentations, which contributes intuition to the proof of the following proposition.
\end{rem}

\begin{prop}
	Every $R$-module has a free resolution, which is necessarily a projective resolution and a flat resolution.
\end{prop}

\begin{proof}
	Let $C$ be an $R$-module. According to \thref{res}, there exists a free $R$-module $F_0$ and a short exact sequence
	\[
		\begin{tikzcd}
			0 \ar[r] & \Ker(\pi_0) \ar[r,hook,"\iota_0"] & F_0 \ar[r,two heads,"\pi_0"] & M \ar[r] & 0.
		\end{tikzcd}
	\]
	Additionally, there is a free $R$-module $F_1$ and a short exact sequence
	\[
		\begin{tikzcd}
			0 \ar[r] & \Ker(\pi_1) \ar[r,hook,"\iota_1"] & F_1 \ar[r,two heads,"\pi_1"] & \Ker(\pi_0) \ar[r] & 0.
		\end{tikzcd}
	\]
	Clearly $\Ker(\iota_0\pi_1) = \Im(\iota_1),$ so we may join the short exact sequences to form an exact sequence
	\[
		\begin{tikzcd}
			0 \ar[r] & \Ker(\pi_1) \ar[r,hook,"\iota_1"] & F_1 \ar[r,"\iota_0\pi_1"] & F_0 \ar[r,"\pi_0"] & M \ar[r] & 0.
		\end{tikzcd}
	\]
	From here, we may create another short exact sequence for $\Ker(\pi_1)$ and adjoin it in a similar way, and then continue to repeat the process indefinitely. To see that other free resolutions exist, we may instead use presentations of modules and adjoin those (albeit with a bit more work to preserve exactness).
\end{proof}

\begin{defin}
	A category $\Cs{C}$ has \textit{enough projectives} if, for each $X \in \Cs{C},$ there exists a projective object $P$ and an epimorphism $P \twoheadrightarrow X.$
\end{defin}

\begin{cor}
	If $\Cs{C}$ is an abelian category with enough projectives, then every $X \in \Cs{C}$ has a projective resolution.
\end{cor}

\begin{proof}
	Copy the prior proof but replace free modules with projective objects guaranteed by the hypothesis that $\Cs{C}$ has enough projectives.
\end{proof}

For this final part of this section, we will return to $\Cb{R\text{-}Mod}$ because we will need tensor products. This does mean that abelian categories with an analogue to tensor products can make use of the upcoming notions and definitions.

\begin{defin}
	Let $A$ and $B$ be right and left $R$-modules, respectively. Given a particular projective resolution $(\Cb{P},d)$ of $A,$ define
	\[
		\Tor{n}{R}{A}{B} := H_n(\Cb{P}_A \otimes_R B) := \frac{\Ker(d_n \otimes \mathds{1}_B)}{\Im(d_{n+1} \otimes \mathds{1}_B)}.
	\]
\end{defin}

\begin{rem}
	Although Tor is defined with reference to a particular projective resolution, we will arrive at the same computation for Tor regardless of what projective resolution is used. In particular, this means that flat resolutions may be used since they are a specific type of projective resolution. Furthermore, fixing $A$ and computing Tor with a resolution of $B$ coincides with the definition.
\end{rem}

\begin{cor}
	If $A$ and $B$ be right and left $R$-modules, then
	\[
		\Tor{0}{R}{A}{B} \cong A \otimes_R B.
	\]
\end{cor}

\begin{exam}
	Here are a couple of sample computations of torsion products.
	\item $\Tor{1}{R}{R/I}{R/J} \cong \frac{I \cap J}{IJ}.$
	\item $\Tor{1}{R}{R/aR}{R/Rb} \cong \frac{aR \cap bR}{aRb}.$
\end{exam}

\begin{theo}
	If a right $R$-module $F$ is flat, then $\Tor{n}{R}{F}{M} = 0$ for all left $R$-modules $M$ and for all $n \geq 1.$ Conversely, if $\Tor{1}{R}{F}{M} = 0$ for all $M,$ then $F$ is flat.
\end{theo}

\begin{theo}
	If $A$ and $B$ are right and left $R$-modules, respectively, then
	\[
		\Tor{n}{R}{A}{B} \cong \Tor{n}{R^{op}}{B}{A}
	\]
	for all $n \geq 0.$ Consequently, if $R$ is a commutative ring, then
	\[
		\Tor{n}{R}{A}{B} \cong \Tor{n}{R}{B}{A}
	\]
	for all $n \geq 0.$
\end{theo}

Just as with homology, there is a long exact sequence for Tor as well.

\begin{theo}
	If
	\[
		\begin{tikzcd}
			0 \ar[r] & A_1 \ar[r] & A_2 \ar[r] & A_3 \ar[r] & 0
		\end{tikzcd}
	\]
	is a short exact sequence of right $R$-modules, then there is an exact sequence for every left $R$-module $B.$
	\[
		\begin{tikzcd}[cramped,sep = small]
			\cdots \ar[r] & \Tor{n+1}{R}{A_3}{B} \ar[r] & \Tor{n}{R}{A_1}{B} \ar[r] & \Tor{n}{R}{A_2}{B} \ar[r] & \Tor{n}{R}{A_3}{B} \ar[r] & \Tor{n-1}{R}{A_1}{B} \ar[r] & \cdots.
		\end{tikzcd}
	\]
\end{theo}

Additionally, Tor preserves some important constructions.

\begin{prop}
	If $\{B_k\}_{k \in K}$ is a family of left $R$-modules then there are natural isomorphisms for all $n \geq 0,$
	\[
		\Tor{n}{R}{A}{\bigoplus_{k \in K} B_k} \cong \bigoplus_{k \in K}\Tor{n}{R}{A}{B_k}.
	\]
	The same again holds if the positions of $A$ and $B$ are swapped.
\end{prop}

\begin{prop}
	If $\{B_i,\varphi_j^i\}$ is a direct system of left $R$-modules over a directed index set $I,$ then there is an isomorphism for all $n \geq 0,$
	\[
		\Tor{n}{R}{A}{\varinjlim B_k} \cong \varinjlim \Tor{n}{R}{A}{B_k}.
	\]
	The same again holds if the roles of $A$ and $B$ are swapped.
\end{prop}

Now let us observe what happens when the coefficient ring is actually a domain.

\begin{lem}
	Let $R$ be a domain.
	\begin{itemize}
		\item If $A$ is a right torsion $R$-module, then $\Tor{1}{R}{\ope{Frac}(R)/R}{A} \cong A.$
		\item For every right $R$-module $A,$ $\Tor{n}{R}{\ope{Frac}(R)/R}{A} = 0$ for all $n \geq 2.$
		\item If $A$ is a right torsion-free $R$-module, then $\Tor{1}{R}{\ope{Frac}(R)/R}{A} \cong 0.$
	\end{itemize}
\end{lem}

\begin{theo}
	If $R$ is a domain, then $\Tor{n}{R}{A}{B}$ is a torsion module for all $A,B$ and all $n \geq 1.$
\end{theo}

Finally, we may approach Tor through an axiomatic framework so that, instead of constructing it, we can show that other constructions coincide (or differ) from it based on whether the axioms are met. Towards the end of this document, the primary result by Robinson relies on these axioms to show that he constructed something which was isomorphic to Tor.

\begin{theo}
	Let $(T_n:\Cb{R\text{-}Mod} \to \Cb{Ab})_{n \geq 0}$ be a sequence of additive covariant functors. If
	\begin{itemize}
		\item for every short exact sequence
		\[
			\begin{tikzcd}
				0 \ar[r] & A \ar[r] & B \ar[r] & C \ar[r] & 0
			\end{tikzcd}
		\]
		of right $R$-modules, there is a long exact sequence with natural connecting homomorphisms
		\[
			\begin{tikzcd}
				\cdots \ar[r] & T_{n+1}(C) \ar[r,"\Delta_{n+1}"] & T_n(A) \ar[r] & T_n(B) \ar[r] & T_n(C) \ar[r,"\Delta_n"] & T_{n-1}(A) \ar[r] & \cdots,
			\end{tikzcd}
		\]
		\item $T_0$ is naturally isomorphic to $- \otimes_R M$ for some left $R$-module $M,$
		\item $T_n(P) = 0$ for all projective right $R$-modules $P$ and all $n \geq 1,$
	\end{itemize}
	then $T_n$ is naturally isomorphic to $\Tor{n}{R}{-}{M}$ for all $n \geq 0.$
\end{theo}




\newpage

\section{Homotopy Theory}

\subsection{Compactly Generated Spaces}
The goal of this section is to construct a category of topological spaces such that, for any $X \in \Cb{Top},$ a functor $(- \times X)$ has right adjoint, i.e. the space $C(X,Y)$ has a ``nice'' topology.

\begin{defin}
	For $X,Y \in \Cb{Top},$ a function $f:X \to Y$ between their underlying sets is called \textit{$k$-continuous} if, for all compact Hausdorff spaces $C$ and all continuous maps $t:C \to X,$ the composition $tf$ is continuous.
\end{defin}

\begin{defin}
	A topological space $X$ is \textit{weak Hausdorff} if all continuous images of compact Hausdorff spaces $C$ in
$X$ are closed in $X.$ Succinctly,
	\[
		f:C \to X \text{ continuous } \Longrightarrow \Im(f) \subseteq X \text{ closed.}
	\]
\end{defin}

\begin{theo} \label{ks}
	Let $X$ be a topological space. Then the following are equivalent:
	\begin{enumerate}[label = \emph{(\roman*)}]
		\item for all spaces $Y$ and functions $f:X \to Y,$ $f$ is continuous if and only if $f$ is $k$-continuous;
		\item the previous statement holds for some family of compact Hausdorff spaces;
		\item $X$ is a quotient space of some disjoint union of compact Hausdorff spaces;
		\item any continuous map $g:C \to X$ from a compact Hausdorff space $C$ is open;
		\item $A \subseteq X$ is closed if and only if $A \cap K \subseteq K$ is closed for any compact $K \subseteq X.$
	\end{enumerate}
\end{theo}

\begin{defin}
	A space $X$ satisfying any of the equivalent conditions of \thref{ks} is called a \textit{$k$-space}. If $X$ is also weak Hausdorff, then we say that $X$ is \textit{compactly generated}.
\end{defin}

\begin{rem}
	The continuous image of a compact Hausdorff space $C$ inside a compactly generated space $X$ must be clopen (i.e. closed and open), so the existence of a non-surjective continuous map $g:C \to X$ implies that $X$ has more than one connected component. With regard to the separation axioms, $T_2 \Longrightarrow \text{weak Hausdorff} \Longrightarrow T_1.$ Recall that $T_2$ spaces are identically Hausdorff and that $T_1$ spaces are defined by the property that any pair of points $x,y$ have open sets $U,V$ such that $x \in U,$ $y \in V,$ $x \notin V,$ and $y \notin U.$
\end{rem}

\begin{exam}
	Some examples of compactly generated spaces include
	\begin{itemize}
		\item compact Hausdorff spaces,
		\item locally compact Hausdorff spaces,
		\item CW-complexes,
		\item first countable spaces, i.e. spaces where each point has a countable base.
	\end{itemize}
\end{exam}

\begin{theo}
	We denote the category of $k$-spaces with continuous maps between them by $\Cb{kTop}.$ Additionally, we denote the category of all topological spaces with $k$-continuous maps between them by $\Cb{Top_k}.$ Then we have
	\[
		\begin{tikzcd}
			\Cb{kTop} \ar[rd,"\approx"'] \ar[r,hook] & \ar[d] \Cb{Top} \\
			& \Cb{Top_k};
		\end{tikzcd}
	\]
	recall that $\approx$ denotes an equivalence of categories.
\end{theo}

\begin{defin}
	For $X \in \Cb{Top},$ we define the \textit{kaonization} (\textit{k-ification}) of $X$ to be a topological space $k(X)$ having the same underlying set as $X$ but with the condition that a subset $A \subseteq k(X)$ is closed if and only if it is closed in any compact subspace $K \subseteq X.$ Succinctly,
	\[
		A \subseteq k(X) \text{ closed } \iff A \cap K \subseteq K \text{ closed, } \forall K \subseteq X \text{ compact.}
	\]
\end{defin}

\begin{rem}
	In finding the kaonization of a topological space$(X,\tau),$ we are enriching $\tau$ by adding more closed sets, and hence open sets, in order to generate the coarsest topology $\tau'$ relative to $\tau$ such that $\tau'$ makes $(X,\tau')$ a $k$-space (which is in fact even compactly generated).
\end{rem}

\begin{prop}
	For any space $X:$
	\begin{itemize}
		\item $k(X) \stackrel{\mathds{1}}{\longrightarrow} X$ is continuous;
		\item $k(X)$ is compactly generated;
		\item $k(X) = X$ if and only if $X$ is compactly generated;
		\item if $Y$ is compactly generated, then
		\begin{align*}
			Y \to k(X) \text{ continuous } & \iff Y \to X \text{ continuous}, \\
			\Mor{\Cb{kTop}}{Y}{k(X)} & \longleftrightarrow \Mor{\Cb{Top}}{Y}{X}, \\
			(f:Y \to k(X)) & \longmapsto (\mathds{1}f:Y \to k(X) \to X);
		\end{align*}
		\item if $\iota:\Cb{CGTop} \to \Cb{Top}$ is an inclusion functor with source the category of all compactly generated spaces with continuous maps between them, then $\iota/X$ has final object $\mathds{1}:k(X) \to X;$
		\item $k(X)$ and $X$ have the same compact subsets;
		\item $k(X)$ and $X$ are weakly homotopy equivalent, i.e. they have isomorphic homotopy groups. Moreover, if $Y = \Delta^n$ is an $n$-simplex, then singular simplices are the ``same.'' Therefore, \emph{(}co\emph{)}homology theories of $k(X)$ and $X$ coincide.
	\end{itemize}
\end{prop}

\begin{rem}
	The identity map $k(X) \to X,$ which you may recall is continuous, acts as the counit map for an adjunction of the inclusion and kaonization functors. More specifically, the functors 
	\[
		\iota:\Cb{kTop} \hookrightarrow \Cb{Top} \text{ and } k:\Cb{Top} \longrightarrow \Cb{kTop}
	\]
	comprise the adjunction $\iota \dashv k.$ Consequently, $k$ preserves limits and $\iota$ preserves colimits. In particular, $k$ preserves products. If we define $X \times_k Y := k(X \times Y),$ then we obtain that
	\[
		k(X) \times k(Y) \cong X \times_k Y.
	\]
	In fact, if $X$ is locally compact and $Y$ is compactly generated, then
	\[
		X \times_k Y = X \times Y.
	\]
\end{rem}

\begin{defin}
	A category $\Cs{C}$ with finite products and terminal object is said to be \textit{cartesian closed} if the functor 
	\[
		- \sqcap X:\Cs{C} \to \Cs{C}
	\]
	has right adjoint, which is called the \textit{exponentiation functor} and denoted by $(-)^X.$ The terminal object $T \in \Cs{C}$ has the property that $A \sqcap T \cong A$ for any $A \in \Cs{C}.$ This further implies that $A^T \cong A,$ which may be shown. Adjunction produces the bijection
	\[
		(A \sqcap T,B) \leftrightarrow (A,B^T) \Longrightarrow (A,B) \leftrightarrow (A,B^T) \Longrightarrow B \cong B^T,
	\]
	where the last implication follows from the Yoneda Lemma.
\end{defin}

\begin{defin}
	Let $X,Y$ be topological spaces. We imbue the set $C(X,Y)$ of continuous maps from $X$ to $Y$ with a topology whose basis consists of the sets
	\[
		U(K,V) := \{f:X \to Y \text{ continuous }|f(K) \subseteq V\},
	\]
	where $K \subseteq X$ compact and $V \subseteq Y$ open. We call this topology the \textit{compact-open topology}.
\end{defin}

\begin{rem}
	$\Cb{kTop} \approx \Cb{Top_k}$ are cartesian closed. Exponential objects are $kC(X,Y),$ i.e. all $k$-continuous (and thereby continuous since $X$ is a $k$-space) maps with the compact-open topology.
\end{rem}

\begin{defin}
	Let $X = \varinjlim X_n$ in $\Cb{Top}.$ We say that $A \subseteq X$ open if $A \cap f_n(X_n) \subseteq f_n(X_n)$ open for each $n,$ where $f_n:X_n \to X$ are structural maps. This defines the \textit{weak topology} on $X.$
\end{defin}


\subsection{Geometric Realization}
The category $\Delta$ generalizes the topological notion of simplicial complexes and provides a completely algebraic method for developing homotopy and homology theories, which is a common theme in abstract homotopy theory.

\begin{defin}
	Let $\Delta$ be a category of finite ordinals $[n]$ (cf. $\Delta[n]$ in \thref{pos}) with order-preserving (equivalently non-decreasing) maps $f:[n] \to [m].$ A contravariant functor $F:\Delta \to \Cb{Set}$ is called a \textit{simplicial set}. From this we obtain a category $\Cb{sSet} = (\Delta^{op},\Cb{Set})$ of simplicial sets with simplicial maps between them.
\end{defin}

Here we should note well that $\Cb{sSet}$ is a functor category, whose objects are simplicial sets. Harkening back to the definition of a functor category, this means that simplicial maps are in fact natural transformations between simplicial sets.

\begin{defin}
	The map $d_i^{n+1}:[n] \to [n+1]$ defined by
	\[
		d_i^{n+1}(j):
		\begin{cases}
			j, & 0 \leq j < i \\
			j+1, & i \leq j \leq n
		\end{cases}
	\]
	is called a \textit{coface map}. On the other hand, the map $s_i^{n}:[n+1] \to [n]$ defined by
	\[
		s_i^{n}(j):
		\begin{cases}
			j, & 0 \leq j \leq i \\
			j-1, & i < j \leq n+1
		\end{cases}
	\]
	is called a \textit{codegeneracy map}.
\end{defin}

The coface map $d_i$ fixes all elements below $i$ and raises those elements at least $i,$ effectively leaving a blank in the position that $i$ used to be in. This indicates $[n]$ as the $i$-th face of $[n+1].$ The codegeneracy map $s_i$ fixes all elements at most $i$ and lowers those elements above $i,$ leading to consolidation of the element $(i+1)$ into $i.$ Considering that $i+1 = i$ indicates a degenerate simplex, $s_i$ acts as though it is converting a degenerate $[n+1]$ simplex into its corresponding $[n]$ simplex.

\begin{theo} \label{simp}
	The coface and codegeneracy maps satisfy the following relations:
	\begin{align*}
		d_i^{n+1}d_j^{n} &= d_{j+1}^{n+1}d_i^n \quad\ \text{ for } i \leq j, \\
		s_j^{n-1}s_i^n &= s_i^{n-1}s_j^n \quad\ \text{ for } i \leq j, \\
		s_j^{n-1}d_i^n &= 
		\begin{cases}
			d_i^{n-1}s_{j-1}^{n-2}, & i < j \\
			\mathds{1}_{[n-1]}, & i = j \text{ or } i = j+1 \\
			d_{i-1}^{n-1}s_j^{n-2}, & i > j+1
		\end{cases}
	\end{align*}
	Moreover, any order-preserving map $f:[n] \to [m]$ may be written uniquely as an expression of the form
	\[
		f = d_{i_k}^m\ldots d_{i_1}^{n-h+1}s_{j_1}^{n-h}\ldots s_{j_h}^{n-1},
	\]
	where $k$ and $h$ are non-negative integers such that
	\begin{itemize}
		\item $n+k-h = m,$
		\item $0 \leq i_1 < \ldots < i_k \leq m,$
		\item $0 \leq j_1 < \ldots < j_h \leq n,$
		\item $i_1,\ldots,i_k$ are the elements of $[m]$ not in the image of $f,$
		\item $j_1,\ldots,j_h$ are the elements of $[n]$ at which $f$ does not increase (as in $f(j) = f(j+1)$).
	\end{itemize}
\end{theo}

\begin{rem}
	Here it is important to note that, while the functions above are for $\Delta,$ a simplicial set is a contravariant functor, so there are dual maps to the coface and codegeneracy maps which manifest in simplicial sets in a manner that may be familiar to those who have worked with simplicial complexes.
\end{rem}

\begin{defin}
	Explicitly, a \textit{simplicial set} $X \in (\Delta^{op},\Cb{Set})$ is a sequence of sets $\{X_k = X[k]\}_{k \in \mathbb{N}}$ together with the maps
	\[
		\partial_i^n = X(d_i^n), \qquad \sigma_j^n = X(s_j^n)
	\]
	which are called \textit{face} and \textit{degeneracy} maps, respectively. The elements of $X_0$ are called \textit{vertices}, of $X_1$ are called \textit{edges}, of $X_n$ are called \textit{$n$-simplices of $X.$}
\end{defin}

By the duality principle, the face and degeneracy maps must satisfy relations dual to those in \thref{simp}. This theorem also implies that simplicial maps may be written uniquely as a composition of degeneracy and face maps.

\begin{rem}
	The process described in this section is reversible. That is, given a sequence of sets $X_0,X_1,\ldots,X_n,\ldots$ and maps $\partial_i^n,\sigma_j^n$ satisfying the condition of \thref{simp}, it is possible to construct a functor $X:\Delta^{op} \to \Cb{Set}$ such that $X[i] = X_i,$ $X(d_i^n) = \partial_i^n,$ and $X(s_j^n) = \sigma_j^n$ for $i,j \geq 0.$
\end{rem}


There is a method of ``topologizing'' categories by means of the following construction.

\begin{rem}
	Given a category $\Cs{C}$, its \textit{classifying space} is a CW-complex $B\Cs{C}$ constructed in the following manner:
	\begin{enumerate}
		\item The 0-cells (vertices) are objects of $\Cs{C}$.
		\item The 1-cells (edges) are morphisms of $\Cs{C}$ (excluding the identities) which attach at the source and target of the morphism.
		\item For each ordered pair $(f,g)$ of composable non-identity morphisms in $\Cs{C}$, attach a 2-simplex $\Delta^2$ as indicated by the following diagram:
		
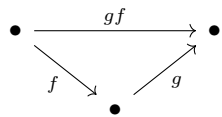
\begin{center}
			\begin{tikzcd}
				\bullet \ar[rr,"gf"] \ar[rd,"f"'] & & \bullet \\
				& \bullet \ar[ru,"g"'] &
			\end{tikzcd}
		\end{center}
		\item Continue this process inductively by attaching an $n$-simplex $\Delta^n$ when given an $n$-tuple of composable non-identity morphisms.
		\item $B\Cs{C}$ is the union of all these spaces with the weak topology, which is the topology in which a subset is open if and only if its intersection with each of the subspaces above is open.
	\end{enumerate}
\end{rem}

\begin{exam}
	Some illustrative examples of classifying spaces include the following.
	\begin{itemize}
		\item $B\{\bullet\}$ is a point.
		\item $B\{0 \leq 1\} = [0,1]$.
		\item $B\{0 \leq 1 \leq 2\} = \Delta^2$.
		\item $B\Delta[n] = \Delta^n$.
		\item $B\{\bullet \rightrightarrows \bullet\} = S^1$.
	\end{itemize}
\end{exam}


\subsection{Homotopy Groups}
Topological spaces formalize and generalize continuity as we intuitively conceive of it. From this, we are able to talk about what it means for two spaces to be topologically equivalent, which is referred to as homeomorphism. As a reminder, a homeomorphism is a continuous bijection with the extra condition that its inverse is also continuous. It turns out this is a fairly strong requirement for something that is colloquially referred to as ``rubber-sheet geometry.'' For example, homeomorphic spaces must share properties such as dimension, cardinality, compactness, and separability even though a rubber sheet should allow us to squish an open ball to a closed ball contained within or even to a point in an obvious fashion, leading to a change in compactness or even dimension.

So what we seek is a weaker means of distinguishing spaces which preserves properties consistent with this rubber sheet perspective. Here we begin with the notion of homotopy and build up to a suitable means of distinguishing spaces. But first it is necessary to outline precisely what types of topological spaces are necessary for homotopy to be well-defined; in essence, we must set out by outlining the operative category.

\begin{defin}
	A \textit{pointed topological space} is a pair $(X,x_0),$ where $X$ is a space and $x_0 \in X$ is a distinguished point. The base point is often not explicitly written when it is clear that the space is based, so we may refer to $(X,x_0)$ as just $X.$ Then a continuous map $f:X \to Y$ between based spaces is a \textit{pointed map} if $f(x_0) = y_0.$ The category of pointed spaces with pointed (continuous) maps is notated as $\Cb{Top_*}.$
\end{defin}

For the remainder of this section, assume that all spaces and maps between them belong within $\Cb{Top_*}$ unless indicated otherwise. Furthermore, we will use $I$ to denote the closed interval $[0,1].$

\begin{defin}
	Let $f,g:X \to Y$ be pointed maps between the spaces $X$ and $Y.$ A continuous map $h:X \times I \to Y$ satisfying the conditions
	\[
		h(x,0) = f(x) \quad\quad\quad h(x,1) = g(x) \quad\quad\quad h(x_0,t) = y_0
	\]
	is called a \textit{homotopy} from $f$ to $g.$ If there exists a homotopy from $f$ to $g,$ we say that $f$ is \textit{homotopic} to $g,$ notated as $f \simeq g.$ A map which is homotopic to the constant map, whose image must be $\qty{y_0}$ in $\Cb{Top_*},$ is called \textit{nullhomotopic}. If the identity $\mathds{1}_X$ is nullhomotopic, then we say that $X$ is \textit{contractible}.
\end{defin}

We may conceive of a homotopy $h:f \simeq g$ between pointed maps $f,g:X \to Y$ as a map which assigns to each point of $x \in X$ a path $p_x:I \to Y$ beginning in the image of $f$ and ending in the image of $g$ with the condition that $p_{x_0}$ be the constant path. With this insight, it becomes clear that path-connectedness is a very convenient property when dealing with homotopy. For one, having path-connectedness makes the choice of base point trivial, and it is also fairly easy to see that contractible spaces must always be path-connected.

\begin{rem}
	The third condition listed for a map to be a homotopy serves the purpose of ensuring that for any choice of $t \in I,$ that the intermediate map $h(-,t):X \to Y$ is a pointed map. However, the map $h(x,-):I \to Y$ need not be pointed unless $x = x_0.$
	
	Also, the alternative view of homotopy as a path-assigning function should remind us of a collection of parameterized univariate maps. Equipped with the motivation given for adjunction, we can reasonably predict that some kind of adjunction is at play, and in fact there is (c.f. loop space/reduced suspension adjunction from \thref{adj}).
\end{rem}

\begin{exam}
	Provided here are some examples in real space which serve to develop geometric intuition for the definitions just given.
	\begin{itemize}
		\item In $\R^n$ with arbitrary base point, any two paths $f,g:I \to \R^n$ are homotopic by way of the ``straight-line'' homotopy, which is given by
		\[
			h:
			\left\{
			\begin{array}{rcl}
				I \times I &\to & \R^n \\
				(x,t) &\mapsto & (1-t)f(x)+tg(x).
			\end{array}
			\right.
		\]
		In some multivariable calculus courses, it may have been mentioned that spaces with this property are called simply connected. This will be defined more precisely further on, but this example serves as a useful heuristic.
		\item The straight-line homotopy on $\R^n$ may actually be used to show that any map is nullhomotopic and, in particular, that $\R^n$ is contractible. As we will soon see, this suggests the curious result that $\R^n$ behaves like a point with respect to homotopy. And because $\C^n \cong \R^{2n}$ are homeomorphic, anything statement about homotopy in real space extends to complex space.
		\item The seemingly similar space $\R^n-\vb{0}$ with arbitrary base point actually exhibits very different properties for $n = 1,$ $n = 2,$ and $n \geq 3.$
		\begin{itemize}
			\item For $n = 1,$ the space consists of 2 connected components (which are also path-components). In fact, $\R-\vb{0}$ is homeomorphic to two copies of $\R$ by way of the maps $e^x$ and $e^{-x}$ defined on each of the copies of $\R.$ Consequently, $\R-\vb{0}$ behaves like 2 points up to homotopy since each copy of $\R$ is contractible.
			\item For $n = 2,$ the space is path-connected, and many maps are nullhomotopic. However, consider the map 
			\[
				f:x \mapsto (\cos(2\pi x),\sin(2\pi x))
			\]
			on $I.$ It is a fact---intuitive but not obvious to show---that this map is not nullhomotopic and that every loop (path with matching endpoints) is either nullhomotopic or homotopic to $f.$ With that said, it is clear then that $\R^2$ is no longer contractible or simply connected if you remove a point from it.
			\item For $n \geq 3,$ the space is again path-connected. In contrast with $\R^2-\vb{0},$ this space is simply connected, but that does not mean that all maps are nullhomotopic as it is not hard to convince ourselves intuitively that this space is still not contractible. The distinction is between a generic map (with arbitrary domain) and a path (a map with domain $I$).
		\end{itemize}
	\end{itemize}
\end{exam}

Homotopies may be defined in $\Cb{Top}$ by omitting the condition that base points be sent to base points since we do not have base points to speak of. However, we will see further on that we run into problems when we ignore the base point entirely. Luckily enough, it is straightforward to show that homotopies of maps $X \to Y$ induce equivalence relations both in the presence and absence of a base point.

\begin{cor} \label{homt}
	The relation $\simeq$ is an equivalence relation on the sets $\Mor{\Cb{Top}}{X}{Y}$ and $\Mor{\Cb{Top_*}}{X}{Y}.$
\end{cor}

\begin{proof}
	It suffices to confirm this within $\Cb{Top}$ since maps between pointed spaces are particular cases of maps between general spaces. Let $f,f',f'':X \to Y$ be continuous maps.
	\begin{itemize}
		\item Reflexivity: The map $h(x,t) = f(x)$ is a homotopy $f \simeq f.$
		\item Symmetry: Let $f \simeq f'$ and $h:f \simeq f'$ be a homotopy between them. Then $h'(x,t) = h(x,1-t)$ is a homotopy $f' \simeq f.$
		\item Transitivity: Let $f \simeq f'$ and $f' \simeq f''$ with corresponding homotopies $h:f \simeq f'$ and $h':f' \simeq f''.$ Then the concatenation function
		\[
			h^\dagger(x,t) =
			\begin{cases}
				h(x,2t) & \text{if } t \in [0,\frac{1}{2}] \\
				h'(x,2t-1)) & \text{if } t \in [\frac{1}{2},1]
			\end{cases}
		\]
		is continuous by the pasting lemma because $h(x,1) = h'(x,0)$ and thereby constitutes a homotopy $f \simeq f''.$
	\end{itemize}
\end{proof}

\begin{defin}
	Let $f:X \to Y$ be a pointed map. The equivalence class
	\[
		[f] := \qty{g:X \to Y\ |\ g \simeq f}
	\]
	is referred to as the \textit{homotopy class} of $f,$ and the collection of all such classes forms the \textit{homotopy set} $[X,Y].$ Alternatively, the homotopy set may be obtained as the quotient
	\[
		[X,Y] := \frac{\Mor{\Cb{Top}_*}{X}{Y}}{\simeq}.
	\]
	Pointed topological spaces with homotopy classes of maps between them form a category $\Cb{HoTop}_*.$ 
\end{defin}

\begin{rem}
	Categories require compositions laws to combine their morphisms. If $f,g$ are composable maps (composition given by $gf$) in $\Cb{Top}_*,$ then the composition $[g][f]$ is given by $[gf].$ As with any representable functor, $[X,-]$ and $[-,Y]$ are covariant and contravariant functors, respectively, from $\Cb{Top_*}$ to $\Cb{Set},$ and the induced maps are defined in the same fashion as with the representable functors.
\end{rem}

Recall that a homeomorphism is a continuous bijection with a continuous inverse. What we mean by inverse is that there exists another function going in the opposite direction which composes on either side to produce the relevant identity map. At the outset of this subsection, we cited a desire to relax the conditions for a homeomorphism, and homotopy seems to be a good fit. However, homotopy is a relation between maps, not spaces. By switching one equivalence relation for another, we arrive at what we sought.

\begin{defin}
	Let $f:X \to Y$ be a continuous map (not necessarily pointed). We call $f$ a \textit{homotopy equivalence} if there exists some $g:Y \to X$ such that
	\[
		fg \simeq \mathds{1}_Y \quad \quad \quad gf \simeq \mathds{1}_X,
	\]
	in which case $X$ and $Y$ are said to be \textit{homotopy equivalent}, denoted by $X \simeq Y.$ If such a $g$ exists, it is also a homotopy equivalence which is called the \textit{homotopy inverse} of $f.$
\end{defin}

To show that homotopy equivalences are closed under composition, it is sufficient to construct the homotopy inverse ($(g'g)(ff') \simeq \mathds{1}$). From this, homotopy equivalence is an equivalence relation on all pointed (and non-pointed) topological spaces, where reflexivity follows from the identity map, symmetry from the existence of a homotopy inverse, and transitivity from closure under compositions.

\begin{rem}
	When it is said that some statement/construction is valid for a topological space $X$ up to homeomorphism, what is really being said is that these statements/constructions are invariants of homeomorphism. In the same way, there are homotopy invariants, which are statements/constructions valid up to homotopy equivalence. Now we are prepared to present the titular construction of this subsection, an important homotopy invariant.
\end{rem}

\begin{defin}
	Let $(X,x_0)$ be a pointed topological space. The homotopy set $[S^1,X]$ forms a group under the associative binary operation of concatenation (c.f. transitivity of $\simeq$ in \thref{homt}), where the inverse of a loop $[f]$ is given by the loop $[f]^{-1}$ going in the opposite direction. Explicitly, a loop $f:I \to X$ has inverse
	\[
		f^{-1}(t) = f(1-t).
	\]
	This group is called the \textit{fundamental group} and is denoted by $\pi_1(X,x_0).$
	
	When $X$ is path-connected, it may be shown that fundamental groups with different base points are isomorphic, so we may instead write $\pi_1(X).$ If $\pi_1(X)$ is trivial, then $X$ is called \textit{simply connected}.
\end{defin}

The fundamental group is an algebraic construct which allows us to detect portions of the topological space which may not be deformed to a point and their ``positions'' relative to each other. However, it is often not very simply to compute, and in many cases is not abelian. Leaving that aside for the time being, the formulation of the fundamental group submits to fairly straightforward generalization.

\begin{defin}
	Let $(X,x_0)$ be a pointed topological space. Then
	\[
		\pi_n(X,x_0) = [S^n,X]
	\]
	is the \textit{$n$-th homotopy group}, where $[S^n,X]$ is the set of homotopy classes of base point-preserving continuous functions $S^n \to X.$ It is an Abelian group for $n \geq 2$ with addition defined by
	\[
		f+g \colon
		\begin{tikzcd}
			S^n \ar[r,"eq"] & S^n \vee S^n \ar[r,"f \vee g"] & X,
		\end{tikzcd}
	\]
	where $eq \colon S^n \to S^n \vee S^n$ is the quotient map which collapses the equator of $S^n$---any copy of $S^{n-1}$ containing the base point---to the wedge point of the two resulting spheres. However for $n = 0,$ there is no canonical group operation in general, so $\pi_0(X)$ is a set in bijection with the connected components of $X.$
\end{defin}

\begin{lem} \label{pizc}
	Let $\Cs{C}$ be a category and $\sim$ be an equivalence relation on $\ope{Obj}(\Cs{C})$ generated by all pairs of objects $(A,B)$ such that there is a morphism $A \to B$. Then $$\pi_0(B\Cs{C}) = \ope{Obj}(\Cs{C})/\sim$$
\end{lem}

This lemma tells us that two objects $A$ and $B$ are equivalent mod $\sim$ if and only if there is a \emph{mixed path} between them. There are four possible forms a mixed path can take:
\begin{align*}
	&
	\begin{tikzcd}[ampersand replacement = \&]
		\bullet \ar[r] \& \bullet \& \ar[l] \cdots \& \ar[l] \bullet \ar[r] \& \bullet
	\end{tikzcd}
	&
	&
	\begin{tikzcd}[ampersand replacement = \&]
		\bullet \ar[r] \& \bullet \& \ar[l] \cdots \ar[r] \& \bullet \& \ar[l] \bullet
	\end{tikzcd}
	\\
	&
	\begin{tikzcd}[ampersand replacement = \&]
		\bullet \& \ar[l] \bullet \ar[r] \& \cdots \& \ar[l] \bullet \ar[r] \& \bullet
	\end{tikzcd}
	&
	&
	\begin{tikzcd}[ampersand replacement = \&]
		\bullet \& \ar[l] \bullet \ar[r] \& \cdots \ar[r] \& \bullet \& \ar[l] \bullet
	\end{tikzcd}
\end{align*}




\newpage

\section{Torsion Products as Homotopy Groups}



\subsection{Category $\Tor{}{R}{A}{B}$}
If $_RP$ is a left $R$-module, then $P_R^* = \Hom{R}{P}{R}$ is a right $R$-module. If $f:P_1\to P_2$ is a homomorphism of left $R$-modules, then $f^*:P_2^* \to P_1^*$ is a homomorphism of right $R$-modules.

\begin{defin}
	Let $R$ be a ring, $A_R$ and $_RB$ be right and left $R$-modules, respectively. Define the \textit{torsion category} $\mathrm{Tor}^R(A,B)$ of modules $A$ and $B$ as a category whose
	\begin{itemize}
		\item objects are triples $(P,\varepsilon,\eta)$, where $P$ is a finitely generated projective right $R$-module with homomorphisms $\varepsilon:P \to A$ and $\eta:P^* \to B$,
		\item morphisms between triples $(P_1,\varepsilon_1,\eta_1)$ and $(P_2,\varepsilon_2,\eta_2)$ are given by the homomorphisms $f:P_1 \to P_2$ such that the following diagrams commute:
		\[
			\begin{tikzcd}
				P_1 \ar[rr,"f"] \ar[rd,"\varepsilon_1"'] & & \ar[ld,"\varepsilon_2"] P_2 \\
				& A &
			\end{tikzcd}
			~~
			\begin{tikzcd}
				P_1^* \ar[rd,"\eta_1"'] & & \ar[ld,"\eta_2"] \ar[ll,"f^*"'] P_2^* \\
				& B, &
			\end{tikzcd}
		\]
		\item composition of morphisms of triples is given by the composition of morphisms between the finitely generated projective modules within the triple,
		\item identity on the triple $(P,\varepsilon,\eta)$ is the identity homomorphism on $P$.
	\end{itemize}
\end{defin}


\begin{defin}
	Define a bifunctor $\Box:\Tor{}{R}{A}{B} \times \Tor{}{R}{A}{B} \to \Tor{}{R}{A}{B}$ as $$(P_1,\varepsilon_1,\eta_1) \boxe (P_2,\varepsilon_2,\eta_2) = (P_1 \oplus P_2,(\varepsilon_1,\varepsilon_2),(\eta_1,\eta_2)^T)$$ for objects. Let $p:(P_1,\varepsilon_1,\eta_1) \to (P_2,\varepsilon_2,\eta_2)$ and $q:(Q_1,\varepsilon_1',\eta_1') \to (Q_2,\varepsilon_2',\eta_2')$ be morphisms in $\Tor{}{R}{A}{B}$. Then $p \boxe q:P_1 \oplus Q_1 \to P_2 \oplus Q_2$ makes the following diagrams commute:
	\[
		\begin{tikzcd}
			P_1 \oplus Q_1 \ar[rr,"p \boxe q"] \ar[rd,"{(\varepsilon_1,\varepsilon_1')}"'] & & \ar[ld,"{(\varepsilon_2,\varepsilon_2')}"] P_2 \oplus Q_2 \\
			& A &
		\end{tikzcd}
		~~
		\begin{tikzcd}
			P_1^* \oplus Q_1^* \ar[rd,"{(\eta_1,\eta_1')}"'] & & \ar[ld,"{(\eta_2,\eta_2')}"] \ar[ll,"p^* \boxe q^*"'] P_2^* \oplus Q_2^* \\
			& B, &
		\end{tikzcd}
	\]
	where
	\[
		p \boxe q = 
		\begin{pmatrix}
			p & 0 \\
			0 & q
		\end{pmatrix}.
	\]
	The object $(0,0,0)$ is an identity for this operation.
\end{defin}


\subsection{Homotopic Relationship Between Torsion Products and Category $\Tor{}{R}{A}{B}$}
The primary result of Robinson's paper \cite{Robinson1981} is as follows.
\begin{theo}
	Let $R$ be a ring, $A_R$ a right $R$-module, and $_RB$ a left $R$-module. Then
	\[
		\pi_n(B\mathrm{Tor}^R(A,B))\cong \Tor{n}{R}{A}{B}
	\]
	as abelian groups.
\end{theo}

This may be used to carry out a computation of tensor products since $\Tor{0}{n}{A}{B} \cong A \otimes_R B.$ We will use the theorem to show that $$\Z/4\Z \otimes_\Z\Z/6\Z \cong \Z/2\Z.$$

\begin{proof}
	The outline of this computation is as follows. Note that projective $\Z$-modules are precisely free abelian groups.
	\begin{enumerate}
		\item Determine generators: Both $P_1$ and $P_2$ may be considered as free $\Z$ modules and so may be written as $\Z^n$ and $\Z^m$, respectively. The maps may then be written as $$\varepsilon_1 = (x_1,x_2,\ldots,x_n) ~~ \text{ and } ~~ \eta_1 = (y_1,y_2,\ldots,y_n)^T$$ for $P_1$. Because $\Hom{\Z}{\Z^n}{\Z} \cong \Z^n$, the maps for $P_1^*$ will be similar (but transposed from $\varepsilon_1$ and $\eta_1$). With this, we may write $$(\Z^n,(x_1,x_2,\ldots,x_n),(y_1,y_2,\ldots,y_n)^T) = \Box_{i,j}x_i y_j(\Z,1,1),$$ so the generators we obtain are of the form $(\Z,1,1)$. The same argument applies to $\Z^m$.
		\item Establish path-connectedness within category: Because all objects are generated by $(\Z,1,1)$ through multiplication, we have the following triangles to satisfy in order to establish path-connectedness to the zero module object $(0,0,0)$ for any object. We will omit the subscripts on the maps $\varepsilon$ and $\eta$ as a result of working with only one generator.
		\begin{align*}
			&
			\begin{tikzcd}[ampersand replacement = \&]
				0 \ar[rd] \ar[rr] \& \& \Z \ar[ld,"\varepsilon"] \\
				\& \Z/4\Z \& 
			\end{tikzcd}
			&
			&
			\begin{tikzcd}[ampersand replacement = \&]
				\Z \ar[rd,"\varepsilon"'] \ar[rr] \& \& 0 \ar[ld] \\
				\& \Z/4\Z \& 
			\end{tikzcd}
			\\
			&
			\begin{tikzcd}[ampersand replacement = \&]
				0 \ar[rd]  \& \& \ar[ll] \Z \ar[ld,"\eta"] \\
				\& \Z/6\Z \& 
			\end{tikzcd}
			&
			&
			\begin{tikzcd}[ampersand replacement = \&]
				\Z \ar[rd,"\eta"']  \& \& \ar[ll] 0 \ar[ld] \\
				\& \Z/6\Z \& 
			\end{tikzcd}
	\end{align*}
	Then the conditions for a morphism existing between $(\Z,\varepsilon,\eta)$ and $(0,0,0)$ are only satisfied when $\varepsilon = 0$ or $\eta = 0$. Therefore, only objects of the form $(\Z^x,0,\eta)$ or $(\Z^x,\varepsilon,0)$ are connected by a path to $(0,0,0)$.
	\item Determine order of generators: We now multiply the generator by specific integers in order to find objects path-connected to $(0,0,0)$: $$4(\Z,1,1) = (\Z,0,1) \sim (0,0,0)\text{ and }6(\Z,1,1) = (\Z,1,0) \sim (0,0,0).$$ Because $2 = \ope{gcd}(4,6)$, there exist integers $p,q \in \Z$ such that $4p+6q = 2$ (namely $q = 1$ and $p = -1$), which implies $$(6-4)(\Z,1,1) = 2(\Z,2,2) \sim (0,0,0).$$ Because $(\Z,1,1)$ is not connected to $(0,0,0)$ by a morphism, we see that $B\Tor{}{\Z}{\Z/4\Z}{\Z/6\Z}$ consists of two path-components. By \ref{pizc}, $$\pi_0(B\Tor{}{\Z}{\Z/4\Z}{\Z/6\Z) \cong (\Z,1,1}/2(\Z,1,1) \cong \Z/2\Z.$$ This aligns with what is obtained by direct computations of the tensor product.
	\end{enumerate}
\end{proof}


\section*{Acknowledgements}
I would like to thank my mentor Alex Sorokin for his time and attention in both preparing me for and working with me on this project. I learned a great deal from him about math and beyond, and his constant patience with me was very much appreciated. Under his guidance, I made significant progress towards being able to identify and postulate about important and ``good'' questions. I would not be upset to be half the math educator he is when I pursue graduate studies.

In addition, I would like to thank Dr. Robin Walters for his advice and feedback on the progress of the project and for providing further opportunities for me to integrate what I learned through meaningful discussions during my progress reports. And I would remiss to not mention Professor Anthony (Tony) Iarrobino for his comments on my presentations, writing, and aspirations which seemed to perfectly balance generality with detail. And overall, I would like to thank the Northeastern RTG group for supporting this very informative and stimulating experience.

\end{document}